\documentclass[preprint,12pt]{elsarticle}

\usepackage{amssymb}
\usepackage{arydshln}
\usepackage{mathrsfs}
\usepackage{amsmath}
\usepackage{amsfonts}
\usepackage[all,cmtip]{xy}
\usepackage[colorlinks]{hyperref}
\usepackage{color}
\usepackage{cleveref}
\usepackage{ulem}
\usepackage[abs]{overpic}
\usepackage{float}
\usepackage{enumitem}
\usepackage{tikz}
\usepackage{lmodern}
\usepackage{anyfontsize}

\usepackage{etoolbox}
\patchcmd{\thebibliography}
  {\settowidth}
  {\setlength{\itemsep}{0pt plus 0.1pt}\settowidth}
  {}{}
\apptocmd{\thebibliography}
  {\small}
  {}{}

\usepackage{tikz-cd}
\usepackage[toc,page]{appendix}

\usepackage{graphicx}
\usepackage{subcaption}
\usepackage[utf8]{inputenc}
\usepackage[english]{babel}
\usepackage{dsfont}
\usepackage{xcolor}
\usepackage{arydshln}
\usepackage{faktor}
\usepackage{mathtools}
\usepackage{scalerel}
\usepackage{booktabs}

\def\endpf{\hbox{\vrule height1.5ex width.5em}}

\usepackage{dsfont}

\makeatletter
\newcommand*{\da@rightarrow}{\mathchar"0\hexnumber@\symAMSa 4B }
\newcommand*{\da@leftarrow}{\mathchar"0\hexnumber@\symAMSa 4C }
\newcommand*{\xdashrightarrow}[2][]{%
	\mathrel{%
		\mathpalette{\da@xarrow{#1}{#2}{}\da@rightarrow{\,}{}}{}%
	}%
}
\newcommand{\xdashleftarrow}[2][]{%
	\mathrel{%
		\mathpalette{\da@xarrow{#1}{#2}\da@leftarrow{}{}{\,}}{}%
	}%
}
\newcommand{\xdashdownarrow}[2][]{%
	\mathrel{%
		\mathpalette{\da@xarrow{#1}{#2}\da@downarrow{}{}{\,}}{}%
	}%
}
\newcommand*{\da@xarrow}[7]{%
	\sbox0{$\ifx#7\scriptstyle\scriptscriptstyle\else\scriptstyle\fi#5#1#6\m@th$}%
	\sbox2{$\ifx#7\scriptstyle\scriptscriptstyle\else\scriptstyle\fi#5#2#6\m@th$}%
	\sbox4{$#7\dabar@\m@th$}%
	\dimen@=\wd0 %
	\ifdim\wd2 >\dimen@
	\dimen@=\wd2 %   
	\fi
	\count@=2 %
	\def\da@bars{\dabar@\dabar@}%
	\@whiledim\count@\wd4<\dimen@\do{%
		\advance\count@\@ne
		\expandafter\def\expandafter\da@bars\expandafter{%
			\da@bars
			\dabar@ 
		}%
	}%  
	\mathrel{#3}%
	\mathrel{%   
		\mathop{\da@bars}\limits
		\ifx\\#1\\%
		\else
		_{\copy0}%
		\fi
		\ifx\\#2\\%
		\else
		^{\copy2}%
		\fi
	}%   
	\mathrel{#4}%
}
\makeatother

\usepackage{graphicx}

\newcommand\tikzmark[1]{%
  \tikz[remember picture,overlay]\coordinate (#1);}

\newcommand{\underbracedmatrixl}[2]{%
  \left(\;
  \smash[b]{\underbrace{
    \begin{matrix}#1\end{matrix}
  }_{#2}}
  \;\right.
  \vphantom{\underbrace{\begin{matrix}#1\end{matrix}}_{#2}}
}

\newcommand{\underbracedmatrixll}[2]{%
  \left(\;\hspace{-.27in}
  \smash[b]{\underbrace{
    \begin{matrix}#1\end{matrix}
  }_{#2}}
  \;\right.
  \vphantom{\underbrace{\begin{matrix}#1\end{matrix}}_{#2}}
}

\newcommand{\underbracedmatrixr}[2]{%
  \left. \;
  \smash[b]{\underbrace{
    \begin{matrix}#1\end{matrix}
  }_{#2}}
  \;\right)
  \vphantom{\underbrace{\begin{matrix}#1\end{matrix}}_{#2}}
}

\newcommand{\underbracedmatrixrr}[2]{%
  \left. \;
  \smash[b]{\underbrace{
    \begin{matrix}#1\end{matrix}
  }_{#2}}
  \;\hspace{-.32in}\right)
  \vphantom{\underbrace{\begin{matrix}#1\end{matrix}}_{#2}}
}

\newcommand{\bigslant}[2]{{\raisebox{.2em}{$#1$}\left/\raisebox{-.2em}{$#2$}\right.}}

\newtheorem{theorem}{Theorem}[section]
\newtheorem{lemma}[theorem]{Lemma}
\newtheorem{corollary}[theorem]{Corollary}
\newtheorem{proposition}[theorem]{Proposition}

\newtheorem{definition}[theorem]{Definition}
\newtheorem{example}[theorem]{Example}
\newtheorem{remark}[theorem]{Remark}

\newtheorem{definitionlemma}[theorem]{Definition-Lemma}

\usepackage{url}

\begin{document}

\begin{frontmatter}

\title{The automorphism groups of generalized Kausz compactifications and spaces of complete collineations}

\author{Hanlong Fang} 

\affiliation{organization={Peking University},           city={Beijing},  postcode={100871},     state={Beijing},     country={China}}

\ead{hlfang@pku.edu.cn}

\begin{abstract}
In this paper, we determine the automorphism groups of generalized Kausz compactifications $\mathcal T_{s,p,n}$. By establishing the (semi-)positivity of the anticanonical bundles of $\mathcal T_{s,p,n}$, we also determine the automorphism groups of generalized spaces of complete collineations $\mathcal M_{s,p,n}$. The results in this paper are partially taken from the author's earlier arxiv post (Canonical blow-ups of grassmann manifolds, arxiv:2007.06200).
% We would like to encourage you to list your keywords within
% the abstract section using the There are K\"ahler-Einstein metrics on $\mathcal M_{s,p,n}$ if $|s+p-n|\gg r$, $|s-p|\gg r$, and there are K\"ahler-Einstein metrics on $\mathcal M_{r,r,2r}$.\keywords{...} command.

\end{abstract}


\end{frontmatter}

% typeset the title of the contribution

%

\tableofcontents

\section{Introduction}
Spaces of complete collineations are classical objects in algebraic geometry compactifying spaces of linear maps of maximal rank. The history of study dates back at least to Study \cite {Stu}, Severi \cite{Sev2}, Van der Waerden \cite{Van}, etc.~on complete conics from the perspective of enumerative geometry. Vainsencher \cite{Va} described the spaces of complete collineations as the closures of rational maps %from projective spaces into products of projective spaces 
via the Pl\"ucker coordinates of Grassmannians, and showed that they are wonderful in the sense of De Concini--Procesi \cite{DP}.  %Given a reductive group $G$,  a $G$-variety is called wonderful if it is smooth, projective, the $G$-stable prime divisors $D_1,\cdots,D_r$ are smooth with simple normal crossings, and  for any $I\subset\{1,\cdots,r\}$ the strata $\cap_{i\in I}D_i$ one-to-one correspond to the closures of the $G$-orbits. 
Wonderful varieties play an important role in enumerative geometry  \cite{DP2,LLT,MMM} and the theory of spherical varieties \cite{Kn,L,BP}, and wonderful compactifications of adjoint semisimple groups are widely used in arithmetic geometry and representation theory \cite{Fa,Sp1,Lu2}. However, despite the importance, wonderful compactifications cease to exist for reductive groups beyond the adjoint type. In 2000, Kausz \cite{Ka} constructed modular compactifications of general linear groups in generalizing the wonderful compactifications. % On the other hand, the natural torus action on the Grassmannian suggests a deep connection between compactifications and the total family associated to the Hilbert/Chow quotient in the sense of \cite{Kap,BS}.

% except the existence of multiple closed orbits. 

The beautiful geometry of the Kausz compactification and its generalization have been investigated in \cite{FW}, where a simple uniform picture has been developed to incorporate the Kausz compactifications and the spaces of complete collineations. % Note that such incorporation is influenced by the study of the Euler vector fields on Grassmannians \cite{Siu,Hw}. Indeed, desingularizations of such foliations encode the centers of general linear groups which are elusive in the framework of wonderful compactifications; moreover, they are total families associated to the Hilbert/Chow quotients of the natural $1$-dimensional torus actions on the Grassmannians.
%Through this alternative construction, we show that a Kausz compactification is a flat family parametrizing orbital rational curves in the sense of Fujiki \cite{Fu,BS}, and the maximal-dimensional connected component of the fixed point scheme is one of the spaces of complete collineations, which in the meantime serves as the base of the family. In particular, together with the corresponding spaces of complete collineations, the Kausz compactifications are toroidal embeddings in the sense of \cite{BM}. As an application, by comparing Kausz's original construction and our new one, we resolve the Landsberg-Manivel birational maps \cite{LM} by explicit blow-ups and blow-downs. Since after the desingularization the source and the sink associated to the torus action are isomorphic, we resolve certain birational maps among projective bundles over Grassmannians. % by the isomorphism among the source and the sink. Recently, our idea of blowing up Grassmannians along the closures of the stable/unstable manifolds according to the   Bia{\l}ynicki-Birula  decomposition has been applied by Ding \cite{Di} to give a resolution of the Landsberg-Manivel birational maps.
Similar ideas have been further developed for other moduli spaces, such as Faltings--Lafforgue spaces \cite{FZ23} and Koll\'ar--Shepherd-Barron--Alexeev  moduli of marked cubic surfaces \cite{FSW}.

%Let $G(p,n)$, $0<p<n$, be the Grassmannian  consisting of $p$-planes in the  $n$-space. The automorphism groups of various moduli spaces have been attracting mathematicians' attention. Bruno and Mella \cite{BM1} determined the automorphism group of the Knudsen-Mumford compactification $\overline M_{0,n}$, while Bolognesi and Massarenti \cite{BM2} computed it for the GIT quotient of configurations of ordered points on $\mathbb P^1$ with respect to the symmetric polarization. For the Deligne-Mumford compactification $\overline M_{g,n}$, Massarenti \cite{M1} first derived the automorphism group over the complex numbers; this result was later extended to positive characteristic by Fantechi and Massarenti \cite{FM1}. Massarenti and Mella \cite{MM1} determined the automorphisms of Hassett’s moduli spaces in characteristic zero, a result generalized to arbitrary fields by Fantechi and Massarenti \cite{FM2}. Furthermore, Massarenti \cite{M2} computed the automorphism group of the Fulton-MacPherson compactification, and Araujo, Fassarella, Kaur, and Massarenti \cite{AFKM} did so for moduli spaces of parabolic vector bundles. In this paper, we will determine the automorphism of the (generalized) Kausz modular compactification and the (generalized) 

The problem of determining automorphism groups of moduli spaces has been a focus of significant research. Key results include the computation for the Knudsen--Mumford compactification $\overline{M}_{0,n}$ by Bruno and Mella \cite{BM1}, and for a GIT quotient of point configurations on $\mathbb{P}^1$ for the standard linearization $\mathcal O(1,\cdots,1)$ by Bolognesi and Massarenti \cite{BM2}.
Subsequent work has addressed other major compactifications: Massarenti established  \cite{M1} the corresponding result for the Deligne--Mumford compactification $\overline{M}_{g,n}$ over $\mathbb{C}$, with Fantechi and Massarenti \cite{FM1} providing the generalization to positive characteristic. The automorphisms of Hassett’s spaces were determined in characteristic zero by Massarenti and Mella \cite{MM1} and in arbitrary fields by Fantechi and Massarenti \cite{FM2}. The automorphism groups of the Fulton--MacPherson compactification and of moduli spaces of parabolic vector bundles were computed by Massarenti \cite{M2}, and by Araujo et al. \cite{AFKM}, respectively. In a series of works, Masserenti \cite{M3,M4} computed 
the automorphisms of spaces of complete forms.

In this paper, we will determine the automorphism of the generalized Kausz modular compactification and the generalized space of complete collineations.

\medskip

Now let us be more precise. Let $E$ be a free $\mathbb Z$-module of rank $n$. Given an integer $p$ such that
$0<p<n$, consider the Grassmannian $Gr(p,E)$ embedded in the projective space $\mathbb P(\bigwedge^pE)$ via the Pl\"ucker embedding. The general linear group ${\rm GL}(E)$ acts linearly on $\bigwedge^pE$ and the induced action on $\mathbb P(\bigwedge^pE)$ stabilizes $Gr(p,E)$.
Choose a non-trivial decomposition $E=E_1\oplus E_2$ and set $s=\dim(E_1)$. Then
%$0<s<n$ and 
we have a decomposition
\begin{equation}\label{cd}
\bigwedge^{p}E:=\bigoplus\limits_{k=0}^p\bigwedge^{k}E_1\otimes\bigwedge^{p-k}E_2.
\end{equation}
The subgroup $G:={\rm GL}(E_1)\times {\rm GL}(E_2)\subset {\rm GL}(E)$ stabilizes
(\ref{cd}) where summands are irreducible $G$-modules. We derive a rational $G$-equivariant map
\begin{equation}\label{kspe}
\mathcal K(s,p,E):Gr(p,E)\xdashrightarrow{\,\,\,\,\,\,\,\,\,\,\,\,\,\,\,}\mathbb P(\bigwedge^{p}E)\times\prod\limits_{k=0}^p\mathbb P(\bigwedge^{k}E_1\otimes\bigwedge^{n-k}E_2).  
\end{equation}
Here we make the convention that when $\bigwedge^{k}E_1\otimes\bigwedge^{n-k}E_2=0$, $\mathbb P(\bigwedge^{k}E_1\otimes\bigwedge^{n-k}E_2)$ is a point. Denote by $\mathcal T(s,p,E)$ the closure of the image of the rational map $\mathcal K(s,p,E)$. Denote by $\mathcal M(s,p,E)$ the image of $\mathcal T(s,p,E)$ under the natural projection
\begin{equation*}%\label{nproj}
\mathbb P(\bigwedge^{p}E)\times\prod\limits_{k=0}^p\mathbb P(\bigwedge^{k}E_1\otimes\bigwedge^{n-k}E_2)\longrightarrow \prod\limits_{k=0}^p\mathbb P(\bigwedge^{k}E_1\otimes\bigwedge^{n-k}E_2).
\end{equation*}

\begin{definition}
We call $\mathcal T(s,p,E)$ {\it the generalized Kausz compactifications} and $\mathcal M(s,p,E)$ {\it the generalized spaces of complete collineations}. When there is no ambiguity, we write $Gr(p,E)$, $\mathcal T(s,p,E)$, $\mathcal K(s,p,E)$, $\mathcal M(s,p,E)$ as $G(p,n)$, $\mathcal T_{s,p,n}$, $\mathcal K_{s,p,n}$, $\mathcal M_{s,p,n}$, respectively.    
\end{definition}

%\begin{example}\label{}$\mathcal T_{s,1,n}$ is the blow-up of $\mathbb{P}^{n-1}$ along the union of the two disjoint linear subspaces $\mathbb{P}^{s-1}$ and $\mathbb{P}^{n-s-1}$. When $n=2s=2p$, $\mathcal T_{p,p,2p}$ is isomorphic to the Kausz compactification of the general linear group ${\rm GL}_p$ (see \cite{FW}). %; $\mathcal K_{s,1,n}$ is given by \begin{equation}\mathcal K_{s,1,n}([x_1,\cdots, x_s,y_1,\cdots,y_{n-s}])=\left([x_1,\cdots, x_s,y_1,\cdots,y_{n-s}], [x_1,\cdots, x_s],[y_1,\cdots,y_{n-s}]\right).\end{equation}\end{example}

According to \cite{FW}, $\mathcal T_{p,p,2p}$ is isomorphic to the Kausz compactification of the general linear group ${\rm GL}_p$, and $\mathcal M_{p,p,n}$ is isomorphic to the space of complete collineations associated with linear maps from a $p$-dimensional vector space to an $(n-p)$-dimensional one.
In general, $\mathcal T_{s,p,n}$ are derived from the Kausz compactifications by parabolic induction as follows. Define a rational map
$\pi_1:Gr(p,E)\dashrightarrow Gr(p,E_1)$ by sending $W\in Gr(p,E)$ with $W\cap E_2=0$ to 
$\pi_1(W)\in Gr(p,E_1)$ such that $W+E_2=\pi_1(W)\oplus E_2$. %, which  
$\pi_1$ extends to an equivariant morphism $f_1:\mathcal T(s,p,E)\rightarrow Gr(p,E_1)$. Viewing $Gr(p,E_1)$ as the quotient of ${\rm GL}(E_1)$ by the stabilizer of a subspace $W\subset E_1$, we 
see that $f_1$ is a locally fibration with fiber at $W$ being $\mathcal T(p,p,W\oplus E_2)$. Then we have 
\begin{proposition}[Proposition 1.5 in \cite{FW}]\label{tfiber}
\begin{enumerate}[label=\rm(\Alph*)]
\item If $p<n-s$ and $p<s$, $\mathcal T_{s,p,n}$ ({resp.}~$\mathcal M_{s,p,n}$) is a locally trivial fibration over $G(p,s)$ with the fiber $\mathcal T_{p,p,n-s+p}$ ({resp.}~$\mathcal M_{p,p,n-s+p}$), and a locally trivial fibration over $G(p,n-s)$ with the fiber  $\mathcal T_{p,p,s+p}$ ({resp.}~$\mathcal M_{p,p,s+p}$).
    
\item If $n-s<p<s$, $\mathcal T_{s,p,n}$ ({resp.}~$\mathcal M_{s,p,n}$) is a locally trivial fibration over  $G(p,s)$ with the fiber $\mathcal T_{n-s,n-s,n-s+p}$ ({resp.}~$\mathcal M_{n-s,n-s,n-s+p}$), and a locally trivial fibration over $G(s+p-n,s)$ with the fiber $\mathcal T_{n-s,n-s,2n-s-p}$ ({resp.}~$\mathcal M_{n-s,n-s,2n-s-p}$).
    
\item If $p=n-s<s$, $\mathcal T_{s,p,n}$ ({resp.}~$\mathcal M_{s,p,n}$) is a locally trivial fibration over $G(p,s)$ with the fiber $\mathcal T_{p,p,2p}$ ({resp.}~$\mathcal M_{p,p,2p}$).
\end{enumerate}
\end{proposition}

The following theorem reveals the beautiful geometry of $\mathcal T_{s,p,n}$ and $\mathcal M_{s,p,n}$.
\begin{theorem}[Theorem 1.3 in \cite{FW}]\label{gwond}
The inverse of (\ref{kspe}) extends to a regular morphism \begin{equation}\label{rspn}
R_{s,p,n}:\mathcal T_{s,p,n}\longrightarrow G(p,n).    
\end{equation} 
$\mathcal T_{s,p,n}$ is smooth and projective over ${\rm Spec}\,\mathbb Z$ with an action of $G$. Set
\begin{equation}\label{rank}
r=\min\{s,n-s,p,n-p\}.
\end{equation} 
Then the complement of the open
$G$-orbit in $\mathcal T_{s,p,n}$ consists of $2r$ smooth prime divisors with simple normal crossings $D^+_1, 
\cdots,D^+_r,D^-_1$, $\cdots,D^-_r$ such that the following holds.
\begin{enumerate}[label={\rm(\Alph*)}]
\item $D^+_1\cong D^-_1\cong\mathcal M_{s,p,n}$ are smooth and projective over ${\rm Spec}\,\mathbb Z$. 

\item There is a $G$-equivariant flat retraction 
\begin{equation}\label{B}
\mathcal P_{s,p,n}:\mathcal T_{s,p,n}\longrightarrow D_1^-\cong\mathcal M_{s,p,n} 
\end{equation} such that the restriction  $\mathcal P_{s,p,n}|_{D^+_1}:D_1^+\rightarrow D_1^-$ %(resp. $\mathcal P_{s,p,n}|_{D^+_1}$) 
is an isomorphism and that for $2\leq i\leq r$,
\begin{equation*}%\label{KGE}
\mathcal P_{s,p,n}(D^-_i)=\mathcal P_{s,p,n}(D^+_{r+2-i})=:\check D_i\,\,    
\end{equation*}
%\item $D^-_1\cong\mathcal M_{s,p,n}$ is a wonderful variety with the $G$-stable divisors $\check D_2,\cdots,\check D_r$, where $\check D_i:=\mathcal P_{s,p,n}(D^-_i)$,  $2\leq i\leq r$.
%\item When $n=2s$ or $n=2p$ there is a holormorphic automoprhism  $\sigma$ of $\mathcal T_{s,p,n}$ such that $\sigma\circ \sigma$ is the identity map  and $\sigma(D^{\pm}_i)=D^{\mp}_i$ for $1\leq i\leq r$.

\item The closures of $G$-orbits in $\mathcal T_{s,p,n}$ are one-to-one given by
\begin{equation*}%\label{inrule}
\left(\bigcap\nolimits_{\,\,i\in I^+}D^+_i\right)\bigcap\left(\bigcap\nolimits_{\,\,i\in I^-}D^-_i\right)%\mathbin{\scaleobj{2.1}{\backslash}} (\bigcup\nolimits_{j\notin I^-}D^-_j\mathbin{\scaleobj{1}{\bigcup}}\,\bigcup\nolimits_{j\notin I^+}D^+_j),
\end{equation*} 
for subsets $I^+,I^-\subset\{1,2,\cdots,r\}$
such that $\min(I^+)+\min(I^-)\geq r+2$ with the convention that $\min(\emptyset)=+\infty$.

\item  $D_1^-$ is wonderful with the $G$-stable divisors $\check D_i$,  $2\leq i\leq r$.
\end{enumerate}
\end{theorem}

In the remainder of the paper, we fix an algebraically closed field $\mathbb K$ of characteristic zero. 

By interpreting $R_{s,p,n}$ as a blow-up, we first derive that
\begin{corollary}\label{picb}
The Picard group of $\mathcal T_{s,p,n}$ is freely generated by
\begin{equation}\label{basis1}
\left\{ \begin{array}{ll}
(R_{s,p,n})^*\mathcal O_{ G(p,n)}(1),D^+_1,\cdots,D^+_{r},D^-_1,\cdots,D^-_{r} &{\rm when}\,\,p<s,\,n-s\neq p\\
(R_{s,p,n})^*\mathcal O_{ G(p,n)}(1),D^+_1,\cdots,D^+_{r},D^-_1,\cdots,D^-_{r-1} &{\rm when}\,\,n-s=p<s\\
(R_{s,p,n})^*\mathcal O_{ G(p,n)}(1),D^+_1,\cdots,D^+_{r-1},D^-_1,\cdots,D^-_{r-1}&{\rm when}\,\,n-s=p=s\,\,\\
    \end{array}\right.,
\end{equation}
where $\mathcal O_{G(p,n)}(1)$ is the hyperplane line bundle on $G(p,n)$.
\end{corollary}
\begin{corollary}\label{mpicb}
The Picard group of $\mathcal M_{s,p,n}$ is freely generated by
\begin{equation*}
\left\{ \begin{array}{ll}
(R_{s,p,n})^*\mathcal O_{G(p,n)}(1)|_{\mathcal M_{s,p,n}},\check D_1,\check D_{2},\cdots,\check D_{r} &{\rm when}\,\,p<s,\,n-s\neq p\\
(R_{s,p,n})^*\mathcal O_{ G(p,n)}(1)|_{\mathcal M_{s,p,n}},\check D_1,\check D_{2},\cdots,\check D_{r-1}&{\rm when}\,\,n-s=p<s\\
\check D_1,\check D_{2},\cdots,\check D_{r-1} &{\rm when}\,\,n-s=p=s\,\,\\
    \end{array}\right..
\end{equation*}
\end{corollary}

Denote by ${\rm Z}(n,\mathbb K)$ the center of ${\rm GL}(n,\mathbb K)$. 
Define a subgroup of the general linear group ${\rm GL}(n,\mathbb K)$ by
\begin{equation*}%\label{sns}
{\rm GL}(s,\mathbb K)\times {\rm GL}(n-s,\mathbb K):=\left\{\left.\left(
\begin{matrix}
g_1&0\\
0&g_2\\
\end{matrix}\right)\right\vert_{} g_1\in {\rm GL}(s,\mathbb K), g_2\in {\rm GL}(n-s,\mathbb K)\right\}.
\end{equation*}
It is well-known that the automorphism group of the Grassmannian $G(p,n)$ is  ${\rm PGL}(n,\mathbb K)$  when $n\neq 2p$, and  ${\rm PGL}(n,\mathbb K)\rtimes\mathbb Z/2\mathbb Z$ when $n=2p$ \cite{C}. %(the discrete symmetry is given by the so-called dual map).  

We determine the symmetry of $\mathcal T_{s,p,n}$ as follows.
\begin{theorem}\label{auto1}
The automorphism group ${\rm Aut}(\mathcal T_{s,p,n})$ of $\mathcal T_{s,p,n}$ is
\begin{equation*}
   \bigslant{{\rm GL}(s,\mathbb K)\times {\rm GL}(n-s,\mathbb K)} {{\rm Z}(n,\mathbb K)}\,\,,
\end{equation*}
except for the following cases.
\begin{enumerate}[label={\rm(\Alph*)}]
	\item {\rm (USD\footnote{USD stands for a discrete symmetry which turns $\mathcal T_{s,p,n}$ upside down.} case)}. If $n=2s$ and $s\neq p$,
	\begin{equation*}
     {\rm Aut}(\mathcal T_{s,p,2s})=\left(\bigslant{{\rm GL}(s,\mathbb K)\times {\rm GL}(n-s,\mathbb K)} {{\rm Z}(n,\mathbb K)}\right)\rtimes\mathbb Z/2\mathbb Z\,\,.
    \end{equation*} 
	\item {\rm (DUAL\footnote{DUAL stands for the discrete symmetry  induced by the dual map of $G(p,2p)$.} case)}. If $n=2p$ and $s\neq p$,
	\begin{equation*}
      {\rm Aut}(\mathcal T_{s,p,2p})=\left(\bigslant{{\rm GL}(s,\mathbb K)\times {\rm GL}(n-s,\mathbb K)} {{\rm Z}(n,\mathbb K)}\right)\rtimes\mathbb Z/2\mathbb Z\,\,.
    \end{equation*} 
	\item {\rm (USD+DUAL\footnote{The two types of discrete symmetries in (A) and (B) coexist in (C).} case)}.  If $n=2s=2p$ and $p\geq 2$,
	\begin{equation*}
      {\rm Aut}(\mathcal T_{p,p,2p})=\left(\bigslant{{\rm GL}(s,\mathbb K)\times {\rm GL}(n-s,\mathbb K)} {{\rm Z}(n,\mathbb K)}\right)\rtimes\mathbb Z/2\mathbb Z\rtimes\mathbb Z/2\mathbb Z\,\,.
    \end{equation*} 
    \item {\rm (Degenerate cases)}. 
    \begin{enumerate}[label={\rm(\alph*).}]
    \item  $\mathcal T_{1,1,2}\cong\mathbb {P}^1$ and  ${\rm Aut}(\mathcal T_{1,1,2})={\rm PGL}(2,\mathbb K)$. 
    \item For $m\geq 2$, $\mathcal T_{m,1,m+1}\cong\mathcal T_{m,m,m+1}\cong\mathcal T_{1,1,m+1}\cong\mathcal T_{1,m,m+1}$. Their  automorphism groups are isomorphic to the following subgroup of ${\rm PGL}(n,\mathbb K)$. \begin{equation}\label{parabolic} \left\{\left. \left( \begin{matrix} V&\eta\\ 0&1\\ \end{matrix}\right)\right\vert_{} \begin{matrix} V\in {\rm GL}(n-1,\mathbb K) \,,\,\,\eta \,\,{\rm is \,\,a\,\,} (n-1)\times 1\,\,{\rm matrix}\\ \end{matrix} \right\}\,. \end{equation}
    \end{enumerate}
\end{enumerate}   
\end{theorem}

To determine the symmetry of $\mathcal M_{s,p,n}$, we first establish the following positivity results, which generalizes that for $\mathcal M_{p,p,2p}$ given by De Concini and Procesi \cite{DP}.
\begin{proposition}\label{mfano}
The anti-canonical bundle $-K_{\mathcal M_{s,p,n}}$  of $\mathcal M_{s,p,n}$ is ample.
\end{proposition}
\begin{proposition}\label{fano}
The anti-canonical bundle $-K_{\mathcal T_{s,p,n}}$ of $\mathcal T_{s,p,n}$ is big and numerical effective. $-K_{\mathcal T_{s,p,n}}$ is ample if and only if the rank $r\leq 2$ (see (\ref{rank})). 
\end{proposition} 

Now, we have
\begin{theorem}\label{mauto} The automorphism group ${\rm Aut}(\mathcal M_{s,p,n})$ of $\mathcal M_{s,p,n}$ is 
\begin{equation*}
    {\rm PGL}(s,\mathbb K)\times {\rm PGL}(n-s,\mathbb K)\,,
\end{equation*} 
except for the following cases. 
\begin{enumerate}[label={\rm(\Alph*)}]
\item {\rm (Usd\footnote{Usd stands for a discrete symmetry induced by a symmetry turns $\mathcal T_{s,p,n}$ upside down.} case)}. If $n=2s$ and $p\neq s$, 
\begin{equation*}
  {\rm Aut}(\mathcal M_{s,p,2s})=\left({\rm PGL}(s,\mathbb K)\times {\rm PGL}(n-s,\mathbb K)\right) \rtimes\mathbb Z/2\mathbb Z\,. 
\end{equation*}
\item {\rm (Dual\footnote{Dual stands for the discrete symmetry induced by the dual map of $G(p,2p)$.}  case)}. If $n=2p$ and $p\neq s$, 
\begin{equation*}
  {\rm Aut}(\mathcal M_{s,p,2p})=\left({\rm PGL}(s,\mathbb K)\times {\rm PGL}(n-s,\mathbb K)\right) \rtimes\mathbb Z/2\mathbb Z\,. 
\end{equation*}
\item {\rm (Usd+Dual\footnote{The two types of discrete symmetries in (A) and (B) coexist in (C).} case)}. If $n=2s=2p$ and $p\geq 3$,
\begin{equation*}
  {\rm Aut}(\mathcal M_{p,p,2p})=\left({\rm PGL}(s,\mathbb K)\times {\rm PGL}(n-s,\mathbb K)\right) \rtimes\mathbb Z/2\mathbb Z\rtimes\mathbb Z/2\mathbb Z\,. 
\end{equation*}
\item {\rm (Degenerate cases)}. \begin{enumerate}[label={\rm(\alph*).}]
      \item  $\mathcal M_{1,1,2}$ is a point.
      \item  $\mathcal M_{2,2,4}\cong\mathbb {P}^3$ and ${\rm Aut}(\mathcal M_{2,2,4})={\rm PGL}(4,\mathbb K)$\,.
\end{enumerate}
\end{enumerate}
 \end{theorem}

%When $s=p$, the manifold $\mathcal M_{s,p,n}$ is isomorphic to a classical object studied in algebraic geometry called the variety of complete collineations.  Study (\cite {Stu}), Severi (\cite{Sev1,Sev2}), and Van der Waerden (\cite{Van}) studied the complete conics in $\mathbb P^2$ from the perspective of enumerative geometry.  Semple (\cite{Se1, Se2}), Alguneid (\cite {Al}), and  Tyrrell (\cite{Ty}) studied the complete collineations in $\mathbb P^n$ of higher dimensions.  Vainsencher (\cite{Va}) established the smoothness and proved that  $\mathcal M_{p,p,n}$ is a wonderful variety in the sense of   De Concini and Procesi (\cite{DP}). 

Note that Brion \cite{Br6} computed the automorphism group of $\mathcal M_{p,p,2p}$, the wonderful compactification of ${\rm PGL}(p,\mathbb K)$. The automorphism group of $\mathcal{M}_{p,p,n}$, the space of complete collineations, was first determined by Massarenti \cite{M3}.
%\begin{remark}\label{dualusd}The above discrete symmetries exist in a more general setting as the following isomorphisms   (see Definitions \ref{tiso}, \ref{USD1}, \ref{musd}.\begin{equation*}\begin{split}     &{\rm DUAL}:\,\mathcal T_{s,p,n}\rightarrow\mathcal T_{s,n-p,n}\,\,\,\,{\rm and}\,\,\,\,     {\rm USD}:\mathcal T_{s,p,n}\rightarrow\mathcal T_{n-s,p,n}\,;\\     &  {\rm Dual}:\mathcal M_{s,p,n}\rightarrow \mathcal M_{s,n-p,n} \,\,\,\,{\rm and}\,\,\,\,  {\rm Usd}:\mathcal M_{s,p,n}\rightarrow\mathcal M_{n-s,p,n}\,.\end{split}\end{equation*}
%We can thus assume that $2p\leq n\leq 2s$ without loss of  generality.\end{remark}

\medskip

We now briefly outline the organization of the paper and the basic ideas for the proof. Our approach utilizes  the Mille Cr\^epes coordinate charts introduced in \cite{FW} to extract geometric information from the spaces $\mathcal T_{s,p,n}$ and $\mathcal M_{s,p,n}$. %, including $\mathbb Z$-bases for their Picard groups, 
Combined with Brion's results that the cones of effective cycles on any irreducible complete spherical variety is generated by the classes of the closures of its $B$-orbits, we thus determine the extremal rays of the effective cones, and the intersection pairings between line bundles and torus-invariant curves. The automorphism of $\mathcal T_{s,p,n}$ follows from such numerical data. The case of $\mathcal M_{s,p,n}$ is more subtle. We begin by proving that the complete linear series of any $B$-stable divisor is generated by the corresponding Pl\"ucker coordinate functions.  Next, we establish the (semi-)positivity of the anticanonical bundles of $\mathcal T_{s,p,n}$ through a detailed verification of all relevant intersection numbers with torus-invariant curves.  Combined these two results with the Kawamata--Viehweg vanishing theorem, the problem is thus reduced to that of determining special automorphisms of certain projective bundles over sub-Grassmannians.
To conclude the proof, we make extensive use of the geometry of the Grassmannians,  particularly the Pl\"ucker coordinates and Pl\"ucker relations.   

The organization of the paper is as follows. After fixing certain notations, in \S \ref{isomt} we realize $\mathcal T_{s,p,n}$ as iterated blow-ups of $G(p,n)$, and then introduce natural group actions on these spaces. In \S \ref{mcc}, we recall the Mille Cr\^epes coordinate charts for $\mathcal T_{s,p,n}$ introduced in \cite{FW}. In \S \ref{foliation}, we recall the explicit Bia{\l}ynicki-Birula decomposition on $\mathcal T_{s,p,n}$ and use it to introduce Usd and Dual isomorphisms for $\mathcal M_{s,p,n}$. In \S \ref{curanddivg}, we define the $G$-stable and $B$-stable divisors on $\mathcal T_{s,p,n}$ and $\mathcal M_{s,p,n}$, and determine the structure of various geometric invariants such as the Picard groups and the canonical bundles. In particular, we prove Lemma \ref{cls} which determines the complete linear series for $B$-stable divisors.
We enumerate all torus-invariant curves and calculate their corresponding intersection numbers in \S \ref{curanddivi}, and use this result to prove Propositions \ref{mfano}, \ref{fano} in \S \ref{curanddivp}. 
Finally, we prove Theorem \ref{auto1} in \S \ref{symt} and Theorem \ref{mauto} in \S \ref{symm}. For the reader's convenience, the intersection numbers from Section \ref{curanddivi} are compiled in Appendix \ref{CIN}.

%the morphism $\mathcal P_{s,p,n}:\mathcal T_{s,p,n}\rightarrow D^-_1$.  In , we introduce the fibration structures on  $\mathcal M_{s,p,n}$;  in \S \ref{basicif}, we prove  Proposition \ref{tfiber}, Theorem \ref{gwond} based on the geometry of the foliation on $\mathcal T_{s,p,n}$ induced by the torus action. 

%We attach five appendices to the paper.
%In Appendix \ref{section:projbc}, we assign certain local coordinate charts for the projective bundles associated with the source and the sink in $G(p,n)$ which are helpful to understand the fibration structures of $\mathcal M_{s,p,n}$ and used in the proof of Lemma \ref{fibpre}. In Appendix \ref{section:rigidmc}, we classify the induced automorphisms of the Picard group of $\mathcal T_{s,p,n}$. In Appendix \ref{section:rigidtc}, we determine the automorphism of $\mathcal M_{s,p,n}$ under the assumption that the induced automorphisms of the Picard group is the identity. We include certain numerical calculations for Delcroix's criterion in Appendix \ref{section:nc}.\medskip
\medskip

{\noindent\bf Acknowledgment.} The author is very grateful to the anonymous referees for their careful reading of the manuscript and their helpful remarks. This work was supported by National Key R\&D Program of China under Grant No.2022YFA1006700 and NSFC-12201012.

\section{Preliminaries}\label{iterated}

We first briefly recall certain notations from \cite{FW}. %For an integer $m\geq1$, we denote by $\mathbb A^{m}$ the affine scheme ${\rm Spec}\,\mathbb Z\left[x_{1},\cdots,x_{m}\right]$, and by $\mathbb P^{m}$ the projective space ${\rm Proj}\,\mathbb Z\left[x_{0},\cdots,x_{m}\right]$. 
In the paper, without loss of generality, we assume that $2p\leq n\leq 2s$, and call $r:=\min\{s,n-s,p,n-p\}$ the rank of $\mathcal T_{s,p,n}$. 

Define an index set
\begin{equation*}
\mathbb I_{p,n}:=\{(i_1,i_2,\cdots,i_p)\in\mathbb Z^p:1\leq i_p<i_{p-1}<\cdots<i_1\leq n\}.
\end{equation*}
For each $I=(i_1,\cdots,i_p)\in\mathbb I_{p,n}$, we define {\it Pl\"ucker coordinate functions} $P_{I}$ on $$\mathbb A^{pn}:={\rm Spec}\mathbb Z\left[x_{ij}(1\leq i\leq p,1\leq j\leq n)\right]$$ to be the $p\times p$-subdeterminant of $(x_{ij})$ consisting of the $i_1$-th,  $\cdots$, $i_p$-th columns.
Define \begin{equation*}\mathcal {G}(p,n):=\left\{\mathfrak p\in\mathbb A^{pn}:P_{I}\notin\mathfrak p\,\,{\rm for\,\,a\,\,certain\,\,}I\in\mathbb I_{p,n}\right\}. \end{equation*} For any partial permutation $\Delta:=(\delta_1,\delta_2,\cdots,\delta_p)$ of $(1,2,\cdots,n)$, define closed subschemes of $\mathcal G(p,n)$ \begin{equation*}U_{\Delta}:={\rm Spec}\left(\mathbb Z\left[\cdots,x_{ij},\cdots\right]/(x_{1\delta_1}-1,x_{1\delta_2},\cdots,x_{1\delta_p},\cdots,x_{p\delta_1},\cdots,x_{p\delta_p}-1)\right) \end{equation*}with the embeddings denoted by $e_{\Delta}:U_{\Delta}\xhookrightarrow{\,\,\,\,\,\,\,}\mathcal G(p,n)$. Identifying $U_{\Delta}$ with their images under $[\cdots,P_{I}\circ e_{\Delta},\cdots]_{I\in\mathbb I_{p,n}}$, we derive an atlas of the Grassmannian $G(p,n)$.  Denote by $\pi:\mathcal G(p,n)\rightarrow G(p,n)$ the natural projection. For $x\in G(p,n)$, we denote by  $\widetilde{x}$ any element in the preimage $\pi^{-1}(x)\subset\mathcal G(p,n)$. For simplicity of notation, we call $\widetilde x$ a matrix representative of $x$.

For convenience, we  write the Pl\"ucker embedding $Gr(p,E)\hookrightarrow\mathbb P(\bigwedge^{p}E)$  as $e:G(p,n)\xhookrightarrow{\,\,\,\,\,\,\,}\mathbb P^{N_{p,n}}$, where $\mathbb P^{N_{p,n}}$ is the projective space of dimension $N_{p,n}:=\frac{n!}{(n-p)!p!}-1$ with homogeneous coordinates  
$[\cdots ,z_I,\cdots]_{\,I\in\mathbb I_{p,n}}$.

For $0\leq k\leq r$ (see \ref{rank}), define index subsets  
\begin{equation}\label{ikk}
\mathbb I_{s,p,n}^{k}:=\{(i_1,\cdots,i_p)\in\mathbb Z^p:{ 1\leq i_p<\cdots<i_{k+1}\leq s\,;s+1\leq i_{k}<\cdots<i_1\leq n}\}.
\end{equation}
Consider the following linear subspace of $\mathbb {P}^{N_{p,n}}$,
\begin{equation}\label{subl}
\{[\cdots ,z_I,\cdots]_{I\in\mathbb I_{p,n}}\in\mathbb {P}^{N_{p,n}}:z_I=0\,,\,\,\forall I\notin\mathbb I_{s,p,n}^k \},\,\,\,\, 0\leq k\leq r,
\end{equation}
For convenience, we denote (\ref{subl}) by $\mathbb {P}^{N^k_{s,p,n}}$ with $N^k_{s,p,n}=|\mathbb I_{s,p,n}^{k}|-1$ and its homogeneous coordinates by  $[\cdots ,z_I,\cdots]_{I\in\mathbb I^k_{s,p,n}}$.  %By dropping the coordinates with indices not in $\mathbb I_{s,p,n}^k$, we have the natural projection $F_s^k:\mathbb {P}^{N_{p,n}}\dashrightarrow\mathbb {P}^{N^k_{s,p,n}}$.

%Under the above notation, it is clear that (\ref{kspe}) takes the form \begin{equation*}\mathcal K_{s,p,n}:=(e, f_s^0,\cdots,f_s^r):\mathbb A^{p(n-p)}\dashrightarrow\mathbb {P}^{N_{p,n}}\times\mathbb {P}^{N^0_{s,p,n}}\times\cdots\times\mathbb {P}^{N^r_{s,p,n}}, \end{equation*} where for $0\leq k\leq r$, $f_s^k:=F_s^k\circ e$. In homogeneous coordinates, \begin{equation*} \begin{split}\mathcal K_{s,p,n}(x)=\big( [\cdots ,P_I(\widetilde x),\cdots]_{I\in\mathbb I_{p,n}},[\cdots,P_I(\widetilde x),\cdots]_{I\in\mathbb I^0_{s,p,n}},\cdots,[\cdots,P_I(\widetilde x),\cdots]_{I\in\mathbb I^r_{s,p,n}}\big)\end{split}.\end{equation*}

\subsection{Isomorphisms among \texorpdfstring{$\mathcal T_{s,p,n}$}{rr}} \label{isomt}
%In this subsection, we will introduce certain isomorphisms among $\mathcal T_{s,p,n}$.

We have the following alternative construction of $\mathcal T_{s,p,n}$ as iterated blow-ups of $G(p,n)$. For $0\leq k\leq r$, define subschemes of $G(p,n)$ by 
\begin{equation}\label{sk}
S_k:=\{ x\in G(p,n):\,P_{I}=0\,\,\forall{I\in \mathbb I_{s,p,n}^k}\}\,. 
\end{equation}
For any permutation $\sigma$ of $\{0,1,\cdots,r\}$, let $g_0^{\sigma}:Y^{\sigma}_0\rightarrow G(p,n)$ be the blow-up of $G(p,n)$ along $S_{\sigma(0)}$, and inductively define $g^{\sigma}_{i+1}:Y^{\sigma}_{i+1}\rightarrow Y^{\sigma}_{i}$, $0\leq i\leq r-1$, to be the blow-up of $Y^{\sigma}_i$ along $(g^{\sigma}_0\circ g^{\sigma}_{1}\circ\cdots\circ g^{\sigma}_i)^{-1}(S_{\sigma(i+1)})$. 
It is clear that $Y_r^{\sigma}$ is independent of the choice of $\sigma$.
Moreover, for any permutation $\sigma$, there is an isomorphism $\nu_{\sigma}:\mathcal T_{s,p,n}\rightarrow Y^{\sigma}_r$ such that the following diagram commutes.
\vspace{-0.05in}
\begin{equation*}
\begin{tikzcd}
&\mathcal T_{s,p,n}\arrow[dashed,swap]{rd}{\hspace{-0.03in}\mathcal K^{-1}_{s,p,n}} \arrow{rr}{\nu_{\sigma}}&&Y^{\sigma}_{r}\arrow{dl}{\hspace{-.03in}(g^{\sigma}_0\circ\cdots\circ g^{\sigma}_{r})} \\
&&\,\,G(p,n)\,\,&\\
\end{tikzcd}\vspace{-20pt}\,.
\end{equation*}
In particular, there is a morphism $R_{s,p,n}:\mathcal T_{s,p,n}\rightarrow G(p,n)$ extending $\mathcal K^{-1}_{s,p,n}$.
It is clear that $R_{s,p,n}$ is also given by the projection 
of $\mathcal T_{s,p,n}$ to the first factor $\mathbb {P}^{N_{p,n}}$ of the ambient space $\mathbb {P}^{N_{p,n}}\times\mathbb {P}^{N^0_{s,p,n}}\times\cdots\times\mathbb {P}^{N^r_{s,p,n}}$.

Recall that every $p$-dimensional subspace $W$ of an $n$-dimensional space  $V$ determines an $(n-p)$-dimensional quotient space $V/W$. Taking the dual yields an inclusion $(V/W)^*\hookrightarrow V^*$. We thus have the canonical isomorphism
\begin{equation}\label{dual}
  {\empty}^*:\,\,\,\,Gr(p, V)\cong Gr(n-p, V^*).  
\end{equation}

In the following, we first extend (\ref{dual}) to an isomorphism between $\mathcal T_{s,p,n}$ and $\mathcal T_{s,n-p,n}$. For $0\leq k\leq p$ and  $I=(i_1,\cdots, i_p)\in\mathbb I^k_{s,p,n}$, let $I^*=(i^*_1,\cdots,i^*_{n-p})\in\mathbb I^{s-k}_{s,n-p,n}$ be such that
\begin{equation*}
    \left\{i_1,\cdots, i_p, i^*_1,\cdots,i^*_{n-p}\right\}=\left\{1,2,\cdots,n \right\}.
\end{equation*}
Define an isomorphism $g_k:\mathbb {P}^{N^{k}_{s,p,n}}\rightarrow\mathbb {P}^{N^{s-k}_{s,n-p,n}}$ by 
\begin{equation*}
\begin{split}
g_k\left(\left[\cdots,z^k_I,\cdots\right]_{I\in\mathbb I^k_{s,p,n}}\right):=\big[\cdots,\underset{\substack{\uparrow\\z^{s-k}_{I^*}}}{(-1)^{\sigma(I)}\cdot z^k_{I},}\cdots\big]_{I^*\in\mathbb I^{s-k}_{s,n-p,n}}\,,
\end{split}
\end{equation*}
where $\sigma(I)$ is the signature of the permutation $\left(\begin{matrix}
		i_1&\cdots&i_p&i^*_1&\cdots&i^*_{n-p}\\
		1&\cdots&p&p+1&\cdots&n\\
	\end{matrix}\right)$, and $\left[\cdots,z^k_I,\cdots\right]_{I\in\mathbb I^k_{s,p,n}}$ and $\left[\cdots,z^{s-k}_{I^*},\cdots\right]_{I^*\in\mathbb I^{s-k}_{s,n-p,n}}$ are the homogeneous coordinates for $\mathbb {P}^{N^{k}_{s,p,n}}$ and  $\mathbb {P}^{N^{s-k}_{s,n-p,n}}$ respectively.
Now we define an isomorphism 
\begin{equation*}
    \widetilde{{\rm DUAL}}:\,\,G(p,n)\times\mathbb {P}^{N^{0}_{s,p,n}}\times\cdots\times\mathbb {P}^{N^p_{s,p,n}}\rightarrow G(n-p,n)\times\mathbb {P}^{N^{0}_{s,n-p,n}}\times\cdots\times\mathbb {P}^{N^{n-p}_{s,n-p,n}}\,
\end{equation*}
as follows. Let $a=\left(x,\left[\cdots,z^0_I,\cdots\right]_{I\in\mathbb I^0_{s,p,n}}, \cdots, \left[\cdots,z^p_I,\cdots\right]_{I\in\mathbb I^p_{s,p,n}}\right)\in G(p,n)\times\mathbb {P}^{N^{0}_{s,p,n}}\times\cdots\times\mathbb {P}^{N^p_{s,p,n}}$. The $G(n-p,n)$-component of $\widetilde{{\rm DUAL}}(a)$ is $x^*$. Define the $\mathbb {P}^{N^k_{s,n-p,n}}$-component of $\widetilde{{\rm DUAL}}(a)$ by
\begin{equation*}
    \left[\cdots,z^{k}_{I^*}\left(\widetilde{{\rm DUAL}}(a)\right),\cdots\right]_{I^*\in\mathbb I^{k}_{s,n-p,n}}=g_{s-k}\left(\left[\cdots,z^{s-k}_I,\cdots\right]_{I\in\mathbb I^{s-k}_{s,p,n}}\right)\,,
\end{equation*}
when $\max\{0,n-s-p\}\leq k\leq \min\{n-s,n-p\}$. If $n-s < k \leq n-p$ or $0 \leq k < n-s-p$, then $\mathbb{I}^{k}{s,n-p,n} = \emptyset$. In this case, by convention, the projective space $\mathbb{P}^{N^k{s,n-p,n}}$ is a single point $c_k$, and we define the $\mathbb{P}^{N^k_{s,n-p,n}}$-component of $\widetilde{\mathrm{DUAL}}(a)$ to be $c_k$.

\begin{definitionlemma}\label{tiso}
The restriction of $\widetilde{{\rm DUAL}}$  induces an isomorphism ${\rm DUAL}$ from  $\mathcal T_{s,p,n}$ to $\mathcal T_{s,n-p,n}$ such that the following diagram commutes.
\vspace{-.03in}
\begin{equation*}%\label{DUAL}
\begin{tikzcd} &\mathcal T_{s,p,n}\arrow{r}{\rm DUAL}&\mathcal T_{s,n-p,n}\arrow{d}{R_{s,n-p,n}} \\ &G(p,n) \arrow{r}{{\empty}^*}\arrow[leftarrow]{u}{R_{s,p,n}}&G(n-p,n)\\ \end{tikzcd}\,\,.\vspace{-20pt} \end{equation*}
\end{definitionlemma}
{\bf\noindent Proof of Lemma \ref{tiso}.}
It suffices to prove that the image of $\mathcal T_{s,p,n}$ under $\widetilde{{\rm DUAL}}$ is contained in $\mathcal T_{s,n-p,n}$. 

For each $x\in Gr(p, V)\cong G(p,n)$, let $\widetilde x$ be a matrix representative of $x$, and $\widetilde{x^*}$ a matrix representative of $x^*\in Gr(n-p, V^*)\cong G(n-p,n)$. Then in terms of matrix multiplication, the dual map is characterized by
\begin{equation*}
    \widetilde x\cdot \left(\widetilde{x^*}\right)^T=0\,,
\end{equation*}
where  $\left(\widetilde{x^*}\right)^T$ is the transpose of the matrix $\widetilde {x^*}$,  $0$ is the $p\times (n-p)$ zero matrix.

Consider the following open set of $G(p,n)$,
\begin{equation*}
   U:=\left\{\left.  \left(I_{p\times p} \hspace{-0.13in}\begin{matrix}
  &\hfill\tikzmark{c1}\\
  &\hfill\tikzmark{d1}
  \end{matrix}\,\,\, A \right)\right\vert_{}
   A\,\,{\rm  is\,\, a\,\,} p\times (n-p)\,\,{\rm matrix}\,\right\}.
  \tikz[remember picture,overlay]   \draw[dashed,dash pattern={on 4pt off 2pt}] ([xshift=0.5\tabcolsep,yshift=7pt]c1.north) -- ([xshift=0.5\tabcolsep,yshift=-2pt]d1.south);
\end{equation*}
Denote by  $U^*$ the image of $U$ under the dual map ${\empty}^*$. Then,
\begin{equation*}
   U^*=\left\{\left.  \left(-A^T \hspace{-0.13in}\begin{matrix}
  &\hfill\tikzmark{c2}\\
  &\hfill\tikzmark{d2}
  \end{matrix}\,\,\, I_{(n-p)\times (n-p)} \right)\right\vert_{}
   A\,\,{\rm  is\,\, a\,\,} p\times (n-p)\,\,{\rm matrix}\,\right\},
  \tikz[remember picture,overlay]   \draw[dashed,dash pattern={on 4pt off 2pt}] ([xshift=0.5\tabcolsep,yshift=7pt]c2.north) -- ([xshift=0.5\tabcolsep,yshift=-2pt]d2.south);
\end{equation*}
and the dual map takes the form
\begin{equation}\label{dtrans}
  { \left(I_{p\times p} \hspace{-0.13in}\begin{matrix}
  &\hfill\tikzmark{c1}\\
  &\hfill\tikzmark{d1}
  \end{matrix}\,\,\, A \right)}^*=\left(-A^T \hspace{-0.13in}\begin{matrix}
  &\hfill\tikzmark{c2}\\
  &\hfill\tikzmark{d2}
  \end{matrix}\,\,\, I_{(n-p)\times (n-p)} \right).
  \tikz[remember picture,overlay]   \draw[dashed,dash pattern={on 4pt off 2pt}] ([xshift=0.5\tabcolsep,yshift=7pt]c1.north) -- ([xshift=0.5\tabcolsep,yshift=-2pt]d1.south);
  \tikz[remember picture,overlay]   \draw[dashed,dash pattern={on 4pt off 2pt}] ([xshift=0.5\tabcolsep,yshift=7pt]c2.north) -- ([xshift=0.5\tabcolsep,yshift=-2pt]d2.south);
\end{equation}

Computing the images of $U$ and $U^*$ under the birational maps $\mathcal K_{s,p,n}$ and  $\mathcal K_{s,n-p,n}$ respectively, we can conclude Lemma \ref{tiso}. \,\,\,$\endpf$
\medskip

In a similar way, we will introduce another isomorphism as follows. 
Let $\mathfrak E=(e_{ij})$ be an $n\times n$ anti-diagonal matrix with $e_{ij}=1$ for $j+i=n+1$ and $0$ otherwise.
For $0\leq k\leq p$ and $I=(i_1,\cdots, i_p)\in\mathbb I^k_{s,p,n}$, define
\begin{equation*}
 I^*:=(n+1-i_p,n+1-i_{p-1},\cdots,n+1-i_2,n+1-i_1)\in\mathbb I^{p-k}_{n-s,p,n}\,.
\end{equation*}
For $0\leq k\leq p$, define an isomorphism $g_k:\mathbb {P}^{N^{k}_{s,p,n}}\rightarrow\mathbb {P}^{N^{p-k}_{n-s,p,n}}$ by
$z^{p-k}_{I^*}=z^k_{I}$, where $[\cdots,z^k_I,\cdots]_{I\in\mathbb I^k_{s,p,n}}$ and $[\cdots,z^{p-k}_{I^*},\cdots]_{I^*\in\mathbb I^{p-k}_{n-s,p,n}}$ are the homogeneous coordinates for $\mathbb {P}^{N^{k}_{s,p,n}}$ and  $\mathbb {P}^{N^{p-k}_{n-s,p,n}}$ respectively..
We define an isomorphism 
\begin{equation*}
    \widetilde{\rm USD}:\,\,G(p,n)\times\mathbb {P}^{N^{0}_{s,p,n}}\times\cdots\times\mathbb {P}^{N^p_{s,p,n}}\rightarrow G(p,n)\times\mathbb {P}^{N^{0}_{n-s,p,n}}\times\cdots\times\mathbb {P}^{N^p_{n-s,p,n}}\,
\end{equation*}
as follows.  Let $a=\left(x,\left[\cdots,z^0_I,\cdots\right]_{I\in\mathbb I^0_{s,p,n}}, \cdots, \left[\cdots,z^p_I,\cdots\right]_{I\in\mathbb I^p_{s,p,n}}\right)\in G(p,n)\times\mathbb {P}^{N^{0}_{s,p,n}}\times\cdots\times\mathbb {P}^{N^p_{s,p,n}}$.
The $G(p,n)$-component of $\widetilde{{\rm USD}}(a)$ is $\mathfrak E\cdot x$. Define the $\mathbb {P}^{N^k_{n-s,p,n}}$-component of $\widetilde{{\rm USD}}(a)$ by
\begin{equation*}
\left[\cdots,z^{k}_{I^*}\left(\widetilde{{\rm USD}}(a)\right),\cdots\right]_{I^*\in\mathbb I^{k}_{n-s,p,n}}=g_{p-k}\left(\left[\cdots,z^{p-k}_I,\cdots\right]_{I\in\mathbb I^{p-k}_{s,p,n}}\right)\,,
\end{equation*}
when $\max\{0,s+p-n\}\leq k\leq \min\{p,s\}$. If $s<k\leq p$ or $0\leq k<s+p-n$, we define the $\mathbb {P}^{N^k_{n-s,p,n}}$-component of $\widetilde{{\rm USD}}(a)$ to be the single point of $\mathbb {P}^{N^k_{n-s,p,n}}$.

Similarly, we have
\begin{definitionlemma}\label{USD1}
$\widetilde {\rm USD}$ induces an isomorphism ${\rm USD}:\mathcal T_{s,p,n}\rightarrow\mathcal T_{n-s,p,n}$ such that the following diagram commutes.
\vspace{-.05in}
\begin{equation*}%\label{USD}
\begin{tikzcd} &\mathcal T_{s,p,n}\arrow{r}{\rm USD}&\mathcal T_{n-s,p,n}\arrow{d}{R_{n-s,p,n}} \\ &G(p,n) \arrow{r}{\mathfrak E}\arrow[leftarrow]{u}{R_{s,p,n}}&G(p,n)\\ \end{tikzcd}\,\,.\vspace{-20pt} \end{equation*}  
\end{definitionlemma}

%\begin{remark}${\rm DUAL}$ is an automorphism of $\mathcal T_{s,p,2p}$ and USD is an automorphism of $\mathcal T_{s,p,2s}$. \end{remark}

\subsection{Mille Cr\^epes Coordinates}\label{mcc}
In this section, we recall the smooth atlas for $\mathcal T_{s,p,n}$ defined in \cite{FW}. 
The idea is to locally represent $\mathcal T_{s,p,n}$ as a sequence of blow-ups such that the transition between successive blow-ups is essentially the Gaussian elimination process. When the coordinate charts are constructed in reverse, all resulting matrices are summed up to consolidate the intermediate steps, which resembles stacking layers of paper-thin crepes and ganache on top of each other (for this reason we call such coordinates {\it Mille Cr\^epes}).

For convenience, 
for any $x\in G(p,n)$, we denote by $P_I(x)$ the functions $P_I(\widetilde x)$, since they are only used in such a way that the choice of $\widetilde x$ is irrelevant.

We define the Mille Cr\^epes coordinate charts up to group actions as follows. For $0\leq l\leq r$, define affine open subscheme $U_l\subset G(p,n)$ by
\begin{equation}\label{ul}
U_l:=\{\,\,\,\,\underbracedmatrixll{Z\\Y}{s-p+l\,\,\rm columns}
  \hspace{-.3in}\begin{matrix}
  &\hfill\tikzmark{a}\\
  &\hfill\tikzmark{b}  
  \end{matrix} \,\,\,\,\,
  \begin{matrix}
  0\\
I_{(p-l)\times(p-l)}\\
\end{matrix}%\hspace{0.1in}
\begin{matrix}
  &\hfill\tikzmark{c}\\
  &\hfill\tikzmark{d}
  \end{matrix}\begin{matrix}
  &\hfill\tikzmark{g}\\
  &\hfill\tikzmark{h}
  \end{matrix}\,\,\,\,
\begin{matrix}
I_{l\times l}\\
0\\
\end{matrix}\hspace{0.01in}
\begin{matrix}
  &\hfill\tikzmark{e}\\
  &\hfill\tikzmark{f}\end{matrix}\hspace{-.3in}\underbracedmatrixrr{X\\W}{(n-s-l)\,\,\rm columns}\,\,\,\,\}
  \tikz[remember picture,overlay]   \draw[dashed,dash pattern={on 4pt off 2pt}] ([xshift=0.5\tabcolsep,yshift=7pt]a.north) -- ([xshift=0.5\tabcolsep,yshift=-2pt]b.south);\tikz[remember picture,overlay]   \draw[dashed,dash pattern={on 4pt off 2pt}] ([xshift=0.5\tabcolsep,yshift=7pt]c.north) -- ([xshift=0.5\tabcolsep,yshift=-2pt]d.south);\tikz[remember picture,overlay]   \draw[dashed,dash pattern={on 4pt off 2pt}] ([xshift=0.5\tabcolsep,yshift=7pt]e.north) -- ([xshift=0.5\tabcolsep,yshift=-2pt]f.south);\tikz[remember picture,overlay]   \draw[dashed,dash pattern={on 4pt off 2pt}] ([xshift=0.5\tabcolsep,yshift=7pt]g.north) -- ([xshift=0.5\tabcolsep,yshift=-2pt]h.south);
\end{equation}
with coordinates
\begin{equation*}%\label{ulx}
\begin{split}
&Z:=(\cdots,z_{ij},\cdots)_{1\leq i\leq l,\,1\leq j\leq s-p-l},\ \ \ \ \ \ X:=(\cdots,x_{ij},\cdots)_{1\leq i\leq l,\,s+l+1\leq j\leq n},\\
&Y:=(\cdots,y_{ij},\cdots)_{l+1\leq i\leq p,\,1\leq j\leq s-p-l},\ \ \ W:=(\cdots,w_{ij},\cdots)_{l+1\leq i\leq p,\,s+l+1\leq j\leq n}.
\end{split}    
\end{equation*}

For $0\leq l\leq r$, we define index sets 
\begin{equation*}
\mathbb J_l:=\left\{\left(
\begin{matrix}
i_1&\cdots&i_r\\
j_1&\cdots&j_r\\
\end{matrix}\right)\rule[-.38in]{0.01in}{.82in}\,\,\begin{matrix}
(i_{1},\cdots,i_{r-l})\,\, {\rm is\,\,a\,\,partial\,\,permutation\,\,of\,\,}(l+1,\cdots,\,p)\\
(i_{r-l+1},\cdots,i_{p})\,\, {\rm is\,\,a\,\,permutation\,\,of\,\,}(1,\cdots,\,l)\\
(j_{1},\cdots,j_{r-l})\,\, {\rm is\,\,a\,\,partial\,\,permutation\,\,of\,\,}(s+l+1,\cdots,\,n)\\(j_{r-l+1},\cdots,j_{r})\,\, {\rm is\,\,a\,\,partial\,\,permutation\,\,of\,\,}(1,\cdots,\,s-p+l)
\end{matrix}\right\}.
\end{equation*}
For each $\tau=\left(\begin{matrix}
i_1&\cdots&i_r\\
j_1&\cdots&j_r\\
\end{matrix}\right)\in\mathbb J_l$, set $\mathbb {A}^{p(n-p)}:={\rm Spec}\,\mathbb Z[\overrightarrow A,\widetilde X,\widetilde Y,\overrightarrow B^{1},\cdots,\overrightarrow B^{r}]$, where %$\left(\widetilde X,\widetilde Y,\overrightarrow B^1,\cdots,\overrightarrow B^p\right)$ defined as follows.
\begin{equation*}%\label{ulu}
\begin{split}
&\overrightarrow A:=\left((b_{i_kj_k})_{1\leq k\leq r-l},(a_{i_kj_k})_{r-l+1\leq k\leq r}\right),\\
&\widetilde X:=(\cdots,x_{ij},\cdots)_{1\leq i\leq l,\,s+l+1\leq j\leq n},\ \ \widetilde 
Y:=(\cdots,y_{ij},\cdots)_{l+1\leq i\leq p,\,1\leq j\leq s-p-l},       
\end{split}
\end{equation*}
%\begin{equation}\label{ulu}\widetilde X:=\left(\begin{matrix}x_{1(s+l+1)}&\cdots &x_{1n}\\\vdots&\ddots&\vdots\\x_{l(s+l+1)}&\cdots &x_{ln}\\\end{matrix}\right)\,\,\,{\rm and}\,\,\,\widetilde Y:=\left(\begin{matrix}y_{(l+1)1}&\cdots& y_{(l+1)(s-p+l)}\\\vdots&\ddots&\vdots\\y_{p1}&\cdots& y_{p(s-p+l)}\\\end{matrix}  \right)  \,;\,\end{equation}
for $1\leq k\leq r-l$,
\begin{equation*}
\overrightarrow B^{k}:=\left( \ (\xi^{(k)}_{i_kj})_{s+l+1\leq j\leq n,\,j\neq j_1,j_2,\cdots,j_k},\ (\xi^{(k)}_{ij_k})_{l+1\leq i\leq p,\,i\neq i_1,i_2,\cdots,i_k}\right),
\end{equation*}
and for $r-l+1\leq k\leq r$,
\begin{equation*}
\overrightarrow B^{k}:=\left(\ (\xi^{(k)}_{i_kj})_{1\leq j\leq s-p+l,\,j\neq j_{r-l+1},j_{r-l+2},\cdots,j_k},\ (\xi^{(k)}_{ij_k})_{1\leq i\leq l,\,i\neq i_{r-l+1},i_{r-l+2},\cdots,i_k}\right). 
\end{equation*}
%\begin{equation}\label{qbpb}\begin{split}&\overrightarrow B^{k}:=\left(a_{i_{k}j_{k}},\xi^{(k)}_{i_{k}1},\xi^{(k)}_{i_{k}2},\cdots,\widehat{\xi^{(k)}_{i_{k}j_{p-l+1}}},\cdots,\widehat{\xi^{(k)}_{i_{k}j_{p-l+2}}},\cdots,\widehat{\xi^{(k)}_{i_{k}j_{k}}},\cdots,\xi^{(k)}_{i_{k}(s-p+l)},\right.\\&\,\,\,\,\,\,\,\,\,\,\,\,\,\,\,\,\,\,\,\,\,\,\,\,\,\,\,\,\,\,\left.\xi^{(k)}_{1j_{k}},\xi^{(k)}_{2j_{k}},\cdots,\widehat{\xi^{(k)}_{i_{p-l+1}j_{k}}},\cdots,\widehat{\xi^{(k)}_{i_{p-l+2}j_{k}}},\cdots,\cdots,\widehat{\xi^{(k)}_{i_{k}j_{k}}},\cdots,\xi^{(k)}_{lj_{k}}\right)\,.\,\,\,\,\,\,\,\,\,\,\,\,\,\\\end{split}\end{equation}
Define a map $\Gamma_l^{\tau}:\mathbb A^{p(n-p)}\rightarrow U_l$ by
\begin{equation*}
\left(
\begin{matrix}
\sum\limits_{k=r-l+1}^r\left(\prod\limits_{t=r-l+1}^{k}a_{i_{t}j_t}\right)\cdot\Xi_k^T\cdot\Omega_k &0_{l\times(p-l)}&I_{l\times l}&\widetilde X\\ \widetilde Y&I_{(p-l)\times(p-l)}&0_{(p-l)\times l}&\sum\limits_{k=1}^{r-l}\left(\prod\limits_{t=1}^{k}b_{i_tj_t}\right)\cdot\Xi_k^T\cdot\Omega_k\\
\end{matrix}\right)\,.     % \end{split}
\end{equation*}
Here, for $1\leq k\leq r-l$, $\Xi_k:=
\left(v_{l+1}^k,\cdots,v_p^k\right)$ with
\begin{equation*}%\label{w3}
     v_t^k=\left\{\begin{array}{ll}
    \xi^{(k)}_{tj_k} \,\,\,\,\,\,\,\,\,\,\,\,\,\,\,\,\,\,\,\,\,\,\,\,\,\,\,\,\,\, & t\in\{l+1,l+2,\cdots,p\}\backslash\{i_1,i_2,\cdots,i_k\}\,\,\,\,\,\,\,\,\,\,\,\,\,\,\,\,\, \,\,\,\,\,\,\,\,\,\,\,\,\,\,\,\,\,\\
    0& t\in\{i_1,i_2,\cdots,i_{k-1}\}\,\,\,\,\,\,\, \\
  
    1 \,\,\,\,\,\,\,\,\,\,\,\,\,\,\,\,\,\,\,\,\,\,\,\,\,\,\,\,\,\, &t=i_k\,\,\,\,\,\,\,\,\,\,\,\,\,\,\,\,\, \,\,\,\,\,\,\,\,\,\,\,\,\,\,\,\,\,\\
    \end{array}\right.,
\end{equation*}
and $\Omega_k:=
\left(w_{s+l+1}^k,\cdots,w_n^k\right)$ with
\begin{equation*}%\label{w4}
     w_t^k=\left\{\begin{array}{ll}
    \xi^{(k)}_{i_kt} \,\,\,\,\,\,\,\,\,\,\,\,\,\,\,\,\,\,\,\,\,\,\,\,\,\,\,\,\,\, & t\in\{s+l+1,s+l+2,\cdots,n\}\backslash\{j_1,j_2,\cdots,j_k\}\,\,\\
    0& t\in\{j_1,j_2,\cdots,j_{k-1}\}\,\,\,\,\,\,\, \\
  
    1 \,\,\,\,\,\,\,\,\,\,\,\,\,\,\,\,\,\,\,\,\,\,\,\,\,\,\,\,\,\, &t=j_k\,\,\,\,\,\,\,\,\,\\
    \end{array}\right.;
\end{equation*}
for $r-l+1\leq k\leq r$,  $\Xi_k:=
\left(v_{1}^k,\cdots,v_l^k\right)$ with
\begin{equation*}%\label{w5}
     v_t^k=\left\{\begin{array}{ll}
    \xi^{(k)}_{tj_k} \,\,\,\,\,\,\,\,\,\,\,\,\,\,\,\,\,\,\,\,\,\,\,\,\,\,\,\,\,\, & t\in\{1,2,\cdots,l\}\backslash\{i_{r-l+1},i_{r-l+2},\cdots,i_k\}\,\,\,\,\,\,\,\,\,\,\,\,\,\,\,\,\, \,\,\,\,\,\,\,\,\,\,\,\,\,\,\,\,\,\\
    0& t\in\{i_{r-l+1},i_{r-l+2},\cdots,i_{k-1}\}\,\,\,\,\,\,\, \\
  
    1 \,\,\,\,\,\,\,\,\,\,\,\,\,\,\,\,\,\,\,\,\,\,\,\,\,\,\,\,\,\, &t=i_k\,\,\,\,\,\,\,\,\,\,\,\,\,\,\,\,\, \,\,\,\,\,\,\,\,\,\,\,\,\,\,\,\,\,\\
    \end{array}\right.,
\end{equation*}
and $\Omega_k:=
\left(w_{1}^k,\cdots,w_{s-p+l}^k\right)$ with
\begin{equation*}%\label{w2}
     w_t^k=\left\{\begin{array}{ll}
    \xi^{(k)}_{i_kt} \,\,\,\,\,\,\,\,\,\,\,\,\,\,\,\,\,\,\,\,\,\,\,\,\,\,\,\,\,\, & t\in\{1,2,\cdots,s-p+l\}\backslash\{j_{r-l+1},j_{r-l+2},\cdots,j_k\} \,\,\,\,\,\,\,\,\,\,\,\,\,\,\,\,\,\\
    0& t\in\{j_{r-l+1},j_{r-l+2},\cdots,j_{k-1}\}\,\,\,\,\,\,\, \\
  
    1 \,\,\,\,\,\,\,\,\,\,\,\,\,\,\,\,\,\,\,\,\,\,\,\,\,\,\,\,\,\, &t=j_k\,\,\,\,\,\,\,\,\,\,\,\,\,\\
    \end{array}\right..
\end{equation*}

We define a rational map   $J_l^{\tau}:\mathbb A^{p(n-p)}\dashrightarrow\mathbb {P}^{N_{p,n}}\times\mathbb {P}^{N^0_{s,p,n}}\times\cdots\times\mathbb {P}^{N^r_{s,p,n}}$ by $J_l^{\tau}:=\mathcal K_{s,p,n}\circ \Gamma_l^{\tau}$.
We can prove that the rational map  $J_l^{\tau}$ is an embedding of $\mathbb A^{p(n-p)}$. Denote by $\mathcal A^{\tau}$ the image of $\mathbb A^{p(n-p)}$ under $J_l^{\tau}$.  It is clear that  $\left\{\left(\mathcal A^{\tau},(J_l^{\tau})^{-1}\right)\right\}_{\tau\in\mathbb J_l}$  is a system of coordinate charts of $R^{-1}_{s,p,n}(U_l)$, where $(J_l^{\tau})^{-1}:\mathcal A^{\tau}\rightarrow \mathbb A^{p(n-p)}$ is the inverse map. We can further prove that $\bigcup_{\tau\in\mathbb J_l}\mathcal A^{\tau}= R_{s,p,n}^{-1}(U_l)$, and hence $\left\{\left(\mathcal A^{\tau},(J_l^{\tau})^{-1}\right)\right\}_{\tau\in\mathbb J_l}$  is an atlas for $R^{-1}_{s,p,n}(U_l)$. \begin{definition}
We call $\left\{\left(\mathcal A^{\tau},(J_l^{\tau})^{-1}\right)\right\}_{\tau\in\mathbb J_l}$ the {\it  Mille Cr\^epes coordinate charts} of $R^{-1}_{s,p,n}(U_l)$, or simply the {\it  Mille Cr\^epes coordinates}.   
\end{definition}

\section{Invariant Divisors in \texorpdfstring{$\mathcal T_{s,p,n}$ }{jj} and \texorpdfstring{$\mathcal M_{s,p,n}$ }{jj}} \label{curanddiv}

\subsection{Explicit Bia{\l}ynicki-Birula decomposition on \texorpdfstring{$\mathcal T_{s,p,n}$}{hh}} \label{foliation}
%Bia{\l}ynicki-Birula \cite{Bi} decomposed a smooth projective scheme $X$ under an algebraic torus action into a disjoint union of stable/unstable subschemes fibered over connected components of the fixed point scheme. Exploiting such idea in the Mille Cr\^epes coordinate charts, we investigate the geometry of $\mathcal T_{s,p,n}$ and $\mathcal M_{s,p,n}$ in this section. 
In this section, we recall the explicit Bia{\l}ynicki-Birula decomposition on $\mathcal T_{s,p,n}$ following \cite{FW}.

Note that there is a unique ${\rm GL}(s,\mathbb K)\times {\rm GL}(n-s,\mathbb K)$-action $\Delta_{s,p,n}$ on $\mathcal T_{s,p,n}$ which is equivariant with respect to 
the natural ${\rm GL}(s,\mathbb K)\times {\rm GL}(n-s,\mathbb K)$ action on $G(p,n)$ by the matrix multiplication.  In particular, 
the $\mathbb G_m$-action on $E=E_1\oplus E_2$ with rank $E=n$ defined by \begin{equation}\label{lt}
t\cdot(e_1,e_2)\mapsto(e_1,t e_2),\ \ 
\,\,(e_1,e_2)\in E_1\oplus E_2 
\end{equation} induces a $\mathbb G_m$-action on $\mathcal T_{s,p,n}$.

For $0\leq k\leq r$, define subschemes of $G(p,n)$ 
\begin{equation*}%\label{vani}
\begin{split}
&\mathcal V_{(p-k,k)} :=\left\{\left. \left(
\begin{matrix}
0&X\\
Y&0\\
\end{matrix}\right)\right\vert_{}\begin{matrix}
X\,\,{\rm is\,\,an\,\,}k\times (n-s)\,\,{\rm matrix\,\,of\,\,rank}\,\,k\,\\
Y\,\,{\rm is\,\,a\,\,}(p-k)\times s\,\,{\rm matrix\,\,of\,\,rank}\,\,(p-k)\\
\end{matrix}
\right\}\,,\\
&\mathcal V_{(p-k,k)}^+:= \left\{\left.\left(
\begin{matrix}
0&X\\
Y&W\\
\end{matrix}\right)\right\vert_{}\begin{matrix}
X\,\,{\rm is\,\,an\,\,}k\times (n-s)\,\,{\rm matrix \,\,of\,\,rank\,\,}k\,\\
Y\,\,{\rm is\,\,a\,\,}(p-k)\times s\,\,{\rm matrix \,\,of\,\,rank\,\,}(p-k)\,\\
\end{matrix}\right\},\\
&\mathcal V_{(p-k,k)}^-:= \left\{\left.\left(
\begin{matrix}
Z&X\\
Y&0\\
\end{matrix}\right)\right\vert_{}{\begin{matrix}
X\,\,{\rm is\,\,an\,\,}k\times (n-s)\,\,{\rm matrix \,\,of\,\,rank\,\,}k\,\\
Y\,\,{\rm is\,\,a\,\,}(p-k)\times s\,\,{\rm matrix \,\,of\,\,rank\,\,}(p-k)\,\\
\end{matrix}}\right\}\,.    
\end{split}
\end{equation*}
The connected components of the fixed point scheme of $\mathcal T_{s,p,n}$ under the $\mathbb G_m$-action defined by (\ref{lt}) are given by
\begin{equation}\label{11}
\mathcal D_{(p-k,k)}:=R_{s,p,n}^{-1}(\mathcal V_{(p-k,k)}), \,\,0\leq k\leq r.
\end{equation}
 By taking the Zariski closures of $\mathcal D^{\pm}_{(p-k,k)}:=R_{s,p,n}^{-1}(\mathcal V^{\pm}_{(p-k,k)})$, we define
\begin{equation*}%\label{boundary}
D_{k}^-:=\overline{\mathcal D_{(p-k+1,k-1)}^-},\,\,D_{k}^+:=\overline{\mathcal D_{(p-r+k-1,r-k+1)}^+},\,\,1\leq k\leq r.  
\end{equation*}
By \cite[Lemma 4.3]{FW}, $D_{1}^+$, $\cdots$, $D_{r}^+$, $D_{1}^-$,  $\cdots$, $D_{r}^-$ are distinct smooth divisors with simple normal crossings. In particular, $D_1^-=\mathcal D_{(p,0)}$,  $D_1^+=\mathcal D_{(p-r,r)}$.
\medskip

Recall \cite[Lemma 4.9]{FW} that (\ref{B}) induces an isomorphism $\mathcal L_{s,p,n}:=\mathcal P_{s,p,n}\big|_{D_1^+}$ from the sink $D_1^+$ to the source $D_1^-$. Then, we have the following two extra isomorphisms of $\mathcal M_{s,p,n}$ induced from that of $\mathcal T_{s,p,n}$. 
\begin{definition}\label{musd}
Define an isomorphism ${\rm Usd}:\mathcal M_{s,p,n}\rightarrow \mathcal M_{n-s,p,n}$ by
\begin{equation*}
{\rm Usd}:=\mathcal L_{n-s,p,n}\circ{\rm USD}\big|_{\mathcal M_{s,p,n}}\,,
\end{equation*}
and an isomorphism ${\rm Dual}:\mathcal M_{s,p,n}\rightarrow \mathcal M_{s,n-p,n}$ by
\begin{equation*}
{\rm Dual}:=\mathcal L_{s,n-p,n}\circ{\rm DUAL}\big|_{\mathcal M_{s,p,n}}\,,
\end{equation*}
where ${\rm USD}$, ${\rm DUAL}$ are defined in \S \ref{isomt}.
\end{definition}

\subsection{The \texorpdfstring{$G$-stable and $B$-stable divisors}{rr}} \label{curanddivg}
For convenience, we denote by
$G$ the group ${\rm GL}(s,\mathbb K)\times {\rm GL}(n-s,\mathbb K)$. Define a Borel subgroup $B$  of $G$ by \begin{equation*}
B:=\left\{\left. \left(
\begin{matrix}
g_1&0\\
0&g_2\\
\end{matrix}\right)\right\vert_{} 
\begin{matrix}
g_1\in {\rm GL}(s,\mathbb K) \,\, {\rm is\,\, a\,\, lower\,\, triangular\,\, matrix\,};\,\,\,\,\\
g_2\in {\rm GL}(n-s,\mathbb K)\,\, {\rm is\,\, an \,\, upper \,\, triangular \,\, matrix}  
\end{matrix}
\right\}.
\end{equation*} Let $T$  be  a maximal torus of $B$ defined by \begin{equation*}%\label{maxt}
\left\{\left. \left(
\begin{matrix}
g_1&0\\
0&g_2\\
\end{matrix}\right)\right\vert_{} 
\begin{matrix}
g_1\in {\rm GL}(s,\mathbb K),\,g_2\in {\rm GL}(n-s,\mathbb K)
\,\, {\rm are\,\, diagonal } \end{matrix}
\right\}.
\end{equation*}
It is clear that $\mathcal T_{s,p,n}$ is a spherical $G$-variety for $B$ has an open orbit.

Define irreducible divisors of $G(p,n)$ by $b_{k}:=\left\{x\in G(p,n) |\,P_{I_k}(x)=0\,\right\}$, $0\leq k\leq r$, where 
\begin{equation}\label{I_k}
I_k:=(s+k,s+k-1,\cdots,s-p+k+1)\in\mathbb I^k_{s,p,n}.   
\end{equation}
Let $B_{k}\subset\mathcal T_{s,p,n}$,
$0\leq k\leq r$,  be the strict transformation of $b_k$ under the blow-up $R_{s,p,n}:\mathcal T_{s,p,n}\rightarrow G(p,n)$. It is easy to verify that
\begin{lemma}\label{gb}
If an irreducible divisor $\mathfrak D$ of $\mathcal T_{s,p,n}$ is $G$-invariant (resp. $B$-invariant), then 
\begin{equation*}
\mathfrak D\in \{D_{1}^-,  \cdots, D_{r}^-, D_{1}^+, \cdots, D_{r}^+\}\,\,(resp.\,\,\{D_{1}^-, \cdots, D_{r}^-, D_{1}^+, \cdots, D_{r}^+, B_0, \cdots, B_r\}\,).
\end{equation*}
\end{lemma}
Hence, we call $D_{1}^-, \cdots, D_{r}^-,  D_{1}^+, \cdots, D_{r}^+$ the $G$-stable divisors, and $B_{0},\cdots, B_r$ the $B$-stable divisors of $\mathcal T_{s,p,n}$. 
\begin{remark}\label{coin}
When $p=n-s$, $B_r=D^-_r$; when  $p=s$, $B_0=D^+_r$. The notion of $B$-stable divisors in this paper differs slightly from that in the theory of spherical varieties. 
$B_0$ (resp. $B_r$) is not a $B$-stable divisor when $p=s$ (resp. $p=n-s$) therein. 
\end{remark}

\begin{lemma}\label{excep}
Let $E$ be the exceptional divisor of the blow-up $R_{s,p,n}$. Then, 
\begin{equation*}%\label{te}
E=\left\{ \begin{array}{ll}
    D_{1}^++D_{2}^++\cdots+D_{r}^++D_{1}^-+D_{2}^-+\cdots+D_{r}^-&{\rm when}\,\,p<s\,\,{\rm and}\,\,n-s\neq p\\
    D_{1}^++D_{2}^++\cdots+D_{r}^++D_{1}^-+D_{2}^-+\cdots+D_{r-1}^-&{\rm when}\,\,n-s=p<s\\
    D_{1}^++D_{2}^++\cdots+D_{r-1}^++D_{1}^-+D_{2}^-+\cdots+D_{r-1}^-&{\rm when}\,\,n-s=p=s\,\,\\
    \end{array}\right..
\end{equation*}
\end{lemma}

{\noindent\bf Proof of Lemma \ref{excep}.} Recall \cite[Lemma 4.2]{FW} that $S_k$ defined in (\ref{sk}), $0\leq k\leq r$, is the scheme-theoretic union of  $\overline{\mathcal V_{(p-k-1,k+1)}^+}$ and  $\overline{\mathcal V_{(p-k+1,k-1)}^-}$, where $\overline{\mathcal V_{(p-k,k)}^{\pm}}$ is the Zariski closure of $\mathcal V_{(p-k,k)}^{\pm}$ and  $\mathcal V_{(p-r-1,r+1)}^+=\mathcal V_{(p+1,-1)}^-=\emptyset$ by convention. By (\ref{11}), it is clear that the support of $E$ is contained in the union of $D_{1}^-, \cdots, D_{r}^-, D_{1}^+, \cdots, D_{r}^+$.

Similarly to \S \ref{iterated}, we can write the blow-up $R_{s,p,n}$ as iterated blow-ups as follows. We first blow up $\overline{\mathcal V_{(p,0)}^-}$, and then the total transform of $\overline{\mathcal V_{(p-1,1)}^-}$, $\overline{\mathcal V_{(p-2,2)}^-}$, $\cdots$, $\overline{\mathcal V_{(p-r+1,r-1)}^-}$ iteratively; then, we blow up the total transform of $\overline{\mathcal V_{(p-r,r)}^+}$, $\overline{\mathcal V_{(p-r+1,r-1)}^+}$, $\cdots$, $\overline{\mathcal V_{(p-1,1)}^+}$ iteratively.  
One can show that the center of each intermediate blow-up is a scheme-theoretic union of a smooth subscheme and a divisor (which can be dropped).

Now computation in the  Mille Cr\^epes coordinates yields the desired multiplicity of $E$ along $D_{l}^{\pm}$. We conclude Lemma \ref{excep}.\,\,\,
$\endpf$
\medskip

{\noindent\bf Proof of Corollary \ref{picb}.} 
First assume that $p<s$ and $p\neq n-s$.
Let $Z\subset \mathcal T_{s,p,n}$ be an irreducible divisor.  If $Z$ is contained in the exceptional divisor of $R_{s,p,n}$, then $Z$ is a $\mathbb Z$-linear combination of $D_j^{\pm}$, $1\leq j\leq r$, by Lemma \ref{excep}. Otherwise,  considering $ R_{s,p,n}^{-1}\left(R_{s,p,n}(Z)\right)$, we can conclude that (\ref{basis1}) generates Pic$(\mathcal T_{s,p,n})$ over $\mathbb Z$. 

Next, assume the following relation in Pic$(\mathcal T_{s,p,n})$.
\begin{equation*}%\label{chow}
    0=h\cdot (R_{s,p,n})^*\left(\mathcal O_{ G(p,n)}(1)\right)+a_r^+\cdot D^+_{r}+\sum_{i=1}^{r-1}\left(a^+_i\cdot D^+_i+a^-_{r-i}D^-_{r-i}\right)+a_r^-\cdot D^-_{r}\,.
\end{equation*}
It is clear that $h=0$. Then, there is a rational function $f$ on  $\mathcal T_{s,p,n}$ such that 
\begin{equation*}
    (f)=a_r^+\cdot D^+_{r}+\sum_{i=1}^{r-1}\left(a^+_i\cdot D^+_i+a^-_{r-i}D^-_{r-i}\right)+a_r^-\cdot D^-_{r}\,,
\end{equation*}
where $(f)$ is the associated principle divisor of $f$. $f$ induces a rational function $\widetilde f$ on $G(p,n)$. Note that $\widetilde f$ is regular outside the center of the blow-up $R_{s,p,n}$ which is of codimension at lest two. Then, $\widetilde f$ extends to a regular function which must be a constant. Hence $a_i^{\pm}=0$ for $1\leq i\leq r$.

We can complete the proof for other cases in the same way.\,\,\,$\endpf$
\medskip

Let $K_{\mathcal T_{s,p,n}}$ be the canonical bundle of $\mathcal T_{s,p,n}$.
We express $K_{\mathcal T_{s,p,n}}$ in terms of the exceptional divisor and  $(R_{s,p,n})^*(\mathcal O_{G(p,n)}(1))$ as follows.
\begin{lemma}\label{kan}
When $r=p$ $(p\leq n-s)$,
    \begin{equation*}%\label{bcp}
    \begin{split}
     K_{\mathcal T_{s,p,n}}=-n&\cdot (R_{s,p,n})^*(\mathcal O_{G(p,n)}(1))+\sum_{i=1}^{r}\big((p-i+1)(n-s-i+1)-1\big)\cdot D_{i}^-\\
     & +\sum_{i=1}^{r}\big((p-i+1)(s-i+1)-1\big)\cdot D_{i}^+\,\,.
    \end{split}
    \end{equation*}
When $r=n-s$ $(n-s\leq p)$,
    \begin{equation*}%\label{bcr}
    \begin{split}
     K_{\mathcal T_{s,p,n}}=-n&\cdot (R_{s,p,n})^*(\mathcal O_{G(p,n)}(1))+\sum_{i=1}^{r}\big((p-i+1)(n-s-i+1)-1\big)\cdot D_{i}^-\\
     & +\sum_{i=1}^{r}\big((n-p-i+1)(n-s-i+1)-1\big)\cdot D_{i}^+\,\,.
    \end{split}
    \end{equation*}
\end{lemma}
{\noindent\bf Proof of Lemma \ref{kan}.} It follows from the iterated blow-up constructed in the proof of Lemma \ref{excep}.  \,\,\,$\endpf$
\medskip

For $I_k$ defined by (\ref{I_k}), computation yields that \begin{equation}\label{te}
    (R_{s,p,n})^*(P_{I_k})=\left\{ \begin{array}{cl}
    1\,\,\,\,\,\,\,\,\,\,\,\,\,\,&{\rm if}\,\,k=l\\
    \prod\nolimits_{t=k+1}^{l} a^{t-k}_{t(s-p+t)}&{\rm if}\,\,k\leq l-1\\
    \prod\nolimits_{t=l+1}^{k}b^{k+1-t}_{t(s+t)}\,\,\,\,\,&{\rm if}\,\,k\geq l+1\,\,\\
    \end{array}\right..
\end{equation}
Hence, for $0\leq k\leq r$, we have
\begin{equation}\label{bst}
\begin{split}
 B_k&=(R_{s,p,n})^*\left(\mathcal O_{G(p,n)}(1)\right)-\sum_{i=1}^{r-k}(r-k+1-i)\cdot D^+_i-\sum_{i=1}^k(k+1-i)\cdot D^-_i.
\end{split}
\end{equation}
If $p=s$ (resp. $p=n-s$) we modify (\ref{bst}) for $B_0$ (resp. $B_r$) as follows (see Remark \ref{coin}). When $p=s$, 
\begin{equation}\label{bstt0}
\begin{split}
B_0=D^+_r&=(R_{s,p,n})^*\left(\mathcal O_{G(p,n)}(1)\right)-\sum_{i=1}^{r-1}(r+1-i)\cdot D^+_i\,;
\end{split}
\end{equation}
when $p=n-s$,
\begin{equation}\label{bsttr}
\begin{split}
  B_r=D^-_r&=(R_{s,p,n})^*\left(\mathcal O_{G(p,n)}(1)\right)-\sum_{i=1}^{r-1}(r+1-i)\cdot D^-_i\,.
\end{split}
\end{equation}

To investigate the invariant divisors in a quantitative way, we define the following convenient local coordinate charts. 
\begin{definition}\label{taul}
For $0\leq l\leq r$, let $\left(A^{\tau_l},\left( J^{\tau_l}_l\right)^{-1}\right)$  be the Mille Cr\^epes coordinates of $R_{s,p,n}^{-1}(U_l)$ associated with the index  $\tau_l\in\mathbb J_l$ given by
\begin{equation*}  \left(\begin{matrix}		l+1&l+2&\cdots&r&l&l-1&\cdots&1\\		s+l+1&s+l+2&\cdots&s+r&s-p+l&s-p+l-1&\cdots&s-p+1\\	\end{matrix}\right)\,. 
\end{equation*}
We call $\left(A^{\tau_l},\left( J^{\tau_l}_l\right)^{-1}\right)$  (or $A^{\tau_l}$ by a slight abuse of notation) the $l$-th main coordinate chart.
\end{definition}

In the following, we give a geometric interpretation of the $B$-stable divisors.
For $0\leq j\leq r$, denote by  $h_j$ the positive generator of the Picard group of $\mathbb {P}^{N^j_{s,p,n}}$. Denote by $\widetilde H_j$ its pullback to $\mathbb {P}^{N_{p,n}}\times\mathbb {P}^{N^0_{s,p,n}}\times\cdots\times\mathbb {P}^{N^p_{s,p,n}}$ under the corresponding projection.  Let $H_j$ be its restriction to $\mathcal T_{s,p,n}$.
\begin{lemma}\label{partialline}For $0\leq j\leq r$,
\begin{equation*}
      H_j=R^*_{s,p,n}(\mathcal O_{G(p,n)}(1))-\sum\limits_{i=1}^{r-j}(r-i+1-j)\cdot D^+_{i}-\sum\limits_{i=1}^{j}(j+1-i)\cdot D^-_{i}\,.\,\,
    \end{equation*}
In particular, $H_j=B_j$, $0\leq j\leq r$, except when $B_0=D^+_r$ or $B_r=D^-_r$.
\end{lemma}
{\bf\noindent Proof of Lemma \ref{partialline}.} 
By computing the Pl\"ucker coordinate functions $P_{I}$ in $A^{\tau_0}$ and $A^{\tau_r}$, $I\in\mathbb I^j_{s,p,n}$, and comparing the result with (\ref{bst}),(\ref{bstt0}),(\ref{bsttr}), we conclude Lemma \ref{partialline}. \,\,\,\,\,$\endpf$

\begin{lemma}\label{basepoint}
The line bundle $\mathcal O_{\mathcal T_{s,p,n}}(B_j)$ is globally generated if
\begin{enumerate}[label=$\bullet$]
    \item $p\neq s$, $p\neq n-s$ and $0\leq j\leq r$,
    \item or $n-s=p<s$ and $0\leq j\leq p-1$,
    \item or $n-s=p=s$ and $1\leq j\leq p-1$.
\end{enumerate}
\end{lemma}
{\bf\noindent Proof of Lemma \ref{basepoint}.} Lemma \ref{basepoint} follows from  Lemma \ref{partialline}.\,\,\,\,\,$\endpf$

\begin{lemma}\label{cls}
The complete linear series of $H_j$ is isomorphic to the  complete linear series of $\mathcal O_{\mathbb {P}^{N^j_{s,p,n}}}(1)$ on $\mathbb {P}^{N^j_{s,p,n}}$.
\end{lemma}
{\bf\noindent Proof of Lemma \ref{cls}.} 
Without loss of generality, we may assume that $p\neq s, n-s$. By Lemma \ref{partialline}, we get that $H^0(\mathcal T_{s,p,n},H_j)$ consists of sections in $H^0(\mathcal T_{s,p,n},R^*_{s,p,n}(\mathcal O_{G(p,n)}(1)))\cong H^0(G(p,n),\mathcal O_{G(p,n)}(1))$ that vanish on $D^+_{1},\cdots,D^+_{r-j},D^-_{1},\cdots,D^-_{j}$ with certain prescribed orders.
Computation in the  Mille Cr\^epes coordinates yields that the all Pl\"ucker coordinate functions $P_I$ with $I\in \mathbb I_{s,p,n}^j$ (see (\ref{ikk}), (\ref{subl}), (\ref{sk})) are in $H^0(\mathcal T_{s,p,n},H_j)$. Notice that the complete linear series of $\mathcal O_{G(p,n)}(1)$ has a basis given by the Pl\"ucker coordinate functions (see \cite{Wey} for instance). Therefore, to prove Lemma \ref{cls}, it suffices to show that any section in $H^0(\mathcal T_{s,p,n},H_j)$ is generated by $P_I$, $I\in \mathbb I_{s,p,n}^j$.

Assume the contrary. Let $s\in H^0(\mathcal T_{s,p,n},H_j)$ be a section with the form
\begin{equation*}
s=P_{I^*}+\sum\nolimits_{I\in\mathbb I_{p,n}\,\,{\rm and}\,\,I\neq I^*}c_{I}\cdot P_I, \,\,c_I\in\mathbb K,
\end{equation*}
where $I^*\in  \mathbb I_{s,p,n}^{j^*}$ with $j^*\neq j$. Since $H^0(\mathcal T_{s,p,n},H_j)$ is invariant under the  ${\rm PGL}(s,\mathbb K)\times {\rm PGL}(n-s,\mathbb K)$-action, we may assume that $I^*=I_{j^*}$ (see (\ref{I_k})). Without loss of generality, we may assume that $j<j^*$. Computation as in the proof of \cite[Lemma 3.2]{FW} yields that, in the local coordinates associated to the $0$-th main coordinate chart $A^{\tau_{0}}$,
\begin{equation*}
s=(-1)^{j^*(p-j^*)}\cdot \prod_{t=1}^{j^*} (b_{t(s+t)})^{j^*+1-t}+\sum_{\alpha}F_{\alpha}\left(\overrightarrow A\right)+\sum_{\beta}G_{\beta}\left(\overrightarrow A,\widetilde X,\widetilde Y,\overrightarrow B^{1},\cdots,\overrightarrow B^{r}\right).
\end{equation*}
Here each $G_{\beta}$ is a monomial divisible by at least one of the variables in $\widetilde X,\widetilde Y,\overrightarrow B^{1},\cdots,\overrightarrow B^{r}$; each $F_{\alpha}$ is a monomial in the variables $b_{1(s+1)},\cdots,b_{p(s+p)}$, and no $F_{\alpha}$ is a constant multiple of $\prod_{t=1}^{j^*} (b_{t(s+t)})^{j^*+1-t}$. However, by the assumption that $s\in H^0(\mathcal T_{s,p,n},H_j)$, $s$ should take the form 
\begin{equation*}
s=c\cdot \prod_{t=1}^{j} (b_{t(s+t)})^{j+1-t},\,\,c\in\mathbb K. 
\end{equation*}
We get the contradiction.
 
We complete the proof of Lemma \ref{cls}.
\,\,\,\,\,$\endpf$
\medskip

We can define invariant divisors of $\mathcal M_{s,p,n}$ in a similar way.
Define 
\begin{equation*}
\check D_i:= D_{1}^-\cap D_{i}^-=\mathcal M_{s,p,n}\cap D_{i}^-,\,\,\,2\leq i\leq r.
\end{equation*} 
We call the divisors $\check D_{2}$, $\check D_{3}$, $\cdots$, $\check D_{r}$ the $G$-stable divisors of $\mathcal M_{s,p,n}$. Define 
\begin{equation*}
\check B_i:=D_{1}^-\cap B_{i}=\mathcal M_{s,p,n}\cap B_{i},\,\,\,0\leq i\leq r.
\end{equation*}
We call  $\check B_i$, $0\leq i\leq r$,  the $B$-stable divisors of  $\mathcal M_{s,p,n}$. It is easy to verify that  $\check B_j$ is an irreducible divisor if $1\leq j\leq r$, or $j=0$ and $p<s$.
\begin{remark}
$\check B_0=D_{1}^-\cap B_{0}=D_{1}^-\cap D^+_{r}=\emptyset$ is trivial when $p=s$; $\check B_r=\check D_r$ when $p=n-s$. By a slightly abuse of notation, we view $\check B_0$ as a $B$-stable divisor even if $p=s$.
\end{remark}

For convenience, we denote by $\check R_{s,p,n}$ the restriction $R_{s,p,n}\big|_{\mathcal M_{s,p,n}}:\mathcal M_{s,p,n}\rightarrow\mathcal V_{(p,0)}$. Since the hyperplane bundle of $\mathcal V_{(p,0)}\cong G(p,s)$ is the restriction of that of $G(p,n)$, we  denote  by $(\check R_{s,p,n})^*(\mathcal O_{G(p,n)}(1))$ the restriction of the line bundle $(R_{s,p,n})^*(\mathcal O_{G(p,n)}(1))$ to $\mathcal M_{s,p,n}$.  Similarly, we can prove that
%\begin{lemma}\label{mpicb}The Picard group of $\mathcal M_{s,p,n}$ is a torsion free abelian group over $\mathbb Z$. Then, Pic$(\mathcal M_{s,p,n})$ has a $\mathbb Z$-basis\begin{equation*}\left\{ \begin{array}{ll}(\check R_{s,p,n})^*(\mathcal O_{ G(p,n)}(1)),\check D_1,\check D_{2},\cdots,\check D_{r} &{\rm when}\,\,p<s,\,n-s\neq p\\(\check R_{s,p,n})^*(\mathcal O_{ G(p,n)}(1)),\check D_1,\check D_{2},\cdots,\check D_{r-1}&{\rm when}\,\,n-s=p<s\\\check D_1,\check D_{2},\cdots,\check D_{r-1} &{\rm when}\,\,n-s=p=s\,\,\\    \end{array}\right..\end{equation*}\end{lemma}

Denote by $K_{\mathcal M_{s,p,n}}$ the canonical bundle of $\mathcal M_{s,p,n}$. Then
\begin{lemma}\label{mkan}
Assume that $p\leq s$. Then,
\begin{equation*}%\label{mkan1}
\begin{split}
K_{\mathcal M_{s,p,n}}=&-n\cdot(\check R_{s,p,n})^*(\mathcal O_{G(p,n)}(1))+p(n-s)\cdot\check D_1\\
+\sum_{i=2}^{r}& \big((p-i+1)(n-s-i+1)-1\big)\cdot\check D_{i}\,.
\end{split}
\end{equation*}
\end{lemma}

{\noindent\bf Proof of  Lemma \ref{mkan}.} 
Notice that $D^+_i\cap D^-_1=\emptyset$ for $1\leq i\leq r$. 
Lemma \ref{mkan} follows form  Lemma \ref{kan} by the adjunction formula. \,\,\,$\endpf$
\medskip

Similarly to Lemma \ref{gb}, we have that
\begin{lemma}\label{checkgb}
Let $\check {\mathfrak D}$ be an irreducible divisor of $\mathcal M_{s,p,n}$. If $\check {\mathfrak D}$ is $G$-invariant (resp. $B$-invariant), then
\begin{equation*}%\label{checkgi}
\check {\mathfrak D}\in \{\check D_{2}, \check D_{3}, \cdots, \check D_{r}\}\,\,(resp.\,\,\{\check D_{2}, \check D_{3}, \cdots, \check D_{r}, \check B_0, \check B_1, \cdots, \check B_r\}).
\end{equation*}
\end{lemma}

Restricting (\ref{bst}), (\ref{bsttr}) to $\mathcal M_{s,p,n}$, we then have
\begin{equation}\label{mb=b}
\begin{split}
&\check B_i=(\check R_{s,p,n})^*(\mathcal O_{G(p,n)}(1))-i\cdot \check D_1-\sum_{k=2}^i (i+1-k)\cdot\check D_k\,,\,\,\,\,0\leq i\leq r\,;
\end{split}
\end{equation}
when $p=n-s$, we modify (\ref{mb=b}) for $\check B_r$ by 
\begin{equation*}
 \check B_r=\check D_r=(\check R_{s,p,n})^*(\mathcal O_{G(p,n)}(1))-r\cdot \check D_1-\sum_{k=2}^{r-1}(r+1-k)\cdot\check D_k\,.
\end{equation*}

\subsection{\texorpdfstring{$T$-invariant curves}{rr}} \label{curanddivi}
In this subsection, we introduce $T$-invariant curves $\gamma_l$, $\zeta^l_j$, $\zeta_{u,v}^{l,k}$, $\delta_{m_1,m_2}^l$, $\Delta_{m_1,m_2}^l$ of $\mathcal T_{s,p,n}$ with various parameters as the closures of certain affine lines in the $l$-th main coordinate chart $A^{\tau_l}$, $0\leq l\leq r$.
Then, we compute their intersection numbers with line bundles $R^*_{s,p,n}((\mathcal O_{G(p,n)}(1))$, $D^{\pm}_1,\cdots,D^{\pm}_{r},-K_{\mathcal T_{s,p,n}}$. Such computation will be used to show the (semi-)positivity of the anti-canonical bundles.

For convenience, we write $R^*_{s,p,n}((\mathcal O_{G(p,n)}(1))$ as $H$, and let  $\delta_{ab}$ be the Kronecker delta.
\smallskip

{\bf\noindent (0). Curves $\gamma_l$, $0\leq l\leq r-1$.}
\smallskip

In terms of the local coordinates, we define an affine line $\mathring{\gamma}_l:\mathbb K\rightarrow A^{\tau_l}$, $0\leq l\leq r-1$,  as follows. $b_{(l+1)(s+l+1)}\left(\mathring{\gamma}_l(t)\right)=t$, and all other variables are constantly zero. %$\widetilde X\left(\mathring{\gamma}_l(t)\right)$ and $\widetilde Y\left(\mathring{\gamma}_l(t)\right)$ are zero matrices for each $t\in\mathbb K$; $\overrightarrow B^{k}\left(\mathring{\gamma}_l(t)\right)$ is a zero vector for each $t\in\mathbb K$ and $1\leq k\leq r$\,; \begin{equation*} \begin{split}\overrightarrow B^{1}\left(\mathring{\gamma}_l(t)\right)&=\left(b_{(l+1)(s+l+1)},\xi^{(1)}_{(l+1)(s+l+2)},\xi^{(1)}_{(l+1)(s+l+3)},\cdots,\xi^{(1)}_{(l+1)n},\right.\\ &\,\,\,\left.\xi^{(1)}_{(l+2)(s+l+1)},\xi^{(1)}_{(l+3)(s+l+1)},\cdots,\xi^{(1)}_{p(s+l+1)}\right)\left(\mathring{\gamma}_l(t)\right)=\left(t,0,0\cdots,0\right)\,.\end{split}\end{equation*}
Define $\gamma_l$ to be the closure of ${\mathring{\gamma}_l}$ in $\mathcal T_{s,p,n}$.

\begin{lemma}\label{i1} For $0\leq l\leq r-1$ and $1\leq i\leq r$, the intersection numbers of the divisors
$H,D_i^\pm$ with the curves $\gamma_l$ are given by
\[
\begin{array}{llll}
H\cdot \gamma_l = 1,&
D_i^-\cdot \gamma_l = \delta_{i,l+1}-\delta_{i,l+2},&
D_i^+\cdot \gamma_l = \delta_{i,r-l}-\delta_{i,r-l+1}.
\end{array}
\]
\end{lemma}
{\bf\noindent Proof of Lemma \ref{i1}.} We view ${\gamma}_l$ as a curve in $\mathbb {P}^{N_{p,n}} \times\mathbb{P}^{N^0_{s,p,n}} \times\cdots\times\mathbb {P}^{N^p_{s,p,n}}$. 
By computing the Pl\"ucker coordinate functions, it is easy to verify that the projection of ${\gamma}_l$  to $\mathbb {P}^{N_{p,n}}$  is a line, and the projection of ${\gamma}_l$  to $\mathbb {P}^{N^j_{s,p,n}}$ is a point for $0\leq j\leq p$.
Therefore,  we have that ${\gamma}_l\cdot H=1$ and ${\gamma}_l\cdot H_j=0$, $0\leq j\leq p$. Notice that $\gamma_l$ is disjoint with the divisors  $D_1^-,D^-_2,\cdots,D^-_l,D^+_1,D^+_2,\cdots,D^+_{r-l-1}$. Combined with Lemma \ref{partialline}, we conclude Lemma \ref{i1} by solving linear equations.\,\,\,\,$\endpf$
\medskip

Combining with Lemma \ref{kan}, a direct computation yields that
\begin{lemma}\label{vki1} When $r=1$,
$-K_{\mathcal T_{s,p,n}}\cdot\gamma_0=2$.
When $r\geq 2$,
\begin{equation*}
-K_{\mathcal T_{s,p,n}}\cdot\gamma_l=\delta_{0,l}+\delta_{r-1,l}.
\end{equation*}
\end{lemma}

{\bf\noindent (1). Curves $\zeta^0_j$, $\zeta_{u,v}^{0,k}$, $\delta_{m_1,m_2}^0$ in $R_{s,p,n}^{-1}\left(\mathcal V_{(p,0)}\right)$.}
\smallskip

For $2\leq j\leq r$, define an affine line $\mathring{\zeta}^0_j\subset A^{\tau_0}$ such that $b_{j(s+j)}\left(\mathring{\zeta}^0_j(t)\right)=t$ and the other variables are constantly zero. Let  $\zeta^0_j$ be the closure of $\mathring{\zeta}^0_j$.
Similarly, we can derive that
\begin{lemma}\label{i2} 
For $2\leq j\leq r$ and $1\leq i\leq r$,
\[
\begin{split}
H\cdot \zeta^0_j = 0,&\,\,\,\,\,
D_i^-\cdot \zeta^0_j= -\delta_{i,j-1}+2\delta_{i,j}-\delta_{i,j+1},\,\,\,\,\,
D_i^+\cdot \zeta^0_j = 0,\\
&-K_{\mathcal T_{s,p,n}}\cdot\zeta^0_j=3-\delta_{r,j}.
\end{split}
\]
\end{lemma}

%{\textcolor{red}{$\zeta_{u,v}^{0,k}$ are useless}}

Assume that $u=k$ and $s+k+1\leq v\leq n$, or $v=s+k$  and $k+1\leq u\leq p$ where  $1\leq k\leq r$.  Define an affine line $\mathring{\zeta}_{u,v}^{0,k}\subset A^{\tau_0}$ such that $\xi^{(k)}_{uv}\left(\mathring{\zeta}_{u,v}^{0,k}(t)\right)=t$ and the other variables are zero.
Let $\zeta_{u,v}^{0,k}$ be its closure. The intersection numbers are given by Lemmas \ref{li1} and \ref{vkli1} with $l=0$.

When $1\leq m_1\leq p$ and $1\leq m_2\leq s-p$, define an affine line $\mathring{\delta}_{m_1,m_2}^0\subset  A^{\tau_0}$ such that $y_{m_1(s-p+1-m_2)}\left(\mathring{\delta}_{m_1,m_2}^0(t)\right)=t$ and the other variables are zero. Let $\delta_{m_1,m_2}^0$  be its closure. The intersection numbers are given by Lemmas \ref{li4} and \ref{vkli4} with $l=0$. 
\smallskip

{\bf\noindent (2). Curves  $\zeta^l_j$, $\zeta_{u,v}^{l,k}$, $\delta_{m_1,m_2}^l$, $\Delta_{m_1,m_2}^l$  in $R_{s,p,n}^{-1}\left(\mathcal V_{(p-l,l)}\right)$, $1\leq l\leq r-1$.}
\smallskip

Assume that $1\leq l\leq r-1$. For $2\leq j\leq r-l$, define an affine line $\mathring{\zeta}^l_j\subset A^{\tau_l}$ such that
$b_{(l+j)(s+l+j)}\left(\mathring{\zeta}^l_j(t)\right)$ and the other variables are constantly zero.
For $r-l+2\leq j\leq r$, define an affine line $\mathring{\zeta}^l_j\subset A^{\tau_l}$ such that
$a_{(r+1-j)(s-p+r+1-j)}\left(\mathring{\zeta}^l_j(t)\right)$ and the other variables are zero. Let $\zeta^l_j$ be the closure of  $\mathring{\zeta}^l_j$ for $2\leq j\leq r-l$ or $r-l+2\leq j\leq r$. The intersection numbers are given by Lemma \ref{i5}.

Assume that one of the following holds for integers $k,u,v$.
\begin{equation*}%\label{l1kuv}
\begin{split}
&(a)\,\,\,\,1\leq k\leq r-l\,; \,u=l+k\,\,{\rm and\,\,}s+l+k+1\leq v\leq n,\,\,{\rm or\,\,}v=s+l+k\,\,{\rm and}\\
&\,\,\,\,\,\,\,\,\,\,\,\,\,\,\,\,\,\,    l+k+1\leq u\leq p.\\
&(b)\,\,\,\,r-l+1\leq k\leq r\,; u=r-k+1\,\,{\rm and\,\,}1\leq v\leq s-p+r-k,\,\,{\rm or\,\,}\\
&\,\,\,\,\,\,\,\,\,\,\,\,\,\,\,\,\,\,  v=s-p+r-k+1\,\,{\rm and\,\,}1\leq u\leq r-k.
    \end{split}
\end{equation*}
We define an affine line  $\mathring{\zeta}_{u,v}^{l,k}\subset A^{\tau_l}$ such that $\xi^{(k)}_{uv}\left(\mathring{\zeta}_{u,v}^{l,k}(t)\right)=t$ and the other variables are zero. Let $\zeta_{u,v}^{l,k}$ be its closure. The intersection numbers are given by Lemmas \ref{li1}, \ref{vkli1}, \ref{li2}, \ref{vkli2}.

For $1\leq m_1\leq p-l$ and $1\leq m_2\leq s-p+l$,  define an affine line 
$\mathring{\delta}_{m_1,m_2}^l\subset A^{\tau_l}$  such that $y_{(l+m_1)(s-p+l+1-m_2)}\left(\mathring{\delta}_{m_1,m_2}^l(t)\right)=t$ and the other variables are zero.  Let $\delta_{m_1,m_2}^l$ be its closure. The intersection numbers are given by Lemmas \ref{li4} and \ref{vkli4}.

For $1\leq m_1\leq n-s-l$ and $1\leq m_2\leq l$,  define  an affine line $\mathring{\Delta}_{m_1,m_2}^l\subset A^{\tau_l}$ such that $x_{(l+1-m_2)(s+l+m_1)}\left(\mathring{\Delta}_{m_1,m_2}^l(t)\right)=t$ and the other variables are zero. Let  $\Delta_{m_1,m_2}^l$ be its closure. The intersection numbers are given by Lemma \ref{li1} with $l=r$.  
\smallskip

{\bf\noindent (3). Curves  $\zeta^r_j$, $\zeta_{u,v}^{r,k}$, $\delta_{m_1,m_2}^r$, $\Delta_{m_1,m_2}^r$ in $R_{s,p,n}^{-1}\left(\mathcal V_{(p-r,r)}\right)$.}
\smallskip

For $2\leq j\leq r$, define an affine line  $\mathring{\zeta}^r_j\subset A^{\tau_r}$  such that  $a_{(r+1-j)(s-p+r+1-j)}\left(\mathring{\zeta}^r_j(t)\right)=t$ and the other variables are zero. Let $\zeta^r_j$ be its closure. The intersection numbers are given by Lemma \ref{i5}.

Let $1\leq k\leq r$. Assume that $u=r-k+1$ and $1\leq v\leq s-p+r-k$, or $v=s-p+r-k+1$ and $1\leq u\leq r-k$. Define an affine line $\mathring{\zeta}_{u,v}^{r,k}\subset  A^{\tau_r}$ such that $\xi^{(k)}_{uv}\left(\mathring{\zeta}_{u,v}^{r,k}(t)\right)=t$ and other variables are zero. Let $\zeta_{u,v}^{r,k}$ be its closure. The intersection numbers are given by Lemmas \ref{li2} and \ref{vkli2} with $l=r$.

When $n-s<p$,   $1\leq m_1\leq s+p-n$, and $1\leq m_2\leq n-p$,  define an affine line $\mathring{\delta}_{m_1,m_2}^r\subset A^{\tau_{r}}$ such that $y_{(r+m_1)(s-p+r+1-m_2)}\left(\mathring{\delta}_{m_1,m_2}^r(t)\right)=t$ and other variables are zero.
Let $\delta_{m_1,m_2}^r$ be its closure. The intersection numbers are given by Lemmas \ref{li4} and \ref{vkli4} with $l=r$.

When  $p<n-s$, $1\leq m_1\leq n-s-r$, and $1\leq m_2\leq r$,    define an affine line $\mathring{\Delta}_{m_1,m_2}^r\subset A^{\tau_l}$ such that $x_{(r+1-m_2)(s+r+m_1)}\left(\mathring{\Delta}_{m_1,m_2}^r(t)\right)=t$ and the other variables are zero. Let $\Delta_{m_1,m_2}^r$ be its closure. The intersection numbers are given by Lemmas \ref{li5} and \ref{vkli5} with $l=r$.

\subsection{(Semi-)positivity of the anti-canonical bundles} \label{curanddivp}
In this subsection, we establish the (semi-)positivity of the anti-canonical bundles of $\mathcal T_{s,p,n}$  and $\mathcal M_{s,p,n}$. 

%Following \cite{Ful}, we first recall the notion of Chow groups. %An algebraic cycle on $X$ means a finite linear combination of closed subvarieties of $X$ with integer coefficients. For a natural number $i$, the group $Z_{i}(X)$ of $i$-cycles on $X$ is the free abelian group on the set of $i$-dimensional subvarieties of $X$. For a variety $W$ of dimension $i+1$ and any rational function $f$ on $W$ which is not identically zero, the divisor of $f$ is the $i$-cycle $(f)=\sum _{Z}\operatorname {ord} _{Z}(f)Z$, where the sum runs over all $i$-dimensional subvarieties $Z$ of $W$ and the integer $ \operatorname {ord} _{Z}(f)$ denotes the order of vanishing of $f$ along $Z$. The group of $i$-cycles rationally equivalent to zero is the subgroup of $ Z_{i}(X)$ generated by the cycles $(f)$ for all $(i+1)$-dimensional subvarieties $W$ of $X$ and all nonzero rational functions $f$ on $W$. 
Denote by $A_{i}(X)$ the Chow group of $i$-cycles on $X$ (see \cite{Ful} for details). %is the quotient group of $Z_{i}(X)$ by the subgroup of cycles rationally equivalent to zero. 
Brion characterized the Chow groups of spherical varieties  as follows.
\begin{lemma}[\cite{Br3}]\label{cone}
Let $X$ be an irreducible, complete spherical variety of dimension $n$.
The cone of effective cycles in $A_{i}(X)\otimes_{\mathbb Z}\mathbb Q$ is a polyhedral convex cone generated by the classes of the closures of the $B$-orbits.
\end{lemma}

{\bf\noindent Proof of Proposition \ref{fano}.} According to \cite[Theorem 1.3]{FZ13}, $-K_{\mathcal T_{s,p,n}}$ is big. By Kleiman's criterion and Lemma \ref{cone}, to prove Proposition \ref{fano}, it suffices to show that $-K_{\mathcal T_{s,p,n}}$ has non-negative degrees on all the irreducible curves  of $\mathcal T_{s,p,n}$ which are the closures of the $B$-orbits of dimension $1$, and that it has strictly positive degrees on such curves if and only if $r\leq 2$.

%Without loss of generality, we may assume that $2p\leq n\leq 2s$.

Notice that for any curve $\gamma$ of $\mathcal T_{s,p,n}$ that is the closure of a $1$-dimensional $B$-orbit, there is a permutation matrix $g\in{\rm GL}(s,\mathbb K)\times {\rm GL}(n-s,\mathbb K)$ and an integer $0\leq l^*\leq r$ such that $A^{\tau_{l^*}}$ contains an open subset of the image curve $g(\gamma)$ (see Definition \ref{taul} for  $A^{\tau_l}$).   
Let $T$ be the maximal torus of $B$. It is clear that the image curve $g(\gamma)$ is invariant under the $T$-action. Consider the local coordinates of a generic point $\mathfrak a$ of $g(\gamma)$ with respect to the $l^*$-th main coordinate chart   $A^{\tau_l^*}$. Then only one of the local coordinates is nonzero, for otherwise $g(\gamma)$ is of dimension at least two. Therefore, $g(\gamma)$ coincides with one of the $T$-invariant curves $\gamma_l,\zeta_j^l,\zeta_{u,v}^{l,k},\delta_{m_1,m_2}^l,  \Delta_{m_1,m_2}^l$ defined in Section \ref{curanddivi}. 

Since $K_{\mathcal T_{s,p,n}}$ is invariant under the holomorphic automorphism of $\mathcal T_{s,p,n}$, we have $-K_{\mathcal T_{s,p,n}}\cdot\gamma=-K_{\mathcal T_{s,p,n}}\cdot g(\gamma)$. By the intersection numbers computed in \S \ref{curanddivi}, we conclude that $-K_{\mathcal T_{s,p,n}}\cdot\gamma\geq 0$, and that $-K_{\mathcal T_{s,p,n}}\cdot\gamma=0$ if and only if $g(\gamma)$ coincides with $\gamma_l$ for a certain $1\leq l\leq r-2$.

We complete the proof of Proposition \ref{fano}.\,\,\,\,\,$\endpf$

\medskip

{\bf\noindent Proof of Proposition \ref{mfano}.} The proof is the same as that in Proposition \ref{fano}. For any curve $\check\gamma\subset\mathcal M_{s,p,n}$ that is the closure of a $1$-dimensional $B$-orbit, there is a permutation matrix $g\in {\rm GL}(s,\mathbb K)\times {\rm GL}(n-s,\mathbb K)$ such that $A^{\tau_0}$ contains an open subset of the image of $g(\check\gamma)$.  Moreover, $g(\check\gamma)$ coincides with one of the $T$-invariant curves $\zeta_j^0,\zeta_{u,v}^{0,k},\delta_{m_1,m_2}^0$ defined in Section \ref{curanddivi}.
 
Notice that for each curve $\gamma\subset\mathcal M_{s,p,n}$,
\begin{equation*}
    -K_{\mathcal M_{s,p,n}}\cdot \gamma=(-K_{\mathcal T_{s,p,n}}-D_1^-)\cdot\gamma
\end{equation*}
Proposition \ref{mfano} follows from the calculation in \S \ref{curanddivi}. \,\,\,\,\,$\endpf$
\medskip

Similarly, we can prove by Lemmas \ref{partialline}, \ref{basepoint} that
\begin{lemma}\label{mcls}
For $0\leq j\leq r$, the line bundle $\mathcal O_{\mathcal T_{s,p,n}}(-K_{\mathcal T_{s,p,n}}-D_1^-+H_j)$ is big and numerical
effective. 
\end{lemma}

\section{Symmetries of \texorpdfstring{$\mathcal T_{s,p,n}$ }{jj} and \texorpdfstring{$\mathcal M_{s,p,n}$ }{jj} } \label{sym}
%In this section, we will determine the automorphisms of $\mathcal T_{s,p,n}$ and $\mathcal M_{s,p,n}$. The proof consists of two parts. Firstly, we show that the induced automorphisms of the Picard group consist of at most two elements. Secondly, we complete the proof under the extra assumption that the induced automorphism of the Picard group is the identity map.  

We first compute the isomorphisms of the Picard groups induced by the special  discrete automorphisms.
\begin{lemma}\label{dis}
Assume that  $1\leq p\leq s$.  The  automorphism  {\rm USD} induces an automorphism of the Picard group of $\mathcal T_{s,p,2s}$ as follows.
\begin{equation*}
\begin{split}
&({\rm USD})^*(D^+_i)=D^-_i,\,\,\, ({\rm USD})^*(D^-_i)=D^+_i,\,\,\, 1\leq i\leq r\,;\\
&({\rm USD})^*\left((R_{s,p,2s})^*(\mathcal O_{G(p,2s)}(1))\right)=(R_{s,p,2s})^*(\mathcal O_{G(p,2s)}(1))\,;\\
&({\rm USD})^*\left(B_i\right)=B_{r-i}\,,\,\,0\leq i\leq r.
\end{split}
\end{equation*}
\end{lemma}
{\noindent\bf Proof of Lemma \ref{dis}.} Since {\rm USD} induces an automorphism of $G(p,2s)$, we have $({\rm USD})^*\left((R_{s,p,2s})^*(\mathcal O_{G(p,2s)}(1))\right)=(R_{s,p,2s})^*(\mathcal O_{G(p,2s)}(1))$. By (\ref{dtrans}), we conclude that {\rm USD} interchanges $D^+_i$ and $D^-_i$ for $1\leq i\leq r$\,. \,\,\,$\endpf$

\begin{lemma}\label{dis2}Assume that  $1\leq p\leq s$. The {\rm DUAL} automorphism  induces an automorphism of the Picard group of 
$\mathcal T_{s,p,2p}$ as follows.
\begin{equation*}
\begin{split}
&({\rm DUAL})^*(D^+_i)=D^-_i,\,\,\, ({\rm DUAL})^*(D^-_i)=D^+_i,\,\,\, 1\leq i\leq r\,;\\
&({\rm DUAL})^*\left((R_{s,p,2p})^*(\mathcal O_{G(p,2p)}(1))\right)=(R_{s,p,2p})^*(\mathcal O_{G(p,2p)}(1))\,;\\
&({\rm DUAL})^*\left(B_i\right)=B_{r-i}\,,\,\,0\leq i\leq r.
\end{split}
\end{equation*}
\end{lemma}
{\noindent\bf Proof of Lemma \ref{dis2}.} The proof is the same as that in Lemma \ref{dis}. \,\,\,$\endpf$
\medskip

Similarly, we can prove that
\begin{lemma}\label{mdis}Assume that  $1\leq p<s$. Then the  automorphism  {\rm Usd} induces an automorphism of the Picard group of 
$\mathcal M_{s,p,2s}$ as follows.
\begin{equation*}
\begin{split}
&({\rm Usd})^*(\check D_i)=\check D_{r+2-i},\,\,\, 2\leq i\leq r\,;\,\,({\rm Usd})^*(\check D_1)=-\sum_{i=1}^r\check D_i\,;\\
&({\rm Usd})^*\left((\check R_{s,p,2s})^*(\mathcal O_{G(p,2s)}(1))\right)=(\check R_{s,p,2s})^*(\mathcal O_{G(p,2s)}(1))-\sum_{i=1}^r(r+1-i)\cdot\check D_i\,;\\
&({\rm Usd})^*\left(\check B_i\right)=\check B_{r-i}\,,\,\,0\leq i\leq r.
\end{split}
\end{equation*}
\end{lemma}

\begin{lemma}\label{mdis2}Assume that  $1\leq p<s$. Then the  automorphism  {\rm Dual} induces an automorphism of the Picard group of 
$\mathcal M_{s,p,2p}$ as follows.
\begin{equation*}
\begin{split}
&({\rm Dual})^*(\check D_i)=\check D_{r+2-i},\,\,\, 2\leq i\leq r\,;\,\,({\rm Dual})^*(\check D_1)=-\sum_{i=1}^r\check D_i\,;\\
&({\rm Dual})^*\left((\check R_{s,p,2p})^*(\mathcal O_{G(p,2p)}(1))\right)=(\check R_{s,p,2p})^*(\mathcal O_{G(p,2p)}(1))-\sum_{i=1}^r(r+1-i)\cdot\check D_i\,;\\
&({\rm Dual})^*\left(\check B_i\right)=\check B_{r-i}\,,\,\,0\leq i\leq r.
\end{split}
\end{equation*}
\end{lemma}

%Recall that when $n=2p=2s$, $(\check R_{p,p,2p})^*(\mathcal O_{G(p,2p)}(1))$ and hence $\sum_{i=1}^r(r+1-i)\cdot\check D_i$ are trivial. We can derive that
\begin{lemma}\label{mdis3}Assume that  $p\geq1$.  The  automorphisms  {\rm Usd} and {\rm Dual} induce the following automorphism of the Picard group of $\mathcal M_{p,p,2p}$.
\begin{equation*}
\begin{split}
&({\rm Usd})^*(\check D_i)=({\rm Dual})^*(\check D_i)=\check D_{r+2-i},\,\,\, 2\leq i\leq r\,;\\
&({\rm Usd})^*(\check D_1)=({\rm Dual})^*(\check D_1)=-\sum_{i=2}^r\frac{i-1}{r}\cdot\check D_i\,;\\
&({\rm Usd})^*\left(\check B_i\right)=({\rm Dual})^*\left(\check B_i\right)=\check B_{r-i}\,,\,\,1\leq i\leq r-1.
\end{split}
\end{equation*}
\end{lemma}

\subsection{The automorphism groups of \texorpdfstring{$\mathcal T_{s,p,n}$}{ff}}  \label{symt} 

\begin{lemma}\label{des}
Let $\sigma\in{\rm Aut}(\mathcal T_{s,p,n})$. Denote by $\Sigma$ the birational self-map of $G(p,n)$ induced by $\sigma$. If the induced automorphism $\sigma^*$ of the Picard group of $\mathcal T_{s,p,n}$ is the identity map, then  $\Sigma$ extends to an automorphism of $G(p,n)$.
\end{lemma}
{\noindent\bf Proof of Lemma \ref{des}.} For convenience, we use the following notation to distinguish the source and the target.
\begin{equation*}%\label{idesent}
\begin{array}{ccc}
\vspace{.03in}
\mathcal T^1_{s,p,n}&\xrightarrow{\,\,\,\,\,\,\sigma\,\,\,\,\,} &\mathcal T^2_{s,p,n}\\
\vspace{-.04in}
\Big\downarrow\llap{$\scriptstyle R^1_{s,p,n}\,\,\,\,\,$}&&\Big\downarrow\rlap{$\scriptstyle R^2_{s,p,n}\,\,\,\,\,$}\\
G^1(p,n)&\xdashrightarrow{\,\,\,\Sigma\,\,\,}&G^2(p,n)\\
\end{array}\,\,\,\,\,\,\,\,\,\,\,\,\,\,\,\,.
\end{equation*}

We next prove that $\Sigma$ extends to a regular map 
from $G^1(p,n)$ to $G^2(p,n)$.  

Take any $x\in G^1(p,n)$. If $(R_{s,p,n}^2\circ\sigma)((R_{s,p,n}^1)^{-1}(x))$ consists of more than one point, we can find a curve $E\subset (R_{s,p,n}^1)^{-1}(x)$ such that  $(R_{s,p,n}^2\circ\sigma)(\sigma(E))$ is a curve, for $\mathcal T_{s,p,n}$ can be realized as a sequence of blow-ups of $G^1(p,n)$ along smooth submanifolds. However, any curve $E\subset\mathcal  T_{s,p,n}^1$ that is contracted by $R_{s,p,n}^1$ must also be contracted by $R_{s,p,n}^2$, as otherwise we would have 
\begin{equation*}
\begin{split}
     0&=\mathcal O_{G^1(p,n)}(1)\cdot R_{s,p,n}^1(E)=(R_{s,p,n}^1)^*(\mathcal O_{G^1(p,n)}(1))\cdot E\\
     &=(\sigma^{-1})^*\left((R_{s,p,n}^1)^*\left(\mathcal O_{G^1(p,n)}(1) \right)\right)\cdot \sigma (E) =(R_{s,p,n}^2)^*(\mathcal O_{G^2(p,n)}(1))\cdot \sigma(E) \\
     &=\mathcal O_{G^2(p,n)}(1)\cdot R_{s,p,n}^2\left(\sigma(E)\right)\geq 1\,.
\end{split}
\end{equation*}

Now $(R_{s,p,n}^2\circ\sigma)((R_{s,p,n}^1)^{-1}(x))$ is a single point. Hence, there is a neighborhood $U$ of $x\in G^1(p,n)$ such that $\Sigma(U)$ is contained in a small neighborhood of $(R_{s,p,n}^2\circ\sigma)((R_{s,p,n}^1)^{-1}(x))$. In local coordinates, we can write $\Sigma$ as a bounded vector-valued rational function which is bounded near $x$. We thus conclude $\Sigma$ extend to a regular map near $x$.

Similarly, $\Sigma^{-1}$ has a regular extension. In conclusion, we conclude that ${\Sigma}$ extends to an automorphism of $G(p,n)$.
\,\,\,$\endpf$
\medskip

The induced automorphisms of the Picard groups are classified by
\begin{lemma}\label{sigma2}
Let $\sigma$ be an automorphism of $\mathcal T_{s,p,n}$ and $\sigma^*$ the induced automorphism of the Picard group. When $n\neq 2s$ and $n\neq 2p$, $\sigma^*$ is the identity map. When $n=2s$, $\sigma^*$ is the identity map or $({\rm USD})^*$. When  $n=2p$, $\sigma^*$ is the identity map or $({\rm DUAL})^*$.
\end{lemma}
{\bf\noindent Proof of Lemma \ref{sigma2}.} Denote by  Eff$(\mathcal T_{s,p,n})\subset A_{n-1}(\mathcal T_{s,p,n})\otimes_{\mathbb Z}\mathbb Q$ the cone  of effective divisors of $\mathcal T_{s,p,n}$. By Lemmas \ref{gb}, \ref{cone}, Eff$(\mathcal T_{s,p,n})$ is generated by $\{B_0,\cdots,B_r,D^+_1,\cdots,D^+_{r},D^-_1,\cdots,D^-_{r} \}$. By (\ref{bst})-(\ref{bsttr}),  one can conclude that the set of extremal rays of Eff$(\mathcal T_{s,p,n})$ is given by
\begin{equation*}
\left\{ \begin{array}{ll}
\{B_0,B_r,D^+_1, \cdots,D^+_{r},D^-_1,\cdots,D^-_{r} \} &{\rm when}\,\,p<s,\,n-s\neq p\\
\{B_0,B_r,D^+_1,\cdots, D^+_{r},D^-_1,\cdots,D^-_{r-1} \}&{\rm when}\,\,n-s=p<s\\
\{B_0,B_r,D^+_1,\cdots, D^+_{r-1},D^-_1,\cdots,D^-_{r-1} \}&{\rm when}\,\,n-s=p=s\,\,\\
    \end{array}\right..
\end{equation*}
Since $\sigma^*$ preserves Eff$(\mathcal T_{s,p,n})$, it induces a permutation of $\mathfrak G$. 

We prove Lemma \ref{sigma2} based on a case by case argument. For simplicity of notation, we denote $(R_{s,p,n})^*\left(\mathcal O_{G(p,n)}(1)\right)$ by $H$ in the following.

\smallskip

{\bf\noindent Case 1 ($r\geq 2$, $p\neq n-s$ and $p<s$).}
By (\ref{bst}) we have that
\begin{equation}\label{b=b}
    \sigma^*(B_0)=\sigma^*(B_r)+\sum_{i=1}^r(r+1-i)\cdot \sigma^*(D^-_i)-\sum_{i=1}^{r}(r+1-i)\cdot \sigma^*(D^+_i)\,.
\end{equation}
%By Lemma \ref{picb}, it is clear that (\ref{b=b}) has integer coefficients with respect to the basis of the Picard group $\left\{H,D^-_1, \cdots, D^-_r,D^+_1,\cdots,D^+_r \right\}$.
If $\sigma^*(B_0)\in\left\{D^-_1, \cdots, D^-_r,D^+_1,\cdots,D^+_r\right\}$, by checking the coefficient of $H$, we have the following possibilities.
\begin{enumerate}[label=(\Alph*)]
    \item For a certain $1\leq j\leq r$, $\sigma^*(D_j^-)=B_0$ and $\sigma^*(D_j^+)=B_r$.
    \item For a certain $1\leq j\leq r$, $\sigma^*(D_j^+)=B_0$ and $\sigma^*(D_j^-)=B_r$.
    \item $\sigma^*(B_r)=B_r$ and $\sigma^*(D_r^+)=B_0$.
    \item $\sigma^*(B_r)=B_0$ and $\sigma^*(D_r^+)=B_r$.
\end{enumerate}

For Case (A), the coefficient of $D^+_1$ is at most $-r(r-j)$ in 
\begin{equation*}
  \sum\nolimits_{i=1}^r(r+1-i)\cdot \sigma^*(D^-_i)-\sum\nolimits_{i=1}^{r}(r+1-i)\cdot \sigma^*(D^+_i)\,.
\end{equation*}
Since $r\geq 2$, we have  $j=r$ and
\begin{equation*}
    \begin{split}
        &\sigma^*(D_r^-)=B_0\,,\,\,\,\sigma^*(D_r^+)=B_r\,,\,\,\,\sigma^*(B_0)=D_r^-\,,\,\,\,\sigma^*(B_r)=D_r^+\,,\\
        &\sigma^*(D_j^+)=D_j^-\,,\,\,\,\sigma^*(D_j^-)=D_j^+\,\,{\rm for\,\,}1\leq j\leq r-1\,.
    \end{split}
\end{equation*}
Since $\sigma^*(K_{\mathcal T_{s,p,n}})=K_{\mathcal T_{s,p,n}}$, we conclude that $n=2s$, and hence $\sigma^*=({\rm USD})^*$.

For Case (B), in the same way as above, the only possibility is that \begin{equation*}
    \begin{split}
        &\sigma^*(D_r^+)=B_0\,,\,\,\,\sigma^*(D_r^-)=B_r\,,\,\,\,\sigma^*(B_0)=D_r^+\,,\,\,\,\sigma^*(B_r)=D_r^-\,,\\
        &\sigma^*(D_j^-)=D_j^-\,,\,\,\,\sigma^*(D_j^+)=D_j^+\,\,{\rm for\,\,}1\leq j\leq r-1\,.
    \end{split}
\end{equation*}
Note that the coefficient of $H$ in $K_{\mathcal T_{s,p,n}}$ is $-n$. When $n-s<p$,  the coefficient of $H$ in $\sigma^*(K_{\mathcal T_{s,p,n}})$ is $-2r=-2(n-s)>-2p>-n$. When $p<n-s$, the coefficient of $H$ in $\sigma^*(K_{\mathcal T_{s,p,n}})$ is $-2p>-s-(n-s)=-n$. Both  lead to a contradiction. By the same argument, we can exclude Cases (C), (D).

Next, we assume that  $\sigma^*(B_0)\in\left\{B_0, B_r\right\}$. By checking the coefficient of $H$ in (\ref{b=b}), we have the following possibilities. \begin{enumerate}[label=(\alph*)]
    \item  $\sigma^*(B_0)=B_r$ and $\sigma^*(D_r^-)=B_0$.
    \item $\sigma^*(B_0)=B_r$ and $\sigma^*(B_r)=B_0$.
    \item $\sigma^*(B_0)=B_0$ and $\sigma^*(D^-_r)=B_r$.
    \item $\sigma^*(B_0)=B_0$ and $\sigma^*(B_r)=B_r$.
\end{enumerate}

For Case (a), by (\ref{b=b}) we can conclude that \begin{equation*}
    \begin{split}
        &\sigma^*(B_0)=B_r\,,\,\,\,\sigma^*(D_r^-)=B_0\,,\,\,\,\sigma^*(D_r^+)=D_r^-\,,\,\,\,\sigma^*(B_r)=D_r^+\,,\\
        &\sigma^*(D_j^+)=D_j^-\,,\,\,\,\sigma^*(D_j^-)=D_j^+\,\,{\rm for\,\,}1\leq j\leq r-1\,.
    \end{split}
\end{equation*}
Then $\sigma^*(H)=H$, and hence $\sigma^*(K_{\mathcal T_{s,p,n}})\neq K_{\mathcal T_{s,p,n}}$. This is a contradiction.

For Case (b), by (\ref{b=b}) we can conclude that \begin{equation*}
    \begin{split}
        &\sigma^*(B_r)=B_0\,,\,\,\,\sigma^*(B_0)=B_r\,,\,\,\sigma^*(D_j^+)=D_j^-\,,\,\,\,\sigma^*(D_j^-)=D_j^+\,\,{\rm for\,\,}1\leq j\leq r\,.
    \end{split}
\end{equation*}
Since $\sigma^*(K_{\mathcal T_{s,p,n}})=K_{\mathcal T_{s,p,n}}$, we conclude that either $n=2p$ and $\sigma^*=({\rm DUAL})^*$, or  $n=2s$ and $\sigma^*=({\rm USD})^*$.

For Case (c),  by (\ref{b=b}) we can show that $\sigma^*$ permutes $D_r^-$ and $B_r$, and fixes all the other divisors. Since $\sigma^*(K_{\mathcal T_{s,p,n}})=K_{\mathcal T_{s,p,n}}$, we have that $p=n-s$ which contradicts the assumption.

For Case (d),  by (\ref{b=b}) we can show that $\sigma^*$ is the identity map.

\smallskip

{\bf\noindent Case 2 ($r\geq 2$ and $n-s=p<s$).} By (\ref{bst}) and (\ref{bsttr}) we have that
\begin{equation}\label{b=b2}
    \sigma^*(B_0)=\sigma^*(B_r)+\sum_{i=1}^{r-1}(r+1-i)\cdot \sigma^*(D^-_i)-\sum_{i=1}^{r}(r+1-i)\cdot \sigma^*(D^+_i)\,.
\end{equation}
By the same argument as above, we can conclude that  $\sigma^*(B_0)\in\left\{B_0, B_r\right\}$. Checking the coefficient of $H$ in (\ref{b=b2}), we have that  $\sigma^*$ is the identity map.

\smallskip

{\bf\noindent Case 3 ($r\geq 2$ and $n-s=p=s$).} By (\ref{bstt0}) and (\ref{bsttr}) we have that
\begin{equation*}%\label{b=b4}
\sigma^*( B_0)=\sigma^*(B_r)+\sum_{i=1}^{r-1}(r+1-i)\cdot \sigma^*(D^-_i)-\sum_{i=1}^{r-1}(r+1-i)\cdot \sigma^*(D^+_i)\,.
\end{equation*}
By the same argument as above, we can conclude that  $\sigma^*$ is the identity map or $\sigma^*=({\rm USD})^*$.
\smallskip

{\bf\noindent Case 4 ($r=p=n-s=1$).}
If $n=2$, Lemma \ref{sigma2} holds trivially. When $n\geq 3$, we have   $B_0=H-D^+_1$, $B_1=H$, and 
$K_{\mathcal T_{s,p,n}}=-nH+(n-2)D^+_1=-nB_0-2D^+_1$.
Hence, $\sigma^*$ is the identity map.

\smallskip

{\bf\noindent Case 5 ($r=p=1$ and $2\leq n-s\leq s$).} 
We have that $B_0=H-D^+_1$, $B_1=H-D^-_1$, $K_{\mathcal T_{s,p,n}}=-nH+(s-1)D^+_1+(n-s-1)D^-_1$, and $\sigma^*$ is a permutation of $\{B_0,B_1,D^+_1, D^-_1\}$. 
Since $\sigma^*(K_{\mathcal T_{s,p,n}})=K_{\mathcal T_{s,p,n}}$, we can conclude that $\sigma^*$ permutes $B_0$ and $B_1$. If $\sigma^*(B_0)=B_0$, then  $\sigma^*$ is the identity map.
If $\sigma^*(B_0)=B_1$ , then   $n=2s$ and $\sigma^*=({\rm USD})^*$.

\smallskip

{\bf\noindent Case 6 ($r=n-s=1$ and $2\leq p<s$).}
Similarly, we can conclude that $\sigma^*$ permutes $B_0$ and $B_1$. If $\sigma^*(B_0)=B_0$, then  $\sigma^*$ is the identity map.
If $\sigma^*(B_0)=B_1$ , then  $n=2p$ and $\sigma^*=({\rm DUAL})^*$.

\medskip

We complete the proof of Lemma \ref{sigma2}.
\,\,\,$\endpf$

\medskip

{\noindent\bf Proof of Theorem \ref{auto1}.}  Without loss of generality, we may assume that $2p\leq n\leq 2s$. We prove Theorem \ref{auto1} based on a case by case argument.
\smallskip

{\bf\noindent Case 1 ($p=s=n-s=1$).} $\mathcal T_{1,1,2}\cong\mathbb {P}^1$.
\smallskip

{\bf\noindent Case 2 ($p=n-s=1$ and $s\geq 2$).} Let $\sigma$ be an automorphism of $\mathcal T_{s,p,n}$. By Lemma \ref{sigma2}, $\sigma^*$ is the identity map. By Lemma \ref{des}, $\sigma$ induces an automorphism  $\Sigma$ of $G(p,n)$. Notice that each irreducible component of the exceptional divisor is invarinat under $\sigma$.   Therefore,  $\Sigma$ maps $\overline{\mathcal V_{(0,1)}^{+}}$ to $\overline{\mathcal V_{(0,1)}^{+}}$.

Recall that the automorphism group of $G(p,n)$ is generated by ${\rm PGL}(n,\mathbb K)$ when $n\neq 2p$. Write $\Sigma$ as a $n\times n$ matrix
\begin{equation}\label{block}
    \Sigma=\left(\begin{matrix}
    A&B\\
    C&D\\
    \end{matrix}\right)\,,
\end{equation}
where $A$, $B$, $C$, $D$ are submatrices of sizes $s\times s$, $s\times (n-s)$, $(n-s)\times s$, $(n-s)\times (n-s)$ respectively.
If  $C$ is not a zero matrix, there is a point $x\in \overline{\mathcal V_{(0,1)}^{+}}$ such that $\Sigma(x)\notin\overline{\mathcal V_{(0,1)}^{+}}$, which is a contradiction. Therefore,  $C$ is a zero matrix and $\Sigma\in P$ where $P$ is defined by (\ref{parabolic}).

\smallskip

{\bf\noindent Case 3 ($p=1$ and $2\leq n-s<s$).}  Let $\sigma\in{\rm Aut}(\mathcal T_{s,p,n})$.  Then $\sigma$ induces an automorphism $\Sigma$ of $G(p,n)$ that maps $\overline{\mathcal V_{(0,1)}^{+}}$ to $\overline{\mathcal V_{(0,1)}^{+}}$ and $\overline{\mathcal V_{(1,0)}^{-}}$ to $\overline{\mathcal V_{(1,0)}^{-}}$.
Write $\Sigma$ as (\ref{block}). If $C$ (resp. $B$) is not zero, there is a point $x\in \overline{\mathcal V_{(0,1)}^{+}}$ (resp. $\overline{\mathcal V_{(0,1)}^{-}}$) such that $\Sigma(x)\notin\overline{\mathcal V_{(0,1)}^{+}}$ (resp. $\overline{\mathcal V_{(0,1)}^{-}}$). Then  $\sigma\in {\rm GL}(s,\mathbb K)\times {\rm GL}(n-s,\mathbb K)$. 
\smallskip

{\bf\noindent Case 4 ($2\leq p$ and $2p<n<2s$).} Let $\sigma\in{\rm Aut}(\mathcal T_{s,p,n})$.  Then $\sigma$ induces an automorphism $\Sigma$ of $G(p,n)$ that maps $\overline{\mathcal V_{(p-r+1,r-1)}^{+}}$ to $\overline{\mathcal V_{(p-r+1,r-1)}^{+}}$ and  $\overline{\mathcal V_{(p-1,1)}^{-}}$ to $\overline{\mathcal V_{(p-1,1)}^{-}}$. Write $\Sigma$ as (\ref{block}). If $C$ (resp. $B$) is not zero, there is a point $x\in \overline{V_{(p-r+1,r-1)}^{+}}$ (resp. $\overline{\mathcal V_{(p-1,1)}^{-}}$) such that $\Sigma(x)\notin\overline{\mathcal V_{(p-r+1,r-1)}^{+}}$ (resp. $\overline{\mathcal V_{(p-r+1,r-1)}^{+}}$). Then 
$\sigma\in {\rm GL}(s,\mathbb K)\times {\rm GL}(n-s,\mathbb K)$.

\smallskip

{\bf\noindent Case 5 ($1\leq p<s=n-s$).} Let $\sigma\in{\rm Aut}(\mathcal T_{s,p,n})$. Then $\sigma$ or ${\rm USD}\circ\sigma$ induces an automorphism  $\Sigma$ of $G(p,n)$ and the identity map on the Picard group. We can conclude that 
$\Sigma\in {\rm GL}(s,\mathbb K)\times {\rm GL}(n-s,\mathbb K)$ similarly.

\smallskip

{\bf\noindent Case 6 ($2\leq p<s$ and $n=2p$).} Let $\sigma\in{\rm Aut}(\mathcal T_{s,p,2p})$. 
Then $\sigma$ or ${\rm DUAL}\circ\sigma$ induces an automorphism $\Sigma$ of $G(p,2p)$ and the identity map on the Picard group. Since the automorphism group of $G(p,2p)$ is generated by ${\rm PGL}(2p,\mathbb K)$ and the dual automorphism ${\empty}^*$ (see \S \ref{isomt}), when $\Sigma\in {\rm PGL}(2p,\mathbb K)$,  we can conclude that $\Sigma\in {\rm GL}(s,\mathbb K)\times {\rm GL}(2p-s,\mathbb K)$ similarly to Case 4.

If $\Sigma\not\in {\rm PGL}(2p,\mathbb K)$, we have that $\widetilde\Sigma:=\Sigma\circ {\empty}^*\in {\rm PGL}(2p,\mathbb K)$ where ${\empty}^*$ takes the form (\ref{dtrans}). Write $\widetilde\Sigma$ as as (\ref{block}). It is easy to verify that $\overline{\mathcal V_{(p,0)}^{-}}$ is invariant under  $\Sigma=\widetilde\Sigma\circ{\empty}^*$ only if $B,C$ are zero matrices, which is a contradiction. %This is a contradiction since $D_1^-$ is invariant under $\sigma$.

\smallskip

{\bf\noindent Case 7 ($2\leq p=s=n-s$).} Let $\sigma\in{\rm Aut}(\mathcal T_{s,p,2p})$. Then $\sigma$ or ${\rm USD}\circ\sigma$ induces an automorphism  $\Sigma$ of $G(p,2p)$ and the identity map on the Picard group. If $\Sigma\in {\rm PGL}(2p,\mathbb K)$, then $\Sigma\in {\rm GL}(s,\mathbb K)\times {\rm GL}(2p-s,\mathbb K)$ in the same way as in Case 4. Otherwise, ${\rm DUAL}\circ\sigma$ or
${\rm DUAL}\circ{\rm USD}\circ\sigma$ induces an automorphism  $\widetilde \Sigma$ of $G(p,2p)$ which is contained in ${\rm PGL}(2p,\mathbb K)$. Hence $\widetilde \Sigma\in {\rm GL}(s,\mathbb K)\times {\rm GL}(2p-s,\mathbb K)$ in the same way as in Case 4.

\smallskip

We complete the proof of Theorem \ref{auto1}.  \,\,\,$\endpf$

\subsection{The automorphism groups of \texorpdfstring{$\mathcal M_{s,p,n}$}{ff}} \label{symm}

Recall that
\begin{example}\label{ccolli}
Let $M_{p\times q}$ be the set of all $p\times q$ matrices and $\mathbb P(M_{p\times q})$ the projectivization of $M_{p\times q}$. It is clear that $\mathbb P(M_{p\times q})$ is a ${\rm GL}(p,\mathbb K)\times {\rm GL}(q,\mathbb K)$-variety under the left multiplication by the group ${\rm GL}(p,\mathbb K)$ and the right multiplication by the group ${\rm GL}(q,\mathbb K)$. Define $r:=\min\left\{p,q\right\}$.
Then $\mathbb P(M_{p\times q})$ has $r$ ${\rm GL}(p,\mathbb K)\times {\rm GL}(q,\mathbb K)$-orbits, whose closures $Z_{r}\supset Z_{r-1}\supset \cdots\supset Z_1$ are given by the condition
that $Z_i$ is the set of points corresponding to matrices of rank at most $i$.  Blowing up $\mathbb P(M_{p\times q})$ successively along the (strict transform of)
$Z_1,\cdots,Z_{r-1}$, we obtain a smooth variety $\widetilde {\mathbb P}(M_{p\times q})$ which is isomorphic to $\mathcal M_{p,p,p+q}$ by \cite[Theorem 1.6]{FW}. 
\end{example}
\begin{lemma}\label{pre}
Assume that $p\neq q$ are positive integers. Let $\sigma$ be an automorphism of $\mathcal M_{p,p,p+q}$ such that $\sigma^*$ is the identity map on the Picard group of $\mathcal M_{p,p,(p+q)}$.  Then $\sigma\in {\rm PGL}(p,\mathbb K)\times {\rm PGL}(q,\mathbb K)$.
\end{lemma}

{\noindent\bf Proof of Lemma \ref{pre}.} 
When $p=1$ or $q=1$, $\mathcal M_{p,p,p+q}$ is isomorphic to a projective space. Without loss of generality, we may assume that $2\leq p<q$.

Similarly to Lemma \ref{des} we can show that $\sigma$ descends to an automorphism $\Sigma$ of $P(M_{p\times q})$ which preserves the ranks of the matrices. Hence, to prove Lemma \ref{pre} it suffices to show that $\Sigma\in {\rm GL}(p,\mathbb K)\times {\rm GL}(q,\mathbb K)$ by viewing  $\Sigma$ as an element of ${\rm GL}(pq,\mathbb K)$.

For $1\leq i\leq p$ and $1\leq j\leq q$, denote by $E_{ij}$ the $p\times q$ matrix such that the $(i,j)^{\rm th}$ entry is $1$ and zero elsewhere. Write $\Sigma$  as a linear transformation
\begin{equation*}
\Sigma(E_{uv})=\sum\nolimits_{i=1}^{p}\sum\nolimits_{j=1}^{q}a_{uv}^{ij}\cdot E_{ij},\,\,\,\,\,\,a_{uv}^{ij}\in\mathbb K\,,\,\,1\leq u\leq p\,,\,\,1\leq v\leq q.
\end{equation*}
Since $\Sigma$ maps any rank-$1$ matrix to a rank-$1$ matrix, there are integers $1\leq i\leq p$ and $1\leq j\leq q$ such that $a_{11}^{ij}\neq 0$. Composing $\Sigma$ with an element of ${\rm GL}(p,\mathbb K)\times {\rm GL}(q,\mathbb K)$, we can assume that $\Sigma( E_{11})=E_{11}$. %$a_{11}^{11}\neq 0$. %; indeed there is an element  $\tau\in {\rm GL}(p,\mathbb K)\times {\rm GL}(q,\mathbb K)$ such that $(\tau\circ\Sigma) E_{11})=E_{11}$. By a slight abuse of notation, we use $\Sigma$ for any composition $\tau\circ\Sigma$, $\tau\in {\rm GL}(p,\mathbb K)\times {\rm GL}(q,\mathbb K)$, in the following.

We claim that there are integers $2\leq i\leq p$, $2\leq j\leq q$ such that $a_{22}^{ij}\neq 0$.  Assume the contrary. Then there are integers $2\leq i\leq p$, $2\leq j\leq q$ such that $a_{22}^{i1}\neq 0$ and $a_{22}^{1j}\neq 0$, for  $\Sigma(E_{11}+E_{22})$ has rank $2$. However, this contradicts the fact that $\Sigma(E_{22})$ has rank $1$. Composing $\Sigma$ with a certain element of ${\rm GL}(p,\mathbb K)\times {\rm GL}(q,\mathbb K)$  we can derive that $\Sigma(E_{11})=E_{11}$ and $\Sigma(E_{22})=E_{22}$. Similarly, we may assume that $\Sigma(E_{ii})=E_{ii}$  for $1\leq i\leq p$.

Since $E_{11}+\lambda E_{12}$ (resp. $E_{22}+\lambda E_{12}$) is of rank-$1$, $\lambda\in\mathbb K$, we conclude that
\begin{equation*}%\label{150}
\Sigma(E_{12})=\sum_{i=1}^{p}a_{12}^{i1}\cdot E_{i1}\,\,\,{\rm or}\,\,\,\sum_{j=1}^{q}a_{12}^{1j}\cdot E_{1j}\,\,\left({\rm resp}.\,\, \Sigma(E_{12})=\sum_{i=1}^{p}a_{12}^{i2}\cdot E_{i2}\,\,\,{\rm or}\,\,\,\sum_{j=1}^{q}a_{12}^{2j}\cdot E_{2j}\right).
\end{equation*}
Then $\Sigma(E_{12})=a_{12}^{12}\cdot E_{12}$  or $a_{12}^{21}\cdot E_{21}$.
We claim that $\Sigma(E_{12})=a_{12}^{12}\cdot E_{12}$. Assume the contrary. Then $\Sigma(E_{1i})=a_{1i}^{i1}\cdot E_{i1}$ with $a_{1i}^{i1}\neq 0$ for $2\leq i\leq p$. Regardless of the value of $\Sigma(E_{1(p+1)})$, it will violate the assumption that the ranks are preserved under $\Sigma$.

By the same argument, we can assume that
\begin{enumerate}[label=(\alph*)]
    \item $\Sigma(E_{ij})=a^{ij}_{ij}\cdot E_{ij}$ for $1\leq i\leq p$ and $1\leq j\leq q$.
    \item $a^{ij}_{ij}=1$ when $i\equiv j\,\, ({\rm mod}\, p)$.
\end{enumerate}
Composing an element of  ${\rm PGL}(p,\mathbb K)$ by the left action, we can fix $a^{21}_{21}=\cdots a^{p1}_{p1}=1$. Then  $a^{ij}_{ij}=1$, $1\leq i\leq p$, $1\leq j\leq q$, for $E_{11}+E_{1j}+E_{i1}+E_{ij}$ has rank $1$. We thus conclude that $\Sigma\in {\rm GL}(p,\mathbb K)\times {\rm GL}(q,\mathbb K)$.

The proof of Lemma \ref{pre} is complete.
\,\,\,$\endpf$
\medskip

\begin{lemma}\label{fibpre}
Assume that $2p\leq n\leq 2s$ and $2\leq p<s$. Let $\sigma$ be an automorphism of $\mathcal M_{s,p,n}$ such that $\sigma^*$ is the identity map on the Picard group of $\mathcal M_{s,p,n}$.  Then 
$\sigma\in {\rm PGL}(s,\mathbb K)\times {\rm PGL}(n-s,\mathbb K)$.
\end{lemma}

Before proceeding, we recall that by \cite[\S 5.2]{FW},
$\mathcal M_{s,p,n}$ is an iterated blow-up of $\mathbb P(N_{\mathcal V_{(p,0)}/G(p,n)})$ the projectivization of the normal bundle $N_{\mathcal V_{(p,0)}/G(p,n)}$ of the sub-Grassmannian $\mathcal V_{(p,0)}\cong G(p,s)$ in $G(p,n)$. In the following, we define convenient local coordinate charts for $\mathbb P(N_{\mathcal V_{(p,0)}/G(p,n)})$. 

For $1\leq i_1<\cdots<i_p\leq s$, we define an open subset $U_{i_1\cdots i_p}$ of $\mathcal V_{(p,0)}$ by
\begin{equation*}
\left\{\underbracedmatrixl{
 \cdots &\widetilde x_{1(i_1-1)}&1& \cdots&\widetilde x_{1(i_2-1)}&0& \cdots&\widetilde x_{1(i_p-1)}&0& \cdots  \\
  \cdots &\widetilde x_{2(i_1-1)}&0& \cdots&\widetilde x_{1(i_2-1)}&1& \cdots&\widetilde x_{2(i_p-1)}&0& \cdots  \\
  \ddots&\vdots&\vdots&\ddots&\vdots&\vdots&\ddots&\vdots&\vdots&\ddots\\
  \cdots &\widetilde x_{p(i_1-1)}&0& \cdots&\widetilde x_{1(i_2-1)}&0& \cdots&\widetilde x_{p(i_p-1)}&1& \cdots \\}{s\,\rm columns}\hspace{-.22in}
\begin{matrix}
  &\hfill\tikzmark{e}\\
  \\
  \\
  \\
  &\hfill\tikzmark{f}\end{matrix}\hspace{-.1in}
  \begin{matrix}
  &\hfill\tikzmark{a}\\
  \\
  \\
  \\
  &\hfill\tikzmark{b}\end{matrix}\,
  \underbracedmatrixr{ 0  &0& \cdots & 0\\
 0  &0 & \cdots&0\\
  \vdots&\vdots&\ddots&\vdots \\
 0  &0 & \cdots& 0\\} {(n-s) \,\rm columns}\right\}.
  \tikz[remember picture,overlay]   \draw[dashed,dash pattern={on 4pt off 2pt}] ([xshift=0.5\tabcolsep,yshift=7pt]a.north) -- ([xshift=0.5\tabcolsep,yshift=-2pt]b.south);\tikz[remember picture,overlay]   \draw[dashed,dash pattern={on 4pt off 2pt}] ([xshift=0.5\tabcolsep,yshift=7pt]e.north) -- ([xshift=0.5\tabcolsep,yshift=-2pt]f.south);
\end{equation*}
After trivializing over $U_{i_1\cdots i_p}$, we obtain an open subset $U_{i_1\cdots i_p}^N\subset N_{\mathcal V_{(p,0)}/G(p,n)}$ with coordinates
\begin{equation*}%\label{trr1}
\left(\left(\cdots,\widetilde x_{uv},\cdots\right)_{\substack{1\leq u\leq p,\,1\leq v\leq s\\v\neq i_1,\cdots,i_p}}\,,\,\, (\cdots, \widetilde a_{ij},\cdots)_{1\leq i\leq p,\,s+1\leq j\leq n}\right)=:(\widetilde X,\widetilde A)_{i_1\cdots i_p}\,,
\end{equation*}
and an open embedding $M^{i_1\cdots i_p}:U_{i_1\cdots i_p}^N\hookrightarrow G(p,n)$ defined by $M((\widetilde X,\widetilde A)_{i_1\cdots i_p}):=$ 
\begin{equation}\label{mnin1}%\label{tri1}
 \underbracedmatrixl{
 \cdots &\widetilde x_{1(i_1-1)}&1& \cdots&\widetilde x_{1(i_p-1)}&0& \cdots  \\
  \cdots &\widetilde x_{2(i_1-1)}&0& \cdots&\widetilde x_{2(i_p-1)}&0& \cdots  \\
  \ddots&\vdots&\vdots&\ddots&\vdots&\vdots&\ddots\\
  \cdots &\widetilde x_{p(i_1-1)}&0& \cdots&\widetilde x_{p(i_p-1)}&1& \cdots \\}{s\,\rm columns}\hspace{-.22in}
\begin{matrix}
  &\hfill\tikzmark{e}\\
  \\
  \\
  \\
  &\hfill\tikzmark{f}\end{matrix}\hspace{-.1in}
  \begin{matrix}
  &\hfill\tikzmark{a}\\
  \\
  \\
  \\
  &\hfill\tikzmark{b}\end{matrix}\,
  \underbracedmatrixr{\widetilde a_{1(s+1)}  & \cdots& \widetilde a_{1n}\\
 \widetilde a_{2(s+1)}   & \cdots& \widetilde a_{2n}\\
  \vdots&\ddots&\vdots \\
 \widetilde a_{p(s+1)}  & \cdots& \widetilde a_{pn}\\} {(n-s) \,\rm columns}.
  \tikz[remember picture,overlay]   \draw[dashed,dash pattern={on 4pt off 2pt}] ([xshift=0.5\tabcolsep,yshift=7pt]a.north) -- ([xshift=0.5\tabcolsep,yshift=-2pt]b.south);\tikz[remember picture,overlay]   \draw[dashed,dash pattern={on 4pt off 2pt}] ([xshift=0.5\tabcolsep,yshift=7pt]e.north) -- ([xshift=0.5\tabcolsep,yshift=-2pt]f.south);
\end{equation}
Let $U_{i_1\cdots i_p}^P$ be the open subset of $\mathbb P(N_{\mathcal V_{(p,0)}/G(p,n)})$ over $U_{i_1\cdots i_p}$ with coordinates
\begin{equation}\label{trr1}  
\left(\left(\cdots,\widetilde x_{uv},\cdots\right)_{\substack{1\leq u\leq p,\,1\leq v\leq s\\v\neq i_1,\cdots,i_p}}\,,\,\, [\cdots, \widetilde a_{ij},\cdots]_{\substack{1\leq i\leq p,\\s+1\leq j\leq n}}\right)=:(\widetilde X,[\widetilde A\,])_{i_1\cdots i_p}\,,
\end{equation}
where $[\cdots, \widetilde a_{ij},\cdots]$ are the homogeneous coordinates for the fiber $\mathbb {P}^{p(n-s)-1}$.

It is easy to verify that $(\widetilde X,[\widetilde A\,])_{i_1\cdots i_p}$ and $(\widetilde X,[\widetilde A\,])_{i^{\prime}_1\cdots i^{\prime}_p}$ represent the same point of $\mathbb P(N_{\mathcal V_{(p,0)}/G(p,n)})$  if and only if the matrices $M_{(\widetilde X,\widetilde A)_{i_1\cdots i_p}}$ and $M_{(\widetilde X,\widetilde A)_{i^{\prime}_1\cdots i^{\prime}_p}}$ represent the same point of $G(p,n)$ up to the $\mathbb G_m$-action, that is,  there exists a matrix $W\in {\rm GL}(p,\mathbb K)$ and  $\lambda\in\mathbb K$ such that 
\begin{equation*}
    M_{(\widetilde X,\widetilde A)_{i_1\cdots i_p}}=W\cdot M_{(\widetilde X,\widetilde A)_{i^{\prime}_1\cdots i^{\prime}_p}}\cdot \left(\begin{matrix} I_{s\times s}&0\\ 0&\lambda\cdot I_{(n-s)\times (n-s)}\\ \end{matrix} \right)\,.
\end{equation*}

{\bf\noindent Proof of Lemma \ref{fibpre}.} %We provide an alternative proof of  Lemma \ref{fibpre} which is based on the semi-positivity of $K_{\mathcal T_{s,p,n}}$ and the Pl\"ucker relations. Let $\sigma$ be an automorphism of $\mathcal M_{s,p,n}$ which induces the identity map on the Picard group. Similarly, we can show that  the line bundle $-K_{\mathcal T_{s,p,n}}+H_j\big|_{\mathcal T_{s,p,n}}-D_1^-$ on $\mathcal T_{s,p,n}$ is nef and big for  $0\leq j\leq r$. Then 
If $p=n-s$, Lemma \ref{fibpre} follows from Lemma \ref{pre}. Without loss of generality,  we may assume that $2p\leq n$, $p\neq n-s$ and $2\leq p<s$ in the following.

By Lemma \ref{mcls} and the Kawamata--Viehweg vanishing theorem, the following natural map is surjective.
\begin{equation*}
\Gamma(\mathcal T_{s,p,n},\,\,H_j)\rightarrow \Gamma(\mathcal M_{s,p,n},\,\,H_j|_{\mathcal M_{s,p,n}})\,.
\end{equation*}
By Lemma \ref{cls}, we conclude that the complete linear series of $H_j|_{\mathcal M_{s,p,n}}$ on $\mathcal M_{s,p,n}$ is isomorphic to $\mathbb {P}^{N^j_{s,p,n}}$ as well. We can thus show that $\sigma$ is induced by an automorphism of the ambient space $\mathbb {P}^{N_{p,n}}\times\mathbb {P}^{N^0_{s,p,n}}\times\cdots\times\mathbb {P}^{N^p_{s,p,n}}$.

Since $\sigma^*$ is the identity map on the Picard group of $\mathcal M_{s,p,n}$, similarly  to Lemma \ref{des} we can conclude that $\sigma$ induces an automorphism $\Sigma$ of $\mathbb P(N_{\mathcal V_{(p,0)}/G(p,n)})$ such that the following diagram commutes.
\begin{equation}
\begin{array}{ccccc}
\vspace{.02in}
\mathcal M_{s,p,n}&\xrightarrow{\,\,\,\,\,\sigma\,\,\,\,\,} &\mathcal M_{s,p,n}&&\\
\vspace{-.03in}
\Big\downarrow\llap{$\scriptstyle \tau_{s,p,n}\,\,\,\,\,$}&&\Big\downarrow\rlap{$\scriptstyle \tau_{s,p,n}\,\,\,\,\,$}&&\\ \vspace{.02in}
\mathbb P(N_{\mathcal V_{(p,0)}/G(p,n)})&\xrightarrow{\,\,\,\,\,\Sigma\,\,\,\,\,}&\mathbb P(N_{\mathcal V_{(p,0)}/G(p,n)})&\xhookrightarrow{\,\,\,i\,\,\,}&\mathbb{P}^{N_{p,n}}\times\mathbb{P}^{N^1_{s,p,n}}\\
\vspace{-.03in}
\Big\downarrow\llap{$\scriptstyle \kappa_{s,p,n}\,\,\,\,\,$}&&\Big\downarrow\rlap{$\scriptstyle \kappa_{s,p,n}\,\,\,\,\,$}&&\\
\mathcal V_{(p,0)}&\xrightarrow{\,\,\,\,\,\widehat\Sigma\,\,\,\,\,}&\mathcal V_{(p,0)}&&\\%&\hookrightarrow&\mathbb{CP}^{N_{p,n}}\times\mathbb{CP}^{N^1_{s,p,n}}\\
\end{array}
\end{equation}
Moreover, $\Sigma$ maps a fiber of $\kappa_{s,p,n}$ to another fiber, and hence $\sigma$ induces an automorphism $\widehat\Sigma$ of the base $\mathcal  V_{(p,0)}$. By  Proposition  \ref{tfiber} and Example \ref{ccolli}, we conclude that each fiber of $\kappa_{s,p,n}$ is isomorphic to $P(M_{p\times (n-s)})$, and each fiber of $\kappa_{s,p,n} \circ\tau_{s,p,n}$ is isomorphic to the variety of complete collineations $\widetilde P(M_{p\times (n-s)})\cong\mathcal M_{p,p,n-s+p}$. Since $\sigma^*$ is the identity map and hence the irreducible components of the exceptional divisor are $\sigma$-invariant, $\Sigma$ preserves the ranks of the matrices in $P(M_{p\times (n-s)})$. Similarly to Lemma \ref{pre}, we can derive that the restriction of $\Sigma$ to each fiber of $\kappa_{s,p,n}$ is given by an element of ${\rm GL}(p,\mathbb K)\times {\rm GL}(n-s,\mathbb K)$. 

We next prove Lemma \ref{fibpre} based on a case by case argument. 
\smallskip

{\bf\noindent Case 1 ($2\leq p<s$ and $s\neq 2p$).} After composing $\sigma$ with an element of ${\rm PGL}(s,\mathbb K)\times {\rm PGL}(n-s,\mathbb K)$, we can assume that the following diagram commutes.
\begin{equation*}%\label{st}
\begin{tikzcd}
&\mathbb P(N_{\mathcal V_{(p,0)}/G(p,n)})\arrow{r}{\,\,\Sigma}\arrow[hookrightarrow]{d}{\,\,i}  &\mathbb P(N_{\mathcal V_{(p,0)}/G(p,n)})\arrow[hookrightarrow]{d}{\,\,i}\\
&\mathbb{P}^{N_{p,n}}\times\mathbb{P}^{N^1_{s,p,n}}\arrow{r}{\,\,({\rm Id},\,\tau)\,\,\,\,\,\,} &\mathbb{P}^{N_{p,n}}\times\mathbb{P}^{N^1_{s,p,n}} \\
\end{tikzcd}\vspace{-20pt}\,,
\end{equation*} 
where $\tau$ is an automorphism of $\mathbb{P}^{N^1_{s,p,n}}$.

In the following, we will show that
$\tau\in {\rm PGL}(s,\mathbb K)\times {\rm PGL}(n-s,\mathbb K)$.

Let $U_{(s-p+1)\cdots s}^P\subset\mathbb P(N_{\mathcal V_{(p,0)}/G(p,n)})$ be the open subset with the local coordinates $(\widetilde X,[\widetilde A\,])_{(s-p+1)\cdots s}$ defined in (\ref{trr1}) by setting $i_1=s,i_2=s-1,\cdots,i_p=s-p+1$. For convenience, we denote $U_{(s-p+1)\cdots s}^P$ by $U^P$ and $(\widetilde X,[\widetilde A\,])_{(s-p+1)\cdots s}$ by $(X, [A])$. Then, the embedding $i$ takes the form
\begin{equation*}%\label{coordinatesemb}
 i((X,[ A\,]))=   [\cdots,P_I(M_{(X,A)}),\cdots]_{I\in\mathbb I_{p,n}}\times[\cdots,P_I(M_{(X,A)}),\cdots]_{I\in\mathbb I^1_{s,p,n}}\,
\end{equation*}
where $P_I$ is the Pl\"ucker coordinate function and $M_{(X, A)}$ is defined by (\ref{mnin1}).  To distinguish the source and target $U^P$, we denote the local coordinates of $U^P$ in the target by 
\begin{equation*}%\label{yuv}
(\cdots,y_{uv},\cdots)_{1\leq u\leq p,\,1\leq v\leq s-p}\times [\cdots, b_{ij},\cdots]_{1\leq i\leq p,\,s+1\leq j\leq n}=:(Y,[B])\,.
\end{equation*}

Since $\tau$ is a projective linear transformation of $\mathbb{P}^{N^1_{s,p,n}}$, we can represent $(Y,[B])=\Sigma((X,[A]))$ in the above local coordinate charts by
\begin{equation}\label{cij}
\begin{split}
&\,\,\,y_{uv}=x_{uv}\,,\,\,\,1\leq u\leq p\,,\,\,1\leq v\leq s-p\,,\\
&[\cdots,P_{I}(M_{(Y,B)}),\cdots]_{I\in\mathbb I^1_{s,p,n}}=[\cdots,\sum_{J\in\mathbb I^1_{s,p,n}}C_{I}^{J}\cdot P_J(M_{(X,A)}),\cdots]_{I\in\mathbb I^1_{s,p,n}}\,,
\end{split}
\end{equation}
where $C^J_{I}\in\mathbb K$. 

We introduce the following notation. 
For $1\leq k\leq p$ and $1\leq t\leq s-p$, 
\begin{equation*}
I_{kt}:=(s,s-1,\cdots,\widehat{s-p+k},\cdots,s-p+2,s-p+1,t)\,.
\end{equation*}
For $1\leq k\leq p$ and $s+1\leq t\leq n$, 
\begin{equation*}
I^*_{kt}:=(t,s,s-1,\cdots,\widehat{s-p+k},\cdots,s-p+2,s-p+1)\,.
\end{equation*}
For $s-p+1\leq \alpha_1,\alpha_2\leq s$ and $\alpha_1\neq\alpha_2$, $1\leq \beta\leq s-p$ and $s+1\leq \gamma\leq n$,  
\begin{equation*}
I^{\alpha_1\alpha_2}_{\beta\gamma}:=(\gamma,s,s-1,\cdots,\widehat{\alpha_1},\cdots,\widehat{\alpha_2},\cdots,s-p+2,s-p+1,\beta)\,.
\end{equation*}

Then, computation yields that 
\begin{equation*}
P_{I_0}(M_{(X,A)})=1\,, \,\,P_{I_{kt}}(M_{(X,A)})=(-1)^{k-1}\cdot  x_{kt}\,, \,\,
 P_{I^*_{kt}}(M_{(X,A)})=(-1)^{p-k}\cdot a_{kt}\,, \,\,
\end{equation*}
and $P_{I^{\alpha_1\alpha_2}_{\beta\gamma}}(M_{(X,A)})=$
\begin{equation*}
\left\{\begin{array}{cc}
   (-1)^{p+\alpha_1-\alpha_2}\left( x_{(\alpha_1-s+p)\beta}\cdot a_{(\alpha_2-s+p)\gamma}-x_{(\alpha_2-s+p)\beta}\cdot a_{(\alpha_1-s+p)\gamma}\right)  &{\rm if\,\,} \alpha_1>\alpha_2 \\
(-1)^{p-1+\alpha_1-\alpha_2}\left( x_{(\alpha_1-s+p)\beta}\cdot a_{(\alpha_2-s+p)\gamma}-x_{(\alpha_2-s+p)\beta}\cdot a_{(\alpha_1-s+p)\gamma}\right)     & {\rm if\,\,}\alpha_1<\alpha_2\\
\end{array}\right..
\end{equation*}

{\bf\noindent Claim.}  Let $C_{I}^J$ be the constants in (\ref{cij}). Assume that $C_{I^*_{ij}}^J\neq0$ for a certain index $J\in\mathbb I^1_{s,p,n}$ and integers $1\leq i\leq p$, $s+1\leq j\leq n$. Then  $J=I^*_{kt}$ where $1\leq k\leq p$ and $s+1\leq t\leq n$.  
\smallskip

{\bf\noindent Proof of Claim.} Suppose  that $C_{I^*_{ij}}^{J^*}\neq0$ for integers $1\leq i\leq p$, $s+1\leq j\leq n$, and an index  $J^*:=(i_1,\cdots,i_p)\in\mathbb I^1_{s,p,n}\backslash\left\{I^*_{rt}|1\leq r\leq p,s+1\leq t\leq n\right\}$.  Then $i_p\leq s-p$. Take an integer $1\leq i^{\prime}\leq p$ such that \begin{equation*}
i^{\prime}\not\in\{i,i_2-s+p,i_3-s+p,\cdots, i_{p-1}-s+p\}.    
\end{equation*}

Recall that the Grassmannian $G(p,n)$ is a subvariety of $\mathbb {P}^{N_{p,n}}$ defined by the Pl\"ucker relations.
In particular, we have 
\begin{equation*}%\label{plcuker}
  P_{I^{(s-p+i^{\prime})(s-p+i)}_{i_pj}}\cdot P_{I_0}=\left\{\begin{array}{cc}
     P_{I_{i^{\prime}i_p}}\cdot P_{I^*_{ij}}- P_{I_{ii_p}}\cdot P_{I^*_{i^{\prime}j}}  & {\rm if\,\,} i<i^{\prime}\\
     -\left(P_{I_{i^{\prime}i_p}}\cdot P_{I^*_{ij}}- P_{I_{ii_p}}\cdot P_{I^*_{i^{\prime}j}} \right) & {\rm if\,\,} i>i^{\prime}
  \end{array}\right.\,.
\end{equation*}
Since the image $\Sigma(U^P)$ satisfies the same equation, substituting (\ref{cij}) into the Pl\"ucker relations, we have that
\begin{equation*}
\begin{split}
 \sum_{I\in\mathbb I^1_{s,p,n}}C_{I^{(s-p+i^{\prime})(s-p+i)}_{i_pj}}^{J}&\cdot P_J(M_{(X,A)})=\pm x_{i^{\prime}i_p}\cdot \left(\sum_{I\in\mathbb I^1_{s,p,n}}C_{I^*_{ij}}^{J} P_J(M_{(X,A)})\right)\\
 &\mp x_{ii_p}\cdot \left(\sum_{I\in\mathbb I^1_{s,p,n}}C_{I^*_{i^{\prime}j}}^{J} P_J(M_{(X,A)})\right)\,.
\end{split}
\end{equation*}
Then the following equality holds as polynomials.
\begin{equation}\label{contr}
\begin{split}
&x_{i^{\prime}i_p}\left(C_{I^*_{ij}}^{J^*}\cdot P_{J^*}(M_{(X,A)})+\sum_{I\in\mathbb I^1_{s,p,n}\,,\,\,J\neq J^*}C_{I^*_{ij}}^{J} P_J(M_{(X,A)})\right)\\
 &=\pm\sum_{J\in\mathbb I^1_{s,p,n}}C_{I^{(s-p+i^{\prime})(s-p+i)}_{i_pj}}^{J} P_I(M_{(X,A)})+ x_{ii_p} \left(\sum_{J\in\mathbb I^1_{s,p,n}}C_{I^*_{i^{\prime}j}}^{J} P_I(M_{(X,A)})\right)\,.
\end{split}
\end{equation}
Notice that the first term on the left hand side of (\ref{contr}) contains a nonzero monomial with a factor $x_{i^{\prime}i_p}^2$, and that there is no nonzero monomial with a factor $x_{i^{\prime}i_p}^2$ on the right hand side. This is a contradiction.

We complete the proof of Claim.\,\,\,$\endpf$
\smallskip

By Claim, we have that
\begin{equation*}
b_{i^{\prime}j^{\prime}}\left(\Sigma(X,[A])\right)=\sum_{i=1}^p\sum_{j=s+1}^n(-1)^{i-i^{\prime}}\cdot C_{I^*_{i^{\prime}j^{\prime}}}^{I^*_{ij}}\cdot a_{ij}\,\,,\,\,\,{\,\,\,\,\,\,C_{I^*_{i^{\prime}j^{\prime}}}^{I^*_{ij}}}\in\mathbb K\,.\\
\end{equation*}
Since $\Sigma$ preserves the ranks of the matrices in $P(M_{p(n-s)})$, similarly to Lemma \ref{pre}, we can show that there is an element $\sigma_l\in {\rm PGL}(p,\mathbb K)$ and an element $\sigma_r\in {\rm PGL}(n-s,\mathbb K)$ such that
\begin{equation*}
\left(b_{i^{\prime}j^{\prime}}\left(\Sigma(X,[A])\right)\right)=\sigma_l\cdot (a_{ij})\cdot \sigma_r\,.\\
\end{equation*}
Composed of a ${\rm PGL}(n-s,\mathbb K)$-action, we can assume that the automorphism $\Sigma$ is given by the left action $\sigma_l$.

We claim that $b_{ij}\left(\Sigma(X,[A])\right)= C_{I^*_{ij}}^{I^*_{ij}}\cdot a_{ij}$ for each  $1\leq i\leq p$, $s+1\leq j\leq n$. Otherwise, there are integers $1\leq i\leq p$, $1\leq i^{\prime}\leq p$, $i\neq i^{\prime}$, $s+1\leq j\leq n$ such that $C_{I^*_{ij}}^{I^*_{i^{\prime}j}}\neq 0$. However, then 
$P_{I^{(s-p+i)(s-p+i^{\prime})}_{(s-p)j}}\left(\Sigma(X,[A])\right)$ contains a monomial with a factor $a_{i^{\prime}j}\cdot x_{i^{\prime}(s-p)}$ which can not be a linear combination of $P_I(M_{(X,[A])})$. This is a contradiction.  
 
Similarly, we can conclude that $C_{I^*_{ij}}^{I^*_{ij}}=C_{I^*_{i^{\prime}j}}^{I^*_{i^{\prime}j}}$ for $1\leq i,i^{\prime}\leq p$. Then $\sigma_l$ is trivial. We complete the proof of Lemma \ref{fibpre} when $2p\neq s$.
\smallskip

{\bf\noindent Case 2 ($2\leq p$ and $s=2p$).} Recall that the automorphism group of $G(p,2p)$ is generated by ${\rm PGL}(2p,\mathbb K)$ and the dual automorphism. Similarly  to Case 2 in the first proof of Lemma \ref{fibpre}, we can show that $\Sigma$ induces an automorphism of $G(p,2p)$ in ${\rm PGL}(2p,\mathbb K)$. Then the same argument in  Case 1 applies here.
\medskip

We complete the proof of Lemma \ref{fibpre}.\,\,\,$\endpf$

\medskip

Similarly to Lemma \ref{sigma2}, we have 
\begin{lemma}\label{msigma}
Assume that $2p\leq n\leq 2s$ and $2\leq p<s$. Let $\sigma$ be an automorphism of $\mathcal M_{s,p,n}$ and $\sigma^*$ the induced automorphism of the Picard group. When $n\neq 2s$ and $n\neq 2p$, $\sigma^*$ is the identity map. When $n=2s$, $\sigma^*$ is the identity map or $({\rm Usd})^*$. When  $n=2p$, $\sigma^*$ is the identity map or $({\rm Dual})^*$.
\end{lemma}
{\noindent\bf Proof of Lemma \ref{msigma}.} By Lemmas \ref{checkgb} and \ref{cone}, Eff$(\mathcal M_{s,p,n})$ the cone  of effective divisors of $\mathcal M_{s,p,n}$ is generated by 
\begin{equation*}
 \{\check D_{2}, \check D_{3}, \cdots, \check D_{r}, \check B_0, \check B_1, \cdots, \check B_r\}\,.
\end{equation*}
By (\ref{mb=b}), the set of extremal rays of Eff$(\mathcal M_{s,p,n})$ is given by
\begin{equation*}
\mathfrak G=\left\{ \begin{array}{ll}
\{\check D_{2}, \check D_{3}, \cdots, \check D_{r}, \check B_0, \check B_r\} &{\rm when}\,\,p\neq n-s\,\,{\rm and}\,\,p<s\\
\{\check D_{2}, \check D_{3}, \cdots, \check D_{r-1}, \check B_0,  \check B_{r}\}&{\rm when}\,\,n-s=p<s\,\,\\
\end{array}\right..
\end{equation*}
Since $\sigma^*$ preserves Eff$(\mathcal M_{s,p,n})$, it induces a permutation of $\mathfrak G$.

We prove Lemma \ref{msigma} based on a case by case argument. For convenience, we denote by $H$ the divisor $(\check R_{s,p,n})^*(\mathcal O_{ G(p,n)}(1))$.
\smallskip

{\bf\noindent Case 1 ($n-s=1$).} 
Eff$(\mathcal M_{s,p,n})$ is generated by $\check B_0,\check B_1$. According to Lemma \ref{mkan} and (\ref{mb=b}), we have $K_{\mathcal M_{s,p,n}}=-\left(n-p\right)\check B_0-p\check B_1$. If $\sigma^*$ is not the identity map, then $n=2p$ and $\sigma^*=({\rm Dual})^*$.
\smallskip

{\bf\noindent Case 2 ($p<s$, $p\neq n-s$, and $n-s\geq 2$).} 
By Lemma \ref{mkan}, we have that \begin{equation}\label{mkan3}
\begin{split}
&K_{\mathcal M_{s,p,n}}=-n\check B_0+\frac{p(n-s)}{r}\left(\check B_0-\check B_`r-\sum_{k=2}^r (r+1-k)\check D_k\right)\\
&\,\,\,\,\,\,\,\,\,\,\,\,\,\,\,\,\,\,+\sum_{i=2}^{r}\big((p-i+1)(n-s-i+1)-1\big) \check D_{i}\\
&=\sum_{i=2}^{r}\left((p-i+1)(n-s-i+1)-1-\frac{p(n-s)(r+1-i)}{r}\right)\check D_{i}\\
&\,\,\,\,\,\,\,\,\,\,\,\,\,\,\,\,\,\,-\left(n-\frac{p(n-s)}{r}\right)\check B_0-\frac{p(n-s)}{r}\check B_r\,.
\end{split}
\end{equation}

We first assume that $r=2$.
If $p=2$,  $K_{\mathcal M_{s,p,n}}=-2\check D_{2}-s\check B_0-(n-s)\check B_r$; if $n-s=2$,  $K_{\mathcal M_{s,p,n}}=-2\check D_{2}-(n-p)\check B_0-p\check B_r$. Then $\sigma^*$ is the identity.

Let $r\geq 3$. By (\ref{mb=b}) we have 
\begin{equation}\label{iden1}
\sigma^*(\check D_1)=\frac{1}{r}\sigma^*(\check B_0)-\frac{1}{r}\sigma^*(\check B_r)-\sum_{i=2}^r\frac{r+1-i}{r} \sigma^*(\check D_i)\,.
\end{equation}
(\ref{iden1}) has integer coefficients with respect to the basis
$\left\{H, \check D_2,\cdots,\check D_r\right\}$ by  Corollary \ref{mpicb}. Checking the coefficient of $H$, we can conclude the following possibilities.
\begin{enumerate}[label=(\alph*)]
\item $\sigma^*$ is the identity map.
\item $\sigma^*(\check B_0)=\check B_r$,  $\sigma^*(\check B_r)=\check B_0$.
\item There is  $2\leq l\leq r$ such that $\sigma^*(\check D_l)=\check B_0$ and $\sigma^*(\check D_{r+2-l})=\check B_r$.
\item $\sigma^*(\check B_0)=\check B_0$ and $\sigma^*(\check D_{r})=\check B_r$, or $\sigma^*(\check B_0)=\check B_r$ and $\sigma^*(\check D_{r})=\check B_0$.
\item $\sigma^*(\check B_r)=\check B_0$ and $\sigma^*(\check D_{2})=\check B_r$, or $\sigma^*(\check B_r)=\check B_r$ and $\sigma^*(\check D_2)=\check B_0$.
\end{enumerate}

For Case (b), by checking the coefficients of $\check D_2,\cdots, \check D_r$  in (\ref{iden1})  we can conclude that $\sigma^*(\check D_i)=\check D_{r+2-i}$ for $2\leq i\leq r$. By (\ref{mkan3}) we have 
$n-\frac{p(n-s)}{r}=\frac{p(n-s)}{r}$.
Hence $n=2s$ and $\sigma^*=({\rm Usd})^*$,  or $n=2p$ and $\sigma^*=({\rm Dual})^*$.

For Cases (c), (d) (e) we can conclude by (\ref{mkan3}) that $p=n-s$ or $p=s$ which contradicts the assumption.

\smallskip

{\bf\noindent Case 3 ($p=n-s<s$).} The proof is similar to that of Case 2. We omit the details for simplicity.

We complete the proof of Lemma \ref{msigma}.
\,\,\,$\endpf$

\medskip

{\noindent\bf Proof of Theorem \ref{mauto}.}   
We  prove Theorem \ref{mauto} based on a case by case argument. Note that we may assume that $2p\leq n\leq 2s$.

\smallskip

{\bf\noindent Case 1 ($p=1$).} If $p=s=1$ and $n=2$, $\mathcal M_{1,1,2}$ is trivially a point.
If $p=n-s=1$ and $n\geq 3$, $\mathcal M_{s,1,s+1}\cong\mathbb {P}^{n-1}$ and  ${\rm Aut}(\mathcal M_{s,1,s+1})={\rm PGL}(n-1,\mathbb K)$. If $p=1$, $2\leq n-s<s$, $\mathcal M_{s,1,2s}\cong\mathbb {P}^{s-1}\times \mathbb {P}^{n-s-1}$ and  
\begin{equation*}
    {\rm Aut}(\mathcal M_{s,1,n})={\rm PGL}(s,\mathbb K)\times {\rm PGL}(n-s,\mathbb K).
\end{equation*}
If $p=1$, $2\leq n-s=s$, $\mathcal M_{s,1,n}\cong\mathbb {P}^{s-1}\times \mathbb {P}^{s-1}$ and  
\begin{equation*}
    {\rm Aut}(\mathcal M_{s,1,2s})={\rm PGL}(s,\mathbb K)\times {\rm PGL}(s,\mathbb K)\rtimes \mathbb Z/2\mathbb Z.
\end{equation*}

{\bf\noindent Case 2 ($2\leq p$, $n\neq 2s$ and $n\neq 2p$).} 
Let $\sigma\in{\rm Aut}(\mathcal M_{s,p,n})$.
By Lemma \ref{msigma}, $\sigma^*$ is the identity map on the Picard group of $\mathcal M_{s,p,n}$. By By Lemma \ref{fibpre}, we have that $\sigma\in {\rm PGL}(s,\mathbb K)\times {\rm PGL}(n-s,\mathbb K)$.

\smallskip

{\bf\noindent Case 3 ($2\leq p<s$ and $n=2s$).}
Let $\sigma\in{\rm Aut}(\mathcal M_{s,p,n})$.  By Lemmas \ref{msigma},  \ref{fibpre}, either $\sigma$ or ${\rm Usd}\circ\sigma$ is in ${\rm PGL}(s,\mathbb K)\times {\rm PGL}(n-s,\mathbb K)$. Then, \begin{equation*}
    {\rm Aut}(\mathcal M_{s,p,2s})={\rm PGL}(s,\mathbb K)\times {\rm PGL}(n-s,\mathbb K)\rtimes \mathbb Z/2\mathbb Z.
\end{equation*}

{\bf\noindent Case 4 ($2\leq p<s$ and $n=2p$).}
Let $\sigma\in{\rm Aut}(\mathcal M_{s,p,n})$. By Lemmas \ref{msigma}, \ref{fibpre}, either $\sigma$ or ${\rm Dual}\circ\sigma$ is in ${\rm PGL}(s,\mathbb K)\times {\rm PGL}(n-s,\mathbb K)$. Then \begin{equation*}
    {\rm Aut}(\mathcal M_{s,p,2p})={\rm PGL}(s,\mathbb K)\times {\rm PGL}(n-s,\mathbb K)\rtimes \mathbb Z/2\mathbb Z.
\end{equation*}

{\bf\noindent Case 5a ($4=2p=2s=n$).} $\mathcal M_{2,2,4}\cong\mathbb {P}^3$, and  ${\rm Aut}(\mathcal M_{2,2,4})={\rm PGL}(4,\mathbb K)$.

\smallskip

{\bf\noindent Case 5b ($6\leq 2p=2s=n$).}
It follows by \cite[Theorem 7.5]{M3}. %Brion \cite{Br6} showed that \begin{equation*}  {\rm Aut}(\mathcal M_{p,p,2p})=\left({\rm PGL}(p,\mathbb K)\times {\rm PGL}(p,\mathbb K)\right) \rtimes\mathbb Z/2\mathbb Z\rtimes\mathbb Z/2\mathbb Z.  \end{equation*}

\smallskip

We complete the proof of Theorem \ref{mauto}.\,\,\,\,$\endpf$

\begin{appendices}
\section{Computation of intersection numbers}\label{CIN}
This appendix provides a summary of the intersection numbers between the $T$-invariant curves and the line bundles. We omit the proof due to its similarity to that of Lemma \ref{i1}. 

In what follows, $l$ is any integer with $0\leq l\leq r$.
\begin{lemma}\label{i5} 
For $2\leq j\leq r-l$, 
\[
\begin{split}
H\cdot \zeta^l_j = 0,&\,\,\,\,\,
D_i^-\cdot \zeta^l_j=-\delta_{i,l+j-1} +2\delta_{i,l+j}-\delta_{i,l+j+1},\,\,\,\,\,
D_i^+\cdot \zeta^l_j= 0.%\\
%&-K_{\mathcal T_{s,p,n}}\cdot\zeta^{0,k}_{u,v}=3-\delta_{r,j}.
\end{split}
\]
For $r-l+2\leq j\leq r$,
\[
\begin{split}
H\cdot \zeta^l_j = 0,&\,\,\,\,\,D_i^-\cdot \zeta^l_j= 0,\,\,\,\,\,
D_i^+\cdot \zeta^l_j=-\delta_{i,j-1} +2\delta_{i,j}-\delta_{i,j+1}.
\end{split}
\]
Moreover, 
\[
-K_{\mathcal T_{s,p,n}}\cdot\zeta^l_j=2+\delta_{r-l,j}+\delta_{r,j}.\]
\end{lemma}

\begin{lemma}\label{li1} 
Assume that $1\leq k\leq r-l$. When $u=l+k$ and $v=s+l+k+1\leq n$, or $v=s+l+k$ and $u=l+k+1\leq p$,
\[
\begin{split}
H\cdot \zeta^{l,k}_{u,v} = 0,&\,\,\,\,\,
D_i^-\cdot \zeta^{l,k}_{u,v}= -\delta_{i,l+k}+2\delta_{i,l+k+1}-\delta_{i,l+k+2},\,\,\,\,\,
D_i^+\cdot \zeta^{l,k}_{u,v} = 0.%\\
%&-K_{\mathcal T_{s,p,n}}\cdot\zeta^{0,k}_{u,v}=3-\delta_{r,j}.
\end{split}
\]
When $u=l+k$ and  $s+l+k+2\leq v\leq n$, 
\[
\begin{split}
H\cdot \zeta^{l,k}_{u,v} = 0,&\,\,\,\,\,
D_i^-\cdot \zeta^{l,k}_{u,v}= -\delta_{i,l+k}+\delta_{i,l+k+1}+\delta_{i,v-s}-\delta_{i,v-s+1},\,\,\,\,\,
D_i^+\cdot \zeta^{l,k}_{u,v} = 0.%\\
%&-K_{\mathcal T_{s,p,n}}\cdot\zeta^{0,k}_{u,v}=3-\delta_{r,j}.
\end{split}
\]
When $v=s+l+k$ and  $l+k+2\leq u\leq p$, 
\[
\begin{split}
H\cdot \zeta^{l,k}_{u,v} = 0,&\,\,\,\,\,
D_i^-\cdot \zeta^{l,k}_{u,v}= -\delta_{i,l+k}+\delta_{i,l+k+1}+\delta_{i,u}-\delta_{i,u+1},\,\,\,\,\,
D_i^+\cdot \zeta^{l,k}_{u,v} = 0.%\\
%&-K_{\mathcal T_{s,p,n}}\cdot\zeta^{0,k}_{u,v}=3-\delta_{r,j}.
\end{split}
\]
\end{lemma}

\begin{lemma}\label{vkli1}
Assume that $1\leq k\leq r-l$. Then, $-K_{\mathcal T_{s,p,n}}\cdot\zeta^{l,k}_{l+k,v}=$
\begin{equation*}
\left\{
\begin{array}{ll}
   2(v-s-l-k)&\,\,{\rm when}\,\, s+l+k+1\leq v\leq r+s-1\\[1em]
   2(r-l-k)+1&\,\,{\rm when}\,\, v=r+s\,\,{\rm and}\,\, k\leq r-l-1\\[1em]
   n-s+p-2(l+k)+1&\,\,{\rm when}\,\, r+s+1\leq v\leq n\,\,{\rm and}\,\, k\leq r-l-1\\[1em]
   n-s+p-2r&\,\,{\rm when}\,\, r+s+1\leq v\leq n\,\,{\rm and}\,\,k=r-l\\
\end{array}\right.
\end{equation*}
and $-K_{\mathcal T_{s,p,n}}\cdot\zeta^{l,k}_{u,s+l+k}=$
\begin{equation*}
\left\{
\begin{array}{ll}
   2(u-l-k)&\,\,{\rm when}\,\, l+k+1\leq u\leq r-1\\[1em]
   2(r-l-k)+1&\,\,{\rm when}\,\, u=r\,\,{\rm and}\,\,k \leq r-l-1\\[1em]
   n-s+p-2(l+k)+1&\,\,{\rm when}\,\, r+1\leq u\leq p\,\,{\rm and}\,\,k\leq r-l-1\\[1em]
   n-s+p-2r&\,\,{\rm when}\,\, r+1\leq u\leq p\,\,{\rm and}\,\,k=r-l\\
\end{array}\right..
\end{equation*}
\end{lemma}

\begin{lemma}\label{li2} 
Assume that $r-l+1\leq k\leq r$. When $u=r-k+1$ and $v=1\leq s-p+r-k$, or $v=s-p+r-k+1$ and $u=r-k\geq 1$, 
\[
\begin{split}
H\cdot \zeta^{l,k}_{u,v} = 0,&\,\,\,\,\,
D_i^-\cdot \zeta^{l,k}_{u,v} = 0,\,\,\,\,\,
D_i^+\cdot \zeta^{l,k}_{u,v}= -\delta_{i,k}+2\delta_{i,k+1}-\delta_{i,k+2}.%\\
%&-K_{\mathcal T_{s,p,n}}\cdot\zeta^{0,k}_{u,v}=3-\delta_{r,j}.
\end{split}
\]
When $u=r-k+1$ and  $1\leq v\leq s-p+r-k-1$, 
\[
\begin{split}
H\cdot \zeta^{l,k}_{u,v} = 0,&\,\,\,\,\,
D_i^-\cdot \zeta^{l,k}_{u,v} = 0,\,\,\,\,\,
D_i^+\cdot \zeta^{l,k}_{u,v}= -\delta_{i,k}+\delta_{i,k+1}+\delta_{i,s-p+r+1-v}-\delta_{i,s-p+r+2-v}.%\\
%&-K_{\mathcal T_{s,p,n}}\cdot\zeta^{0,k}_{u,v}=3-\delta_{r,j}.
\end{split}
\]
When $v=s-p+r-k+1$ and  $1\leq u\leq r-k-1$, 
\[
\begin{split}
H\cdot \zeta^{l,k}_{u,v} = 0,&\,\,\,\,\,
D_i^-\cdot \zeta^{l,k}_{u,v} = 0,\,\,\,\,\,
D_i^+\cdot \zeta^{l,k}_{u,v}= -\delta_{i,k}+\delta_{i,k+1}+\delta_{i,r+1-u}-\delta_{i,r+2-u}.%\\
%&-K_{\mathcal T_{s,p,n}}\cdot\zeta^{0,k}_{u,v}=3-\delta_{r,j}.
\end{split}
\]
\end{lemma}

\begin{lemma}\label{vkli2}
Assume that $r-l+1\leq k\leq r$. Then, $-K_{\mathcal T_{s,p,n}}\cdot\zeta^{l,k}_{r-k+1,v}=$
\begin{equation*}
\left\{
\begin{array}{ll}
   2(s-p+r-k+1-v)&\,\,{\rm when}\,\, s-p+2\leq v\leq s-p+r-k\,\,\\[1em]
   2(r-k)+1&\,\,{\rm when}\,\, v=s-p+1\,\,{\rm and}\,\, k\leq r-1\\[1em]
   2(r-k)+s-p+1&\,\,{\rm when}\,\, 1\leq v\leq s-p\,\,{\rm and}\,\, k\leq r-1\\[1em]
   s-p&\,\,{\rm when}\,\, 1\leq v\leq s-p\,\,{\rm and}\,\,k=r\\
\end{array}\right.
\end{equation*}
and $-K_{\mathcal T_{s,p,n}}\cdot\zeta^{l,k}_{u,s-k+r-p+1}=$
\begin{equation*}
\left\{
\begin{array}{ll}
   2(r-k+1-u)&\,\,{\rm when}\,\, 2\leq u\leq r-k\\[1em]
   2(r-k)+1&\,\,{\rm when}\,\, u=1\,\,{\rm and}\,\, k\leq r-1\\
\end{array}\right..
\end{equation*}
\end{lemma}

\begin{lemma}\label{li4}
For $1\leq m_1\leq p-l$ and $1\leq m_2\leq s-p+l$, 
\[
\begin{split}
H\cdot \delta_{m_1,m_2}^l = 1,&\,\,\,\,\,
D_i^-\cdot \delta_{m_1,m_2}^l= \delta_{i,l+m_1}-\delta_{i,l+m_1+1},\,\,\,\,\,
D_i^+\cdot \delta_{m_1,m_2}^l = \delta_{i,r-l+m_2}-\delta_{i,r-l+m_2+1}.%\\
%&-K_{\mathcal T_{s,p,n}}\cdot\zeta^{0,k}_{u,v}=3-\delta_{r,j}.
\end{split}
\]
\end{lemma}

\begin{lemma}\label{vkli4}
For $1\leq m_1\leq p-l$ and $1\leq m_2\leq s-p+l$, $-K_{\mathcal T_{s,p,n}}\cdot\delta_{m_1,m_2}^l=$
\begin{equation*}
\left\{
\begin{array}{ll}
2m_1+2m_2-2&\,\,{\rm when}\,\,1\leq m_1\leq r-l-1\,\,{\rm and}\,\,1\leq m_2\leq l-1\\[1em]
2m_1+2l-1&\,\,{\rm when}\,\,1\leq m_1\leq r-l-1\,\,{\rm and}\,\, m_2=l\\[1em]
2m_1+2l-1+s-p&\,\,{\rm when}\,\,1\leq m_1\leq r-l-1\,\,{\rm and}\,\,m_2\geq l+1\\[1em]
2(r-l+m_2)-1&\,\,{\rm when}\,\,m_1=r-l\,\,{\rm and}\,\,1\leq m_2\leq l-1\\[1em]
2r&\,\,{\rm when}\,\,m_1=r-l\,\,{\rm and}\,\, m_2=l\\[1em]
2r+s-p&\,\,{\rm when}\,\,m_1=r-l\,\,{\rm and}\,\,m_2\geq l+1\\[1em]
2(r-l+m_2)-1+s+p-n&\,\,{\rm when}\,\,m_1\geq r-l+1\,\,{\rm and}\,\,1\leq m_2\leq l-1\\[1em]
n-s+p&\,\,{\rm when}\,\,m_1\geq r-l+1\,\,{\rm and}\,\, m_2=l\\[1em]
n&\,\,{\rm when}\,\,m_1\geq r-l+1\,\,{\rm and}\,\,m_2\geq l+1\\
\end{array}\right..
\end{equation*}
\end{lemma}

\begin{lemma}\label{li5}
For $1\leq m_1\leq n-s-l$ and $1\leq m_2\leq l$, 
\[
\begin{split}
H\cdot \Delta_{m_1,m_2}^l = 1,&\,\,\,\,\,
D_i^-\cdot \Delta_{m_1,m_2}^l= \delta_{i,l+m_1}-\delta_{i,l+m_1+1},\,\,\,\,\,
D_i^+\cdot \Delta_{m_1,m_2}^l = \delta_{i,r-l+m_2}-\delta_{i,r-l+m_2+1}.%\\
%&-K_{\mathcal T_{s,p,n}}\cdot\zeta^{0,k}_{u,v}=3-\delta_{r,j}.
\end{split}
\]
\end{lemma}

\begin{lemma}\label{vkli5}
For $1\leq m_1\leq n-s-l$ and $1\leq m_2\leq l$, $-K_{\mathcal T_{s,p,n}}\cdot\Delta_{m_1,m_2}^l=$
\begin{equation*}
\left\{
\begin{array}{ll}
2m_1+2m_2-2&\,\,{\rm when}\,\,1\leq m_1\leq r-l-1\,\,{\rm and}\,\,1\leq m_2\leq l-1\\[1em]
2m_1+2l-1&\,\,{\rm when}\,\,1\leq m_1\leq r-l-1\,\,{\rm and}\,\, m_2=l\\[1em]
2(r-l+m_2)-1&\,\,{\rm when}\,\,m_1=r-l\,\,{\rm and}\,\,1\leq m_2\leq l-1\\[1em]
2r&\,\,{\rm when}\,\,m_1=r-l\,\,{\rm and}\,\, m_2=l\\[1em]
2(r-l+m_2)-1+n-s-p&\,\,{\rm when}\,\,m_1\geq r-l+1\,\,{\rm and}\,\,1\leq m_2\leq l-1\\[1em]
n-s+p&\,\,{\rm when}\,\,m_1\geq r-l+1\,\,{\rm and}\,\, m_2=l\\
\end{array}\right..
\end{equation*}
\end{lemma}

\end{appendices}


\begin{thebibliography}{00}

%% For numbered reference style
%% \bibitem{label}
%% Text of bibliographic item

%\bibitem{AB} V. Alexeev and M. Brion, {\it Boundedness of spherical Fano varieties}, The Fano Conference, 69-80, Univ. Torino, Turin, 2004.

%\bibitem[A]{A} Ahiezer, D., {\it Equivariant completions of homogeneous algebraic varieties by homogeneous divisors}, Ann. Global Anal. Geom. 1 (1983), no. 1, 49-78.

\bibitem{Al}A. Alguneid, {\it Degeneration of space collineations}, Proc. Egyptian Acad. Sci. 7 (1951), 1-17 (1952). 


\bibitem{AFKM} C. Araujo, T. Fassarella, I. Kaur, and A. Massarenti, {\it On automorphisms of moduli spaces of parabolic vector bundles}, Int. Math. Res. Not. IMRN 2021, no. 3, 2261-2283.
 
%\bibitem[B2]{B2} Bialynicki-Birula, A.,  {\it Some properties of the decompositions of algebraic varieties determined by actions of a torus}, Bull. Acad. Polon. Sci. S\'er. Sci. Math. Astronom. Phys., 1976.
%\bibitem[B3]{B3} Bialynicki-Birula, A.,  {\it On action of $SL(2)$ on complete algebraic varieties}, Pacific J. Math. 86 (1980), no. 1, 53-58.

%\bibitem[Be]{Be} Berman, R., {\it K-polystability of $Q$-Fano varieties admitting K\"ahler-Einstein metrics}, Invent. Math., 203(3): 973–1025, 2016.

%\bibitem[Be]{Be} A. Besse,  {\it Einstein manifolds}, Classics in Mathematics. Springer-Verlag, Berlin, 2008.

%\bibitem[BDP]{BDP} Bifet, E.,  De Concini, C., and  Procesi, C., {\it Cohomology of regular embeddings}, Adv. Math. 82 (1990), 1–34.

%\bibitem{Bi} A. Bia{\l}ynicki-Birula, {\it Some theorems on actions of algebraic groups}, Ann. of Math. (2) 98 (1973), 480-497.

\bibitem{BM1}  A. Bruno and M. Mella, {The automorphism group of $\overline M_{0,n}$}, J. Eur. Math. Soc. 15:3 (2013), 949-968.

\bibitem{BM2}  M. Bolognesi and A. Massarenti, {\it Birational geometry of moduli spaces of configurations of points on the line}, Algebra Number Theory 15 (2021), no. 2, 513-544.
 
%\bibitem[BK]{BK} Brion, M. and Kumar, S., {\it Frobenius splitting methods in geometry and representation theory}, Progress in Mathematics, 231. Birkhäuser Boston, Inc., Boston, MA, 2005.

%\bibitem[BLV]{BLV} Brion, M., Luna, D. and Vust, T., {\it Espaces homog\`enes sph\'eriques}, Invent. Math. 84 (1986), no. 3, 617-632.

%\bibitem[Bo1]{Bo1} Bott, R.,  {\it Nondegenerate critical manifolds}, Ann. of Math. (2) 60, (1954). 248-261. 



%\bibitem[Bo2]{Bo2} Bott, R.,  {\it A residue formula for holomorphic vector fields},  J. Diff. Geom. 1 (1967). 311-330.

\bibitem{BP} P. Bravi and G. Pezzini, {\it Primitive wonderful varieties}, Math. Z., 282 (3-4), 1067-1096 (2016)

%\bibitem{BPa}  M. Brion and F. Pauer,  {\it Valuations des espaces homog\`enes sph\'eriques}, Comment. Math. Helv. 62 (1987), no. 2, 265-285.

%\bibitem[Bor]{Bor} Borot, G.,   Lecture notes, {\it  An introductory walk in integrable woods},  \href{http://people.mpim-bonn.mpg.de/gborot/files/Hannover-14jul2015.pdf}{http://people.mpim-bonn.mpg.de/gborot/files/Hannover-14jul2015.pdf}, 2015.

%\bibitem[Br1]{Br1}  Brion, M., {\it Quelques propri\'et\'es des espaces homog\`enes sph\'eriques}, Manuscripta Math. 55 (1986), no. 2, 191–198.

%\bibitem{Br1}  M. Brion, {\it Sur l'image de l'application moment}, S\'eminaire d'alg\`ebre Paul Dubreil et Marie-Paule Malliavin (Paris, 1986), 177-192, Lecture Notes in Math., 1296, Springer, Berlin, 1987.

\bibitem{Br3}  M. Brion,  {\it Classification des espaces homogènes sph\'eriques},  Compositio Math. 63 (1987), no. 2, 189–208.

%\bibitem[Br4]{Br4}  Brion, M., {\it Spherical varieties: an introduction}, Topological methods in algebraic transformation groups (New Brunswick, NJ, 1988), 11-26, Progr. Math., 80, Birkhäuser Boston, Boston, MA, 1989.

%\bibitem{Br2} M. Brion,  {\it Groupe de Picard et nombres caract\'eristiques des vari\'et\'es sph\'eriques}, Duke Math. J. 58 (1989), no. 2, 397-424.

%\bibitem{Br3} M. Brion, {\it Vari\'et\'es sph\'eriques et th\'eorie de Mori}, Duke Math. J. 72 (1993), no. 2, 369-404.

%\bibitem{Br4}  M. Brion, {\it The behaviour at infinity of the Bruhat decomposition}, Comment. Math. Helv. 71 (1998), 137–174.

%\bibitem{Br5} M. Brion, {\it Curves and divisors in spherical varieties}, Algebraic groups and Lie groups, Austral. Math. Soc. Lect. Ser., vol. 9, Cambridge Univ. Press, Cambridge, 1997, pp. 21-34.

%\bibitem[Br9]{Br9}  Brion, M., {\it  Group completions via Hilbert schemes}, J. Algebraic Geom. 13 (2003), 603–626.

\bibitem{Br6} M. Brion, {\it The total coordinate ring of a wonderful variety}, J. Algebra 313 (2007), no. 1, 61-99. 


%\bibitem{BS} A. Bia{\l}ynicki-Birula and A. Sommese,  {\it Quotients by $\mathbb C^*$ and $SL(2,\mathbb C)$ actions}, Trans. Amer. Math. Soc. 279 (1983), no. 2, 773-800.  

%\bibitem[BS2]{BS2} Bialynicki-Birula, A., and Swiecicka, J.,{\it Complete quotients by algebraic torus actions}, Group actions and vector fields (Vancouver, B.C., 1981), 10-22, Lecture Notes in Math., 956, Springer, Berlin, 1982. 

%\bibitem[BV]{BV} W. Bruns and U. Vetter,  {\it Determinantal rings}, Lecture Notes in Mathematics, 1327. Springer-Verlag, Berlin, 1988.

%\bibitem[C]{C} Carrell, J.,  {\it Torus actions and cohomology}, Algebraic quotients. Torus actions and cohomology. The adjoint representation and the adjoint action, 83-158, Encyclopaedia Math. Sci., 131, Invariant Theory Algebr. Transform. Groups, II, Springer, Berlin, 2002. 



\bibitem{C} W. Chow,  {\it On the geometry of algebraic homogeneous spaces}, Ann. of Math. (2) 50 (1949), 32-67.

%\bibitem[Co]{Co} Coskun, I., {\it A Littlewood-Richardson rule for two-step flag varieties}, Invent. Math. 176 (2009), no. 2, 325-395.

%\bibitem[CS1]{CS1} Carrell, J. and Sommese, A.,  {\it $C^*$ actions}, Math. Scand. 43 (1978/79), no. 1, 49-59. 
%\bibitem[CS2]{CS2} Carrell, J. and Sommese, A.,  {\it Some topological aspects of $C^*$ actions on compact K\"ahler manifolds}, Comment. Math. Helv. 54 (1979), no. 4, 567-582.

%\bibitem[C]{C} W. Chow,  {\it On the geometry of algebraic homogeneous spaces}, Ann. of Math. (2) 50 (1949), 32-67.

%\bibitem[CDS1]{CDS1} Chen, X.,  Donaldson, S. and Sun, S., {\it K\"ahler-Einstein metrics on Fano manifolds. I: Approximation of metrics with cone singularities}, J. Amer. Math. Soc., 28(1): 183–197, 2015.

%\bibitem[CDS2]{CDS2} Chen, X.,  Donaldson, S. and Sun, S., {\it K\"ahler-Einstein metrics on Fano manifolds. II: Limits with cone angle less than $2\pi$}, J. Amer. Math. Soc., 28(1): 199–234, 2015. 

%\bibitem[CDS3]{CDS3} Chen, X.,  Donaldson, S. and Sun, S., {\it K\"ahler-Einstein metrics on Fano manifolds. III: Limits as cone angle approaches $2\pi$ and completion of the main proof }, J. Amer. Math. Soc., 28(1): 235–278, 2015. 
%\bibitem[Da]{Da} V. Danilov, {\it Algebraic varieties and schemes}, Algebraic geometry, I, 167–297, Encyclopaedia Math. Sci., 23, Springer, Berlin, 1994.

%\bibitem[De1]{De1} T. Delcroix, {\it K-Stability of Fano spherical varieties}, Ann. Sci. \'Ec. Norm. Sup\'er. (4) 53 (2020), no. 3, 615-662.

%\bibitem[De2]{De2} T. Delcroix, {\it K\"ahler-Einstein metrics on group compactifications}, Geom. Funct. Anal. 27 (2017), no. 1, 78-129. 

%\bibitem[De3]{De3} T. Delcroix, {\it Examples of K-unstable Fano manifolds}, Ann. Inst. Fourier (Grenoble) 72 (2022), no. 5, 2079-2108.

%\bibitem[De4]{De4} Delcroix, T., {\it Log canonical thresholds on group compactifications}, Algebr. Geom. 4 (2017), no. 2, 203–220.


%\bibitem[Dem]{Dem}  Demazure, M., {\it Automorphismes et d\'eformations des vari\'et\'es de Borel}, Invent. Math. 39 (1977), no. 2, 179–186.

%\bibitem[Do]{Do}  Donaldson, S., {\it Scalar curvature and stability of toric varieties}, J. Differential Geom. 62 (2002), no.2, 289-349.

\bibitem{DP} C. De Concini and C. Procesi,  {\it Complete symmetric varieties I}, Invariant theory (Montecatini, 1982), 1-44, Lecture Notes in Math., 996, Springer, Berlin, 1983. 

\bibitem{DP2} C. De Concini and C. Procesi, {\it Complete symmetric varieties II}, Adv. Stud. Pure Math. Vol. 6, 481–513, North-Holland, Amsterdam, 1985.

%\bibitem[FH]{FH}  Fioresi, R. and Hacon, C., {\it On infinite-dimensional Grassmannians and their quantum deformations}, Rend. Sem. Mat. Univ. Padova 111 (2004), 1-24. 

%\bibitem[FH]{FH} Fang, H. and Huang, X., {\it Geometry of grassmannians at infinity}, in preparation. 

\bibitem{Fa} G. Faltings, {\it Explicit resolution of local singularities of moduli-spaces}, J. Reine Angew. Math., 483:183-196, 1997.

%\bibitem{Fi} G. Fischer, {\it Complex analytic geometry}, Lecture Notes in Mathematics, Vol. 538. Springer-Verlag, Berlin-New York, 1976.



%\bibitem{Fan} H. Fang, {\it Canonical blow-ups of Grassmann manifolds}, arXiv: 2007.06200, 2020.

\bibitem{FM1}  B. Fantechi and A. Massarenti, {\it On the rigidity of moduli of curves in arbitrary characteristic}, Int. Math. Res. Not. IMRN 2017, no. 8, 2431-2463.

\bibitem{FM2}  B. Fantechi and A. Massarenti, {\it On the rigidity of moduli of weighted pointed stable curves}, J. Pure Appl. Algebra 222 (2018), no. 10, 3058-3074.

 
\bibitem{FSW} H. Fang, L. Schaffler, and X. Wu, {\it Fineness and smoothness of a KSBA moduli of marked cubic surfaces}, Accepted by Proceeding of American Mathematical Society.


\bibitem{Ful} W. Fulton,  {\it Intersection theory}, Springer-Verlag, Berlin, 1984.

\bibitem{FW} H. Fang and X. Wu, {\it Canonical blow-ups of Grassmannians I: how canonical is a Kausz compactification?} Int. Math. Res. Not. IMRN 2025, no. 11, Paper No. rnaf138, 25 pp.


\bibitem{FZ23} H. Fang and M. Zhang, {\it Canonical blow-ups of Grassmannians II}, arXiv: 2310.17367, 2023.

\bibitem{FZ13} B. Fu and D. Zhang, {\it A characterization of compact complex tori via automorphism groups}, Math. Ann. 357 (2013), no. 3, 961-968.
%\bibitem[FMSS]{FMSS} Fulton, W.,  MacPherson, R., Sottile, F. and Sturmfels, B., {\it Intersection theory on spherical varieties}, J. Algebraic Geom. 4 (1995), no. 1, 181-193. 


%\bibitem[Fr]{Fr} Frankel, T.,  {\it Fixed points and torsion on K\"ahler manifolds}, Ann. of Math. (2) 70 1959 1-8. 




%\bibitem[Fu]{Fu} A. Fujiki,   {\it Fixed points of the actions on compact K\"ahler manifolds}, Publ. Res. Inst. Math. Sci. 15 (1979), no. 3, 797-826.




%\bibitem[FZ]{FZ}   Fu, B. and Zhang, D., {\it A characterization of compact complex tori via automorphism groups}, Math. Ann. 357(2013), no. 3, 961-968.
 

%\bibitem[GH]{GH}  Griffiths, P. and Harris, J., {\it Principles of algebraic geometry}, Pure and Applied Mathematics, Wiley-Interscience, New York, 1978.

%\bibitem[GH]{GH} Gagliardi, G. and Hofscheier, J., {\it Gorenstein spherical Fano varieties}, Geom. Dedicata 178 (2015), 111–133.


%\bibitem[H1]{H1} Hwang, Jun-Muk,  Geometry of minimal rational curves on Fano manifolds. School on Vanishing Theorems and Effective Results in Algebraic Geometry (Trieste, 2000), 335�C393, ICTP Lect. Notes, 6, Abdus Salam Int. Cent. Theoret. Phys., Trieste, 2001. 
%\bibitem[Ha]{Ha} Hartshorne, R., {\it Algebraic geometry}, Graduate Texts in Mathematics, No. 52. Springer-Verlag, New York-Heidelberg, 1977.

%\bibitem[Hat]{Hat}  Hatcher, A., {\it Algebraic topology}, Cambridge University Press, Cambridge, 2002.

%\bibitem[He]{He} He, X., {\it The G-stable pieces of the wonderful compactification}, Trans. Amer. Math. Soc. 359 (2007), no. 7, 3005–3024.

%\bibitem[He2]{He2} He, X., {\it Unipotent variety in the group compactification}, Adv. Math. 203 (2006), no. 1, 109–131.

%\bibitem[Hi]{Hi} Hirschowitz, A., {\it Le groupe de Chow équivariant}, C. R. Acad. Sci. Paris Sér. I Math. 298 (1984), no. 5, 87-89.

%bibitem[HM]{HM} Hwang, J. and Mok, N., {\it Rigidity of irreducible Hermitian symmetric spaces of the compact type under Kähler deformation}, Invent. Math. 131 (1998), no. 2, 393–418.

%\bibitem[Ho]{Ho} Hochster, M., {\it Grassmannians and their Schubert subvarieties are arithmetically Cohen-Macaulay}, J. Algebra 25 (1973), 40-57. 

%\bibitem{Hw} J. Hwang,  {\it Nondeformability of the complex hyperquadric}, Invent. Math. 120 (1995), no. 2, 317-338. 

%\bibitem[HT]{HT} He, X. and Thomsen, J., {\it Closures of Steinberg fibers in twisted wonderful compactifications}, Transform. Groups 11 (2006), no. 3, 427–438.

%\bibitem[Kl]{Kl} Kleiman, S., {\it The transversality of a general translate}, Compositio Math. 28 (1974), 287-297. 

%\bibitem{He2} X. He, {\it Normality and {C}ohen-{M}acaulayness of local models of Shimura varieties}, Duke Math. J., 162(13):2509-2523, 2013.


\bibitem{Ka} I. Kausz,  {\it
A modular compactification of the general linear group}, Doc. Math. 5 (2000), 553-594.

\bibitem{Kn} F. Knop,  {\it The Luna-Vust theory of spherical embeddings}, Proceedings of the Hyderabad Conference on Algebraic Groups (Hyderabad, 1989), 225-249, Manoj Prakashan, Madras, 1991. 

%\bibitem[Kn2]{Kn2} Knop, F., {\it Automorphisms, root systems, and compactifications of homogeneous varieties}, J. Amer. Math. Soc., 9(1), 153–174 (1996).

%\bibitem[Ko]{Ko} S. Kov\'acs, S. {\it Rational Singularities}, \href{https://arxiv.org/pdf/1703.02269.pdf}{arXiv:1703.02269}.

%\bibitem[Ku]{Ku} Kumar, S., {\it Kac-Moody groups, their flag varieties and representation theory}, Progress in Mathematics, 204. Birkhäuser Boston, Inc., Boston, MA, 2002.

%\bibitem[La]{La} Laksov, D.,{\it The arithmetic Cohen-Macaulay character of Schubert schemes}, Acta Math. 129 (1972), no. 1-2, 1-9.


%\bibitem[Li]{Li} Liu, R., {\it An algorithmic Littlewood-Richardson rule}, J. Algebraic Combin. 31 (2010), no. 2, 253-266. 

%\bibitem[Li]{Li} Lieberman, D., {\it Compactness of the Chow scheme: applications to automorphisms and deformations of Kähler manifolds}, Fonctions de plusieurs variables complexes, III (S\'em. François Norguet, 1975–1977), pp. 140-186, Lecture Notes in Math., 670, Springer, Berlin, 1978. 

%\bibitem[LS]{LS}Lakshmibai, V. and Seshadri, C., {\it Singular locus of a Schubert variety}, Bull. Amer. Math. Soc. (N.S.) 11 (1984), no. 2, 363-366. 

%\bibitem[La1]{La1} Laksov, D., {\it Notes on the evolution of complete correlations},  Enumerative geometry and classical algebraic geometry (Nice, 1981), pp. 107–132, Progr. Math., 24, Birkhäuser Boston, Boston, MA, 1982.

%\bibitem[La2]{La2} Laksov, D., {\it Completed quadrics and linear maps}, Algebraic geometry, Bowdoin, 1985 (Brunswick, Maine, 1985), 371–387, Proc. Sympos. Pure Math., 46, Part 2, Amer. Math. Soc., Providence, RI, 1987.


\bibitem{L} D. Luna,  {\it Toute vari\'et\'e magnifique est sph\'erique}, Transform. Groups, 1(3), 249-258 (1996).

%\bibitem[Lu2]{Lu2}  Luna, D., {\it Vari\'et\'es sph\'eriques de type A}, Inst. Hautes ’Etudes Sci. Publ. Math., 94, 161-226 (2001).

\bibitem{LLT} D. Laksov, A. Lascoux, and A. Thorup, {\it On {G}iambelli’s theorem on complete correlations}. Acta Math., 162(3-4):143-199, 1989.

\bibitem{Lu2} G. Lusztig, {\it Parabolic character sheaves. {I}},  Mosc. Math. J., 4(1):153-179, 311, 2004.

%\bibitem{LV} D. Luna and T. Vust,  {\it Plongements d'espaces homogènes}, Comment. Math. Helv. 58 (1983), no. 2, 186-245. 

%\bibitem[Ma]{Ma} Manivel, L., {\it Symmetric functions, Schubert polynomials and degeneracy loci}, 3. American Mathematical Society, Providence, RI; Soci\'et\'e Mathématique de France, Paris, 2001.

%\bibitem[Mab]{Mab} T. Mabuchi,  {\it Einstein-K\"ahler forms, Futaki invariants and convex geometry on toric Fano varieties}, Osaka J. Math., 24(4): 705-737, 1987.

\bibitem{M1} A. Massarenti, {\it The automorphism group of
$\overline M_{g,n}$}, J. Lond. Math. Soc. (2) 89:1 (2014), 131-150.
 
\bibitem{M2}  A. Massarenti, {\it On the biregular geometry of the Fulton-MacPherson compactification}, Adv. Math. 322 (2017), 97-131. 
 
\bibitem{M3} A. Massarenti, {\it On the birational geometry of spaces of complete forms I: collineations and quadrics}, Proc. Lond. Math. Soc. (3) 121 (2020), no. 6, 1579-1618.

\bibitem{M4} A. Massarenti, {\it On the birational geometry of spaces of complete forms, II: Skew-forms}, J. Algebra 546 (2020), 178-200. 


%\bibitem{MM1} A. Massarenti and M. Mella, {\it On the automorphisms of moduli spaces of curves}, pp. 149-167 in Automorphisms in birational and affine geometry, edited by I. Cheltsov et al., Springer Proc. Math. Stat. 79, Springer, 2014.

\bibitem{MM1} A. Massarenti and M. Mella, {\it On the automorphisms of Hassett’s moduli spaces}, Trans. Amer. Math. Soc. 369:12 (2017), 8879-8902.

 
\bibitem{MMM} L. Manivel, M. Micha{\l}ek, L. Monin, T. Seynnaeve, and M. Vodi\v{c}ka, {\it Complete quadrics: {S}chubert calculus for {G}aussian models and semidefinite programming}, J. Eur. Math. Soc. (JEMS), 26(8):3091-3135, 2024.

%\bibitem[Mac]{Mac} Macdonald, I., {\it Symmetric functions and Hall polynomials}, Second edition. The Clarendon Press, Oxford University Press, New York, 1995.

%\bibitem[Mo]{Mo} Mok, N., Lecture notes, {\it From Rational Curves to Complex Structures on Fano Manifolds},  \href{https://hkumath.hku.hk/~nmok/HKU-Jan2005II.pdf}{https://hkumath.hku.hk/~nmok/HKU-Jan2005II.pdf}, 2005.
%\bibitem[Z]{Z} Zhuang, Xiaobo, {\it Poincar�� polynomials of moduli spaces of stable maps into flag manifolds }, 	arXiv:1601.04243 .

%\bibitem[Mul]{Mul} Mulase, M., {\it Algebraic theory of the KP equations}, Perspectives in mathematical physics, 151-217, Conf. Proc. Lecture Notes Math. Phys., III, Int. Press, Cambridge, MA, 1994.

%\bibitem[Mus]{Mus} Musili, C., {\it Some properties of Schubert varieties}, J. Indian Math. Soc. (N.S.) 38 (1974), no. 1-4, 131-145 (1975). 

%\bibitem{MT} J. Martens and M. Thaddeus, {\it Compactifications of reductive groups as moduli stacks of bundles}, Compos. Math. 152 (2016), no. 1, 62-98.

%\bibitem[Pe]{Pe} Perrin, N., {\it On the geometry of spherical varieties}, Transform. Groups 19 (2014), no. 1, 171-223. 


%\bibitem{Pez}  G. Pezzini, {\it Lectures on wonderful varieties}, Acta Math. Sin. (Engl. Ser.) 34 (2018), no. 3, 417-438.

%\bibitem[Pez2]{Pez2}  Pezzini, G., {\it Automorphisms of wonderful varieties}, Transform. Groups 14 (2009), no. 3, 677-694. 

%\bibitem[PS]{PS} Pressley, A. and Segal, G.,  {\it Loop groups}, Oxford Mathematical Monographs. Oxford Science Publications. The Clarendon Press, Oxford University Press, New York, 1986.


%\bibitem[Pu]{Pu} Putcha, M., {\it Linear algebraic monoids}, London Mathematical Society Lecture Note Series, 133. Cambridge University Press, Cambridge, 1988.

%\bibitem[Re]{Re} Renner, L., {\it Linear algebraic monoids}, Encyclopaedia of Mathematical Sciences, 134. Invariant Theory and Algebraic Transformation Groups, V. Springer-Verlag, Berlin, 2005.

%\bibitem{Ri1} A. Rittatore, {\it Algebraic monoids and group embeddings}, Transform. Groups 3 (1998), no. 4, 375-396.

%\bibitem{Ri2} A. Rittatore, {\it  Reductive embeddings are Cohen-Macaulay}, Proc. Amer. Math. Soc. 131 (2003), no. 3, 675-684.



%\bibitem[S1]{S1} M. Sato,  {\it The KP hierarchy and infinite-dimensional Grassmann manifolds}, Theta functions---Bowdoin 1987, Part 1 (Brunswick, ME, 1987), 51-66, Proc. Sympos. Pure Math., 49, Part 1, Amer. Math. Soc., Providence, RI, 1989.

%\bibitem[S2]{S2} M. Sato,  {\it Soliton equations as dynamical systems on an infinite dimensional Grassmann manifold}, RIMS Kokyuroku 439 (1981), 30-46.

%\bibitem[Sc]{Sc} Schubert, H., {\it Kalk\"ul der abz\"ahlenden Geometrie}, reprint of the 1879 original, with an introduction by Steven L. Kleiman. Springer-Verlag, Berlin-New York, 

%\bibitem{Se1} J. Semple,  {\it On complete quadrics}, J. London Math. Soc. 23 (1948), 258-267. 

%\bibitem{Se2} J. Semple, {\it The variety whose points represent complete collineations of $S_r$ on $S_{r^{\prime}}$}, Univ. Roma, Ist. Naz. Alta Mat. Rend. Mat. e Appl., 10 (5), 201-208 (1951).

%\bibitem{Sev1} F. Severi,  {\it Sui fondamenti della geometria numerativa e sulla teoria delle caratteristiche}, Atti del R. Ist. Veneto, 75 (1916), 1122-1162.

\bibitem{Sev2} F. Severi,   {\it I fondamenti della geometria numerativa}, Ann. di Mat. (4) 19 (1940),151-242.

%\bibitem{Siu} Y. Siu, {\it Nondeformability of the complex projective space}, J. Reine Angew. Math., 399:208-219, 1989.


\bibitem{Sp1} T. A. Springer, {\it  Intersection cohomology of $B\times B$-orbit closures in group compactifications}, with an appendix by Wilberd van der Kallen, J. Algebra 258 (2002), no. 1, 71-111.

%\bibitem[Sh]{Sh}  Shafarevich, I., {\it On some infinite-dimensional groups. II}, Izv. Akad. Nauk SSSR Ser. Mat. 45 (1981), no. 1, 214-226, 240. 

%\bibitem[So]{So} Sommese, A.,  {\it Some examples of $C^*$ actions}, Group actions and vector fields (Vancouver, B.C., 1981), 118-124, Lecture Notes in Math., 956, Springer, Berlin, 1982. 

%\bibitem[Sp]{Sp} Springer, T., {\it  Intersection cohomology of $B\times B$-orbit closures in group compactifications}, with an appendix by Wilberd van der Kallen, J. Algebra 258 (2002), no. 1, 71–111.

%\bibitem[Sp2]{Sp2} Springer, T., {\it Some results on compactifications of semisimple groups}, Proceedings of the International Congress of Mathematicians, vol. II, pp. 1337–1348, Eur. Math. Soc., Z\"urich, 2006.

%\bibitem[SS]{SS} Sato, M. and Sato, Y., {\it Soliton equations as dynamical systems on infinite-dimensional Grassmann manifold}, Nonlinear partial differential equations in applied science (Tokyo, 1982), 259-271, North-Holland Math. Stud., 81, Lecture Notes Numer. Appl. Anal., 5, North-Holland, Amsterdam, 1983. 

%\bibitem[Sta]{Sta} Stanley, R. {\it Enumerative combinatorics}, Vol. 1, Cambridge Studies in Advanced Mathematics, 49. Cambridge University Press, Cambridge, 1997.

%\bibitem[Str]{Str} Strickland, E., {\it A vanishing theorem for group compactifications}, Math. Ann. 277 (1987), 165–171.

\bibitem{Stu} E. Study, {\it \"Uber die Geometrie der Kegelschnitte, insbesondere deren Charakteristiken-problem}, Math. Ann., 26 (1886), 58-101.


\bibitem{Wey} J. Weyman, {\it Cohomology of vector bundles and syzygies}, Cambridge Tracts in Mathematics, 149, 2003. 

%\bibitem[Ti1]{Ti1}Tian, G., {\it K\"ahler-Einstein metrics with positive scalar curvature}, Invent. Math. 130 (1997), 239-265.

%\bibitem[Ti2]{Ti2}Tian, G., {\it K-stability and K\"ahler-Einstein metrics}, Comm. Pure Appl. Math.,68(7): 1085–1156, 2015.

%\bibitem[Ti]{Ti} Timashev, D., {\it Equivariant compactifications of reductive groups}, Mat. Sb. 194 (2003), no. 4, 119–146; translation in Sb. Math. 194 (2003), no. 3-4, 589–616.


\bibitem{Ty} J. Tyrrell,  {\it Complete quadrics and collineations in $S_n$}, Mathematika 3 (1956), 69-79. 

%\bibitem[V]{V} R. Vakil, Lecture notes, {\it  Foundations of Algebraic Geometry}, November 18, 2017 version,  \href{http://math.stanford.edu/~vakil/216blog/FOAGnov1817public.pdf}{http://math.stanford.edu/~vakil/216blog/FOAGnov1817public.pdf}, 2017.

\bibitem{Va} I. Vainsencher,  {\it Complete collineations and blowing up determinantal ideals}, Math. Ann. 267 (1984), no. 3, 417-432. 

%\bibitem[Va2]{Va2} Vainsencher, I., {\it Schubert calculus for complete quadrics},  Enumerative geometry and classical algebraic geometry (Nice, 1981), pp. 199-235, Progr. Math., 24, Birkhäuser, Boston, Mass., 1982. 

\bibitem{Van} B. van der Waerden, {\it Zur algebraischen Geometrie. XV, L\"osung des Charakteristikenproblems f\"ur Kegelschnitte}, Math. Ann. 115 (1938), no. 1, 645-655. 

%\bibitem[Vi]{Vi} Vinberg, E., {\it On reductive algebraic semigroups}, Lie groups and Lie algebras: E. B. Dynkin's Seminar, 145–182, Amer. Math. Soc. Transl. Ser. 2, 169, Adv. Math. Sci., 26, Amer. Math. Soc., Providence, RI, 1995.

%\bibitem[WZ]{WZ}  Wang, X. and Zhu, X. {\it K\"ahler-Ricci solitons on toric manifolds with positive first Chern class}, Adv. Math., 188(1): 87–103, 2004.

%bibitem[Zh]{Zh} Zhuang, Z., {\it Product theorem for K-stability}, arXiv:1904.09617. 


\end{thebibliography}
\end{document}